\documentclass[10pt]{amsart}

\usepackage{graphicx}        % standard LaTeX graphics tool
\usepackage{multicol}        % used for the two-column index
\usepackage[bottom]{footmisc}% places footnotes at page bottom
\usepackage{amscd,amsmath,amssymb,amsfonts}
\usepackage[cmtip, all]{xy}

\usepackage{bbm}
\usepackage{tikz-cd}
\usepackage{enumerate}
\usepackage{mathrsfs}
\usepackage{amsthm}

\newtheorem{theorem}{Theorem}[section]
\newtheorem{proposition}[theorem]{Proposition}
\newtheorem{lemma}[theorem]{Lemma}
\newtheorem{corollary}[theorem]{Corollary}
\newtheorem{conjecture}[theorem]{Conjecture}

\theoremstyle{definition}
\newtheorem{definition}[theorem]{Definition}
\theoremstyle{remark}
\newtheorem{remark}[theorem]{Remark}

\newtheorem{ex}[theorem]{Example}

\numberwithin{equation}{section}

\newcommand{\Pic}{{\rm Pic}}
\newcommand{\End}{{\rm End}}
\newcommand{\Hom}{{\rm Hom}}

\newcommand{\Spec}{{\rm Spec\,}}

\newcommand{\Char}{{\rm char}}
\newcommand{\Tr}{{\text{Tr}}}

\newcommand{\0}{\emptyset}
\newcommand{\sA}{{\mathcal A}}
\newcommand{\sB}{{\mathcal B}}
\newcommand{\sC}{{\mathcal C}}
\newcommand{\sD}{{\mathcal D}}
\newcommand{\sE}{{\mathcal E}}
\newcommand{\sF}{{\mathcal F}}

\newcommand{\sH}{{\mathcal H}}
\newcommand{\sI}{{\mathcal I}}

\newcommand{\sK}{{\mathcal K}}
\newcommand{\sL}{{\mathcal L}}

\newcommand{\sN}{{\mathcal N}}
\newcommand{\sO}{{\mathcal O}}

\newcommand{\sT}{{\mathcal T}}
\newcommand{\sU}{{\mathcal U}}
\newcommand{\sV}{{\mathcal V}}

\newcommand{\sX}{{\mathcal X}}

\newcommand{\A}{{\mathbb A}}

\newcommand{\C}{{\mathbb C}}

\newcommand{\G}{{\mathbb G}}
\renewcommand{\H}{{\mathbb H}}

\renewcommand{\P}{{\mathbb P}}

\newcommand{\R}{{\mathbb R}}

\newcommand{\Z}{{\mathbb Z}}

\newcommand{\Cplx}{\operatorname{Cplx}}
\newcommand{\lcm}{\operatorname{lcm}}

\renewcommand{\det}{\operatorname{det}}

\newcommand{\id}{{\operatorname{\rm Id}}}

\newcommand{\Sch}{{\operatorname{\mathbf{Sch}}}}

\newcommand{\<}{\langle}
\renewcommand{\>}{\rangle}

\newcommand{\tensor}{\otimes}
\newcommand{\ext}{\bigwedge}

\renewcommand{\dim}{{\operatorname{\rm dim}}}

\newcommand{\del}{\partial}

\newcommand{\Sm}{{\mathbf{Sm}}}

\newcommand{\Proj}{{\operatorname{Proj}}}

\newcommand{\Ext}{{\operatorname{Ext}}}

\newcommand{\rnk}{{\operatorname{\text{rnk}}}}

\newcommand{\GW}{{\operatorname{GW}}} 
 
\newcommand{\SH}{{\operatorname{SH}}} 
\newcommand{\Th}{{\operatorname{Th}}}

\newcommand{\Kos}{{\operatorname{Kos}}} 
\newcommand{\Aut}{{\operatorname{Aut}}}

\newcommand{\topo}{{\operatorname{top}}}

\newcommand{\GL}{\operatorname{GL}}

\newcommand{\dg}{e}
\newcommand{\Hess}{\operatorname{Hess}}
\newcommand{\Der}{{\operatorname{Der}}}

\newcommand{\can}{{can}}

\newcommand{\ev}{\text{\it ev}}

\newcommand{\ind}[1]{}

\newcommand{\inp}[1]{}

\newcommand{\res}{{\operatorname{res}}}
\newcommand{\Log}{{\operatorname{log}}}

\newcommand{\coker}{{\operatorname{coker}}}

\newcommand{\Hodge}{{\operatorname{Hodge}}}
\newcommand{\Jac}{{\operatorname{Jac}}}
\newcommand{\Eul}{{\operatorname{Eul}}}

\newcommand{\spc}{{\operatorname{sp}}}
\newcommand{\Cone}{{\operatorname{Cone}}}

\newcommand{\one}{\mathbbm{1}}
\newcommand{\SP}{\operatorname{SP}}
\newcommand{\spe}{\operatorname{sp}}
\newcommand{\sgn}{\operatorname{sgn}}

\makeatletter
 \@namedef{subjclassname@2020}{\textup{2020} Mathematics Subject Classification}
\makeatother

\begin{document}

\title[Euler characteristics of hypersurfaces]{Euler characteristics of homogeneous and weighted-homogeneous hypersurfaces}

\date{ \today}

\author[M.~Levine]{Marc~Levine}
\address{Universit\"at Duisburg-Essen,
Fakult\"at Mathematik, Campus Essen, 45117 Essen, Germany}
\email{marc.levine@uni-due.de}

\author[S.~Pepin Lehalleur]{Simon Pepin Lehalleur}
\address{Radboud University, Heyendaalseweg 135, 6525 AJ Nijmegen Netherlands}
\email{simon.pepin.lehalleur@gmail.com}

\author[V.~Srinivas]{Vasudevan Srinivas}
\address{TIFR, Mumbai, India}
\email{srinivas@math.tifr.res.in}

%\keywords{}

\setcounter{tocdepth}{1}

\begin{abstract} Let $k$ be a perfect field and let $\GW(k)$ be the Grothendieck-Witt  ring  of (virtual) non-degenerate symmetric bilinear forms over $k$. We develop methods for computing the quadratic Euler characteristic $\chi(X/k)\in \GW(k)$ for $X$ a smooth hypersurface in a projective space and in a weighted projective space. We raise the question of a quadratic refinement of classical conductor formulas and find such a formula for the degeneration of a smooth hypersurface $X$ in $\P^{n+1}$ to the cone over a smooth hyperplane section of $X$; we also find a similar formula in the weighted homogeneous case. We formulate a conjecture that generalizes these computations to similar types of degenerations. Finally, we give an interpretation of the quadratic conductor formulas in terms of Ayoub's nearby cycles functor. 
\end{abstract}
\thanks{S.P.L. is supported by the Netherlands Organisation for Scientific Research (NWO), under project number 613.001.752. \\ V.S. is supported by a J. C. Bose Fellowship of the Department
of Science and Technology, India. V.S. also recognizes
support of the Department of Atomic Energy, Government of India, under
project number 12-R\&D-TFR-RTI4001.   \\
M.L. gratefully acknowledges support from the DFG through the SPP 1786 and from the ERC through the project QUADAG.  This paper is part of a project that has received funding from the European Research Council (ERC) under the European Union's Horizon 2020 research and innovation programme (grant agreement No. 832833).\\ 
\includegraphics[scale=0.08]{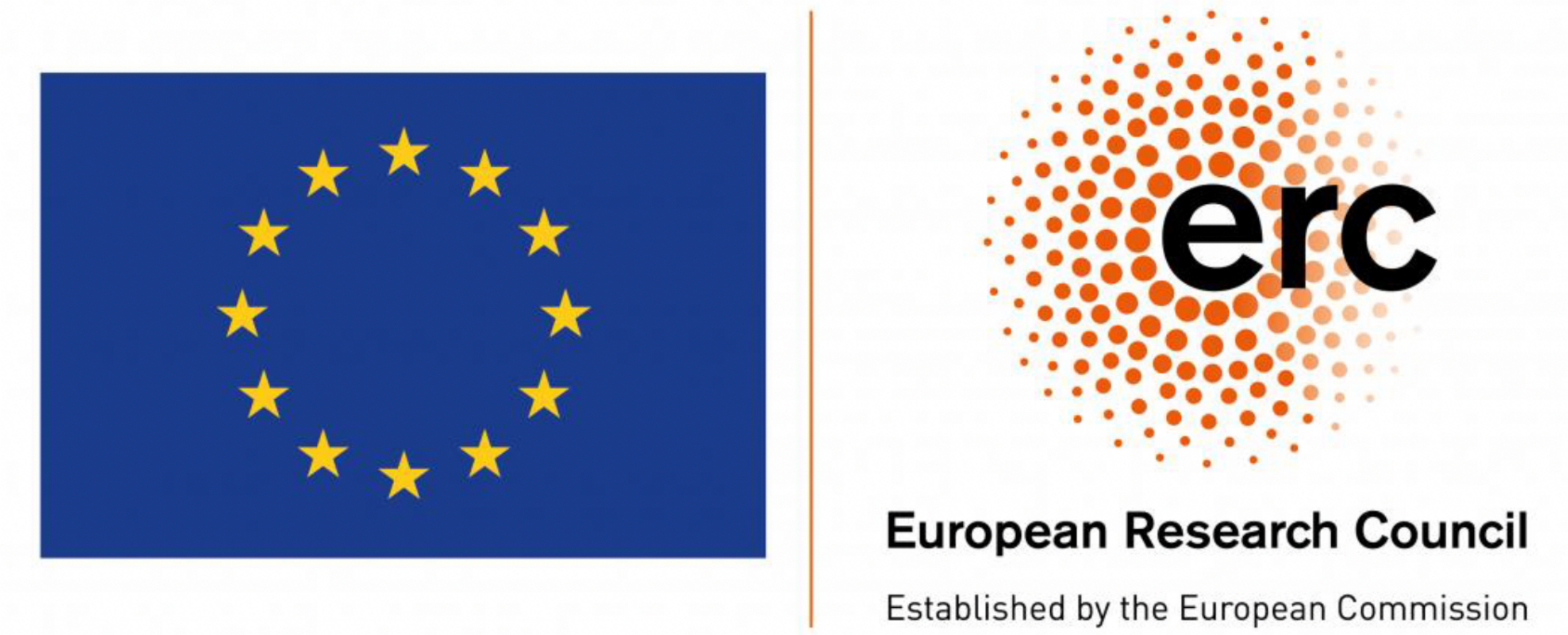}}

\maketitle

\tableofcontents

\section{Introduction}  Let $k$ be a perfect field. For a smooth projective $k$-scheme $X$, the corresponding suspension spectrum $\Sigma^\infty_{\P^1}(X_+)$ in the motivic stable homotopy category $\SH(k)$ is strongly dualizable (see \cite[Theorem 5.22]{Hoyois6}, \cite[Appendix A]{Hu},  \cite{Riou}  and   \cite[\S2]{Voev})  hence has a well-defined Euler characteristic $\chi_{\SH(k)}(\Sigma^\infty_{\P^1}(X_+))$ in the endomorphism ring of the unit object $\one_k\in \SH(k)$. Morel \cite[Lemma 6.3.8, Theorem 6.4.1]{MorelICTP} has defined an isomorphism of rings $\GW(k)\xrightarrow{\sim} \End_{\SH(k)}(\one_k)$  of the Grothendieck-Witt ring of virtual non-degenerate symmetric bilinear forms to this endomorphism ring, giving the ``quadratic'' Euler characteristic 
\[
\chi(X/k)\in \GW(k). 
\]

The quadratic Euler characteristic $\chi(X/k)$ for a smooth projective $k$-scheme is a refinement to the Grothendieck-Witt ring of the classical Euler characteristic, defined using \'etale cohomology, de Rham cohomology or singular cohomology, as appropriate. Carrying  additional information, $\chi(X/k)$  is naturally more difficult to compute than its classical counterpart. With A. Raksit \cite{LR}, the first named author has given a description of $\chi(X/k)$ as the bilinear form on Hodge cohomology $\oplus_{p,q}H^q(X, \Omega_{X/k}^p)$ given by cup product followed by the trace map $\Tr_{X/k}:H^n(X, \Omega^n_{X/k})\to k$, with $n=\dim_kX$. Although this gives a quite concrete expression for $\chi(X/k)$, explicit computations are still not easy; it is the purpose of this note to give other descriptions more amenable to computation.

As is well known, the primitive Hodge cohomology of a smooth hypersurface $X\subset \P^{n+1}$ is computable as suitable graded pieces of the Jacobian ring $J(F):=k[X_0,\ldots, X_{n+1}]/(\ldots, \del F/\del X_i,\ldots)$ of the defining equation $F$, at least if $k$ has characteristic zero. This goes back to results of Griffiths \cite{Griffiths} and was pursued further in Carlson-Griffiths \cite{CG}, where they showed a certain compatibility of this identification with respect to product structures. We will review these results here. This material has been treated elsewhere, for instance, Steenbrink \cite{Steenbrink} handled the case of hypersurfaces in a weighted projective space over a characteristic zero field, without an explicit discussion of products, and the additive theory in the generality discussed here may be found  in the article of Dolgachev \cite{Dolgachev}.  

Combined with the results of \cite{LR} computing the quadratic Euler characteristic in terms of Hodge cohomology, relating the multiplicative structure on the Jacobian ring with that of the Hodge cohomology gives a quite explicit description of the quadratic Euler characteristic of a smooth hypersurface in a projective space, as well as for  smooth hypersurfaces in a weighted projective space $\P(a_0,\ldots, a_{n+1})$ assuming that the weighted degree $e$ and all the $a_i$ are prime to the characteristic and that $e$ is divisible by the lcm of the $a_i$.

 See Theorem~\ref{thm:Main} and Corollary~\ref{cor:HdgJacQ} for our results for a smooth hypersurface in a projective space and Theorem~\ref{thm:HdgJacQH} for  smooth hypersurface in a weighted projective space.

A motivating problem underlying these computations is the search for a quadratic replacement of the conductor formulas of Milnor, Deligne, Bloch, Saito, Kato-Saito, and others. One considers a flat proper morphism $f:\sX\to \Spec \sO$ with $\sO$ a dvr having parameter $t$, with $\sX$ a regular scheme and with generic fiber $X_K$ smooth over the quotient field $K$ of $\sO$.  A ``conductor formula'' is an expression for the difference $\chi^{top}_c(X_K)-\chi^{top}_c(X_0)$ in terms of algebraic invariants of $f$; here $\chi^{top}_c(-)$ is the Euler characteristic of compactly supported cohomology. For a morphism of relative dimension zero, this is local ramification theory, where in positive and mixed characteristic, the Swan conductor makes its appearance.

Suppose first that $X_{0}$ has isolated singularities,  the case of interest for this paper.  In characteristic zero, the conductor formula in this case  goes back to Milnor \cite{Milnor}, where the answer is given by $\pm$ the sum of the Milnor numbers (dimensions of the Jacobian rings at the critical points of $f$), this being the same as the sum of local indices at the critical  points of $f$ of the section of the cotangent bundle of $\sX$ given by $f^*(dt)$. In positive characteristic, one needs to add information about wild ramification, contained in the Swan conductor; this case has been treated by Deligne \cite[Expos\'e XVI, Th\'eor\`eme 2.4]{DK}. The resulting ``Deligne-Milnor formula'' is still open in mixed characteristic. Note that in this situation of isolated singularities, one can use vanishing cycle functors to express  $\chi^{top}_c(X_K)-\chi^{top}_c(X_0)$ in a purely local way around the critical points of $f$.

The mixed characteristic case for a relative curve was treated by Bloch \cite{Bloch}, who also formulated a general conjecture, now known as the ``Bloch conductor formula'', without assuming that $X_{0}$ has isolated singularities (the precise relation with the Deligne-Milnor formula is worked out by Orgogozo \cite{Orgogozo}). The Bloch conductor formula has now been established in equal characteristic zero \cite{Kapranov} and positive equal characteristic \cite{Saito}. The case of mixed characteristic is still open in general; it has been settled when the reduced special fiber has normal crossings singularities \cite{KatoSaito} and there is a recent approach by To\"en-Vezzosi \cite{TV} using $K$-theory and dg-categories and assuming only that the inertia acts unipotently.

One can ask for an analogue of the Deligne-Milnor formula for the Euler characteristics of the real points of a map of $\R$-schemes. This has been treated by Eisenbud-Levine \cite{EL} and Khimshiasvili \cite{Kh}, where they use the signature of a quadratic form on the Jacobian ring defined by Scheja-Storch \cite{SchejaStorch}. The Scheja-Storch form has been recognized as giving the local quadratic Euler characteristic of the section of $\Omega_{\sX/k}$ defined by $df$ over an arbitrary field (of characteristic $\neq2$) by Kass-Wickelgren \cite{KassWickelgren} and Bachmann-Wickelgren \cite{BachmannWickelgren}.

One important point that distinguishes the quadratic case from the topological one is that the invariants $\chi_c(X_K/K)$ and $\chi_c(X_0/k)$ live in different rings, namely $\GW(K)$ and $\GW(k)$. To compare them, one needs to apply a specialization map to $\chi_c(X_K/K)$, with the added complication that specialization depends on the choice of parameter. Fortunately, the Scheja-Storch form also has this dependence, so one can hope that the two are still compatible. Another distinction is that given  a globally defined morphism $\pi:\bar\sX\to \P^1$, there does not appear to be any nice relation of the difference of the Euler characteristics of the fibers near the critical values and the global Euler characteristic $\chi(\bar\sX/k)$.  Although one does have a quadratic version of the Riemann-Hurwitz formula in this global case (see \cite[Corollary 10.4]{LevineEuler}), the lack of a local-global relation means that the method of Deligne for proving the local formula via the 
 Riemann-Hurwitz formula  and an explicit computation in the case of ordinary double points does not have an obvious extension to  the quadratic case.

Another important difference with the classical situation is that, in situations of wild ramification in positive characteristic, we do not know what the quadratic refinement of the Swan conductor should be.

We use our computation of the quadratic Euler characteristic of hypersurfaces to look at the simplest case of a quadratic conductor formula, namely the case of a degeneration of a smooth hypersurface to the cone over a hyperplane section, and an analog of this in the case of a hypersurface in a weighted projective space. This gives examples of conductor formulas for those homogeneous and quasi-homogeneous singularities that correspond to smooth hypersurfaces in a (weighted) projective space. 

Perhaps surprisingly, we find that the direct generalization of the Deligne-Milnor formula (replacing topological Euler characteristics with quadratic ones and Milnor numbers with Scheja-Storch forms) does not hold even in this simple case. We find that it is necessary to replace occurrences of 1's and -1's in the classical formulas with certain rank 1 and virtual rank -1 quadratic forms, with coefficients depending on the hypersurface degree (in the weighted homogeneous case, the product of the weights also comes in). See Theorem~\ref{thm:ConductorFormula1} and
Theorem~\ref{thm:ConductorFormula2} for precise statements. We do not have a conceptual interpretation of these ``correction terms''.

M.L. would like to thank Daniel Huybrechts for very helpful discussions on the multiplicative properties detailed in the work of Griffiths and Carlson-Griffiths.

\section{Bott's theorem and the residue sequence}  
We fix a field $k$; unless explicitly stated otherwise, all schemes are separated and of finite type over $k$. For $Y$ a $k$-scheme, we will write $\Omega^1_Y$ for the sheaf of K\"ahler differentials $\Omega_{Y/k}$ and  $\Omega^m_Y$ for the sheaf $\bigwedge^m_{\sO_Y}\Omega_{Y/k}$.

\subsection{Bott's theorem}
We recall Bott's theorem on the cohomology of $ \Omega^q_{\P^n}(m)$.

\begin{theorem}[Bott's theorem for projective space]  The cohomology of $\Omega^q_{\P^n}(m)$ satisfies
\[
H^p(\P^n, \Omega^q_{\P^n}(m))=\begin{cases} 0&\text{ for }p>0\text{ and }m\ge q-n, m\neq0,\\0& \text{ for }p>0, m=0\text{ and }p\neq q\\
0&\text{ for }p=0,   m\le q, \text{ except for }m=p=q=0.
\end{cases}
\]
Moreover $H^q(\P^n, \Omega^q_{\P^n})\cong k$ for $0\le q\le n$ and for $m>q\ge0$, 
\[
\dim_kH^0(\P^n, \Omega^q_{\P^n}(m))=\begin{pmatrix}m+n-q\\m\end{pmatrix}\cdot\begin{pmatrix}m-1\\q\end{pmatrix}.
\]
\end{theorem}

\begin{proof} See \cite[Theorem (2.3.2)]{Dolgachev}, where the case of a weighted projective space is also treated.  
\end{proof}

\subsection{The residue sequence} Since we will need to keep careful track of signs in our computations, we make explicit here our conventions. 

\begin{remark}[Cones, cone sequences and the sign in the coboundary map] \label{rem:Sign} We fix the sign of the coboundary map arising from a termwise exact sequence of complexes.  For simplicity, we work in an abelian category $\sA$ that admits elements. For $?\in\{b,+, -,\0\}$ we let $\Cplx^?(\sA), K^?(\sA), D^?(\sA)$ denote the category of complexes, the homotopy category and the derived category of $\sA$, with the corresponding boundedness condition, respectively. We use throughout cohomological complexes, so the differential has degree +1.

For a complex $(A^*, d_A)$, we have the shift $(A[1], d_{A[1]})$ with $A[1]^n=A^{n+1}$ and $d_{A[1]}^n=-d_A^{n+1}$; given a morphism of complexes $f^*:A^*\to B^*$, the shift $f[1]^*:A[1]^*\to B[1]^*$ is defined by $f[1]^n=f^{n+1}$. We identify $H^n(A[1]^*, d_{A[1]})$ with 
$H^{n+1}(A^*, d_A)$ via the identity map on $A^{n+1}$.

 Given a map of complexes $f^*:(A^*, d_A)\to (B^*, d_B)$, we have the cone $\Cone(f)$ with
\[
\Cone(f)^n=A^{n+1}\oplus B^n
\]
and with differential
\[
d^n_\Cone(a,b)=(-d_A(a), d_B(b)+f^{n+1}(a))
\]
Let $i_B^*:B^*\to \Cone(f)^*$ and $p_{A[1]}^*:\Cone(f)^*\to A[1]^*$ be the inclusion of $B^*$ and projection to $A[1]^*$, both maps of complexes.

The cone sequence associated to $f$ is
\[
A^*\xrightarrow{f^*} B^*\xrightarrow{i_B^*}\Cone(f)^*\xrightarrow{p_{A[1]}^*}A[1]
\]
and these are the sequences that give rise to the triangulated structure in the homotopy category $K^?(\sA)$ and derived category $D^?(\sA)$, $?=+,-, b,\0$.

Given a termwise exact sequence of complexes
\begin{equation}\label{eqn:ExactSeqComplex}
0\to C^*_1\xrightarrow{f} C^*_2\xrightarrow{g}  C^*_3\to 0
\end{equation}
we have two natural quasi-isomorphisms
\[
\alpha^*:C^*_1\to \Cone(g)^*[-1],\ \beta^*: \Cone(f)^*\to C_3^*
\]
with $\alpha^n(c_1)=(f^n(c_1),0)$, $\beta^n(c_1, c_2)=g^n(c_2)$. Note that the sequence
\[
\Cone(g)^*[-1]\xrightarrow{p^*_{C_2^*}}C^*_2\xrightarrow{g}C^*_3\xrightarrow{-i^*_{C_3}}
\Cone(g)
\]
is isomorphic in $K(\sA)$  to the cone sequence for $p^*_{C_2^*}$, hence defines a distinguished triangle. 

We have the commutative diagrams in $C(\sA)$
\[
\xymatrix{
C^*_1\ar[r]^{f}\ar@{=}[d]& C^*_2\ar[r]^{g}\ar@{=}[d] & C^*_3\\
C^*_1\ar[r]^{f}& C^*_2\ar[r]^-{i_{C_2}^*}& \Cone(f)^*\ar[r]^{p_{C_1[1]}}\ar[u]^{\beta^*}&C_1[1]^*
}
\]
and
\[
\xymatrix{
C^*_1\ar[r]^{f}\ar[d]^{\alpha^*}& C^*_2\ar[r]^{g}\ar@{=}[d] & C^*_3\ar@{=}[d]\\
\Cone(g)[-1]^*\ar[r]^-{p_{C_2}}& C^*_2\ar[r]^{g}& C^*_3\ar[r]^-{-i_{C_3}}&\Cone(g)^*
}
\]
Via the quasi-isomorphisms $\beta^*$ and $\alpha[1]^*$, these give us the maps
\[
H^n(C_3^*)\xrightarrow{H^n(p_{C_1[1]})\circ H^n(\beta)^{-1}}H^n(C_1[1]^*)=
H^{n+1}(C_1)
\]
and
\[
H^n(C_3^*)\xrightarrow{H^n(\alpha[1]^*)^{-1}\circ H^n(-i_{C_3})}
H^n(C_1[1]^*)=H^{n+1}(C_1^*)
\]
An direct computation shows that these are both define the same map
$\delta^n:H^n(C^*_3)\to H^{n+1}(C^*_1)$, given as follows.

Take an element $z\in \ker(d^n_{C_3})$ representing a cohomology class $[z]\in H^n(C_3^*)$. Lift $z$ to an element $y\in C^n_2$, $g^n(y)=z$. Then $d_{C_2}^n(y)=f^{n+1}(x)$ for a unique $x\in C_1^{n+1}$ and $d_{C_1}(x)=0$. Then $\delta^n([z])=-[x]$ (note the minus sign!).
\end{remark}

\begin{remark}[Derived functors and $\delta$-functors]\label{rem:DerivedFunctorDeltaFunctor} Let $\sA$ and $\sB$ be abelian categories with enough injectives and let $F:\sA\to \sB$ be a left-exact functor, giving the (right) derived functors $R^nF:D^+(\sA)\to \sB$, defined as
\[
R^nF(A^*)=H^n(F(I^*))
\]
where $A^*\to I^*$ is quasi-isomorphism of $A^*$ with a complex of injectives $I^*$. 

Given a short exact sequence $0\to a\to b\to c\to 0$ in $\sA$, one forms a termwise short exact sequence of injective resolutions
\[
\xymatrix{
0\ar[r]&I^*_a\ar[r]&I^*b\ar[r]&I^*_c\ar[r]&0\\
0\ar[r]&a\ar[r]\ar[u]&b\ar[r]\ar[u]&c\ar[r]\ar[u]&0.
}
\]
Applying $F$ gives the termwise exact sequence of complexes in $\sB$
\[
\xymatrix{
0\ar[r]&F(I^*_a)\ar[r]&F(I^*b)\ar[r]&F(I^*_c)\ar[r]&0 
}
\]
and the associated boundary maps $\delta^n:R^nF(c)\to R^{n+1}F(a)$ defined as in Remark~\ref{rem:Sign} define a universal $\delta$-functor $\{R^nF, \delta^n\}_{n\ge0}$.

Take an object $x\in \sA$ concentrated in degree zero, and let
\[
0\to x\xrightarrow{\epsilon} I^0\to\ldots\to I^n\to \ldots
\]
be an injective resolution of $x$. We can break up the resolution into short exact sequences
\[
0\to K^j\xrightarrow{i_j} I^j\xrightarrow{p_{j+1}} K^{j+1}\to 0
\]
with $K^0\to I^0$ the augmentation $x\to I^0$. Let $F:\sA\to \sB$ be a left-exact functor, as above. Then the resolution $x\to I^*$ and the map $K^n\to I^n$ define the canonical surjective map
\[
R^0F(K^n)=F(K^n)\to   H^n(F(I^*))=R^nF(x).
\]
with kernel $R^nF(p_n)(R^nF(I^{n-1}))$. 

We also have the sequence of boundary maps given by the short exact sequences above, 
\[
R^0F(K^n)\xrightarrow{\delta^0}R^1F(K^{n-1})\xrightarrow{\delta^1}\ldots
\xrightarrow{\delta^{n-1}}R^nF(K^0)=R^nF(x).
\]
The fact that we have normalized the boundary maps by using the triangulated structure in $D^+(\sB)$ implies these differ by the factor $(-1)^{n(n-1)/2}$.

More generally, given an arbitrary exact sequence
\[
0\to x\xrightarrow{\epsilon} N^0\to\ldots\to N^n\to 0
\]
an injective resolution $\alpha: N^*\to I^*$ defines an injective resolution $\alpha\circ\epsilon:x\to I^*$, so we have the isomorphism 
the quasi-isomorphism $\epsilon$ induces an isomorphism
\[
RF^n(\epsilon)^*:RF^n(N^*)\to RF^n(x).
\]
Breaking up the above sequence into short exact sequences as before gives the sequence of boundary maps
\[
RF^0(N^n)\xrightarrow{\delta^0}\ldots\xrightarrow{\delta^{n-1}} RF^n(x).
\]
The identity map on $N^n$ gives the map $N^n[-n]\to N^*$ and $(-1)^{n(n-1)/2}\cdot \delta^{n-1}\circ\ldots\circ\delta^0$ agrees with the map 
\[
RF^0(N^n)= RF^n(N^n[-n])\to RF^n(N^*)\xrightarrow{RF^n(\epsilon)^*}RF^n(x).
\]
\end{remark}

Let $X\subset \P^n$ be a smooth hypersurface defined by $F\in k[X_0,\ldots, X_n]$, homogeneous of degree $\dg$, with inclusion $i:X\to \P^n$. We have the sheaf of differential with log poles along $X$, $\Omega_{\P^n}^p(\Log X)$ and the residue sequence
\begin{equation}\label{eqn:ResSeq}
0\to \Omega_{\P^n}\to \Omega^p_{\P^n}(\Log X)\xrightarrow{\res_X}i_*\Omega_X^p\to0
\end{equation}
This gives the coboundary map $\delta^{p,q}_X:H^q(X, \Omega^p_X)\to H^{q+1}(\P^n,\Omega^{p+1}_{\P^n})$, following Remark~\ref{rem:Sign}
 for our sign convention on coboundary maps.

For $L$ a line bundle on some smooth $k$-scheme $Y$, we have the 1st Chern class $c_1(L)\in H^1(Y,\Omega_Y^1)$, defined as the image of $[L]\in H^1(Y, \sO_Y^\times)$ via the dlog map $dlog:\sO_Y^\times\to \Omega^1_Y$.

In \cite{SriHodge}, the third named author has constructed functorial push-forward maps in Hodge cohomology,
 $f_*:H^q(X, \Omega^pX_)\to H^{q+d}(Y, \Omega^{p+d}_Y)$ for each projective morphism $f:X\to Y$ of relative codimension $d$, where $X$ and $Y$ are in $\Sm_k$. Combined with the usual pullback maps, these satisfy the usual properties of a Bloch-Ogus theory on $\Sm_k$, such as a projection formula and a base-change identity in transverse cartesian squares. Moreover, for $p:X\to \Spec k$ smooth and proper of dimension $d$ over $k$, the canonical trace map $\Tr_{X/k}:H^d(X,\Omega^d_X)\to k$ is equal to the push-forward map $p_*$.

\begin{lemma}\label{lem:ResiduePushforward} Let $i:X\to \P^n$ be a hypersurface. Then $\delta^{p,q}_X=i_*:H^q(X, \Omega^p_X)\to H^{q+1}(\P^n,\Omega^{p+1}_{\P^n})$.
\end{lemma}

\begin{proof} We first recall the construction of $i_*$ from \cite{SriHodge} in the case of a regular closed immersion $i:X\to Y$ of codimension one. We have the standard exact sequence
\[
0\to \sO_Y(-X)\to \sO_Y\to i_*\sO_X\to 0.
\]
For $\sE$ a locally free  $\sO_X$-module and $\sF$ a locally free $\sO_Y$-module, this sequence induces an isomorphism
\[
\alpha_{\sE, \sF}:\Hom_{\sO_X}(\sE, i^*\sF(X))\xrightarrow{\sim} \Ext^1_{\sO_Y}(i_*\sE, \sF)
\]
as follows: given a map $\phi:\sE\to  i^*\sF(X)$ of $\sO_X$-modules, we get the map of $\sO_Y$-modules $\phi':i_*\sE\to i_*i^*\sF(X)$. Tensoring the standard exact sequence above with $\sF(X)$ gives the exact sequence
\[
0\to \sF\to \sF(X)\to  i_*i^*\sF(X)\to 0
\]
and pulling back by $\phi'$ gives the exact sequence 
\[
0\to \sF\to \sF_\phi\to i_*\sE\to 0.
\]

Now suppose $X\subset \P^n$ is defined by some homogeneous polynomial $F$.  Take $\sE=\Omega^p_X$, $\sF=\Omega^{p+1}_Y$ and $\phi:\Omega^p_X\to i^*\Omega^{p+1}_{\P^n}(X)$ the map $dF/F\wedge -$. We have the commutative diagram with exact rows
\[
\xymatrix{
0\ar[r]&\Omega^{p+1}_{\P^n}\ar@{=}[d]\ar[r]&\Omega^{p+1}_{\P^n}(\Log X)\ar@{_(->}[d]\ar[r]&
i_*\Omega^p_X\ar[r]\ar[d]^{dF/F\wedge-}&0\\
0\ar[r]&\Omega^{p+1}_{\P^n}\ar[r]&\Omega^{p+1}_{\P^n}(X)\ar[r]&
i_*i^*\Omega^{p+1}_{\P^n}(X)\ar[r]&0
}
\]
which shows that the class of the residue sequence in $\Ext^1_{\sO_{\P^n}}(i_*\Omega^p_X, \Omega^{p+1}_{\P^n})$ is equal to $\alpha_{\Omega^p_X, \Omega^{p+1}_{\P^n}}(dF/F\wedge-)$. Since 
$dF/F\cdot -$ is precisely the map $\Omega^p_X\to i^*\Omega^{p+1}_{\P^n}(X)$ induced from the short exact sequence
\[
0\to \sO_X(-X)\xrightarrow{dF/F\cdot-} i^*\Omega^1_{\P^n}\xrightarrow{\pi}\Omega^1_X\to 0
\]
we see that $\alpha_{\Omega^p_X, \Omega^{p+1}_{\P^n}}(dF/F\wedge-)$ is exactly the Ext-element denoted $\alpha_p(i)$ used in \cite[pg. 3]{SriHodge}  to construct the map $i_*$.

Given a sheaf $\sF$ on a $k$-scheme $Y$, the identity $H^0(Y, \sF)=\Hom_{\sO_Y}(\sO_Y, \sF)$ extends canonically to   the isomorphism
\[
\beta_{q,\sF}:\Ext^q_{\sO_Y}(\sO_Y, \sF)\xrightarrow{\sim}
H^q(Y, \sF)
\]
of derived functors.  $i_*$ is the map given  by taking the product with $\alpha_p(i)$:
\begin{multline*}
H^q(X, \Omega^p_X)\xrightarrow{\beta_q^{-1}}\Ext^q_{\sO_X}(\sO_X, \Omega^p_X)\xrightarrow{\alpha_p(i)\circ-}
\Ext^{q+1}_{\sO_{\P^n}}(\Omega_{\P^n}^{p+1}, i_*\sO_X)\\\xrightarrow{\pi^*}
\Ext^{q+1}_{\sO_Y}(\Omega_{\P^n}^{p+1}, \sO_{\P^n})\xrightarrow{\beta^{q+1}}H^{q+1}({\P^n}, \Omega_{\P^n}^{p+1})
\end{multline*}
The fact that $\alpha_p(i)=\alpha_{\Omega^p_X, \Omega^{p+1}_{\P^n}}(dF/F\wedge-)$, which in turn is the class of the residue sequence, shows that $i_*=\delta^{p,q}_X$.
\end{proof}

\section{The Jacobian ring and primitive Hodge cohomology}

Let $X\subset \P^{n+1}$ be a hypersurface defined by $F=0$, $F\in k[X_0,\ldots, X_{n+1}]$ of degree $\dg>1$, with closed immersion $i:X\to \P^{n+1}$. Let $F_i:=\del F/\del X_i$. The Jacobian ring $J(F)$ is defined as
\[
J(F):=k[X_0,\ldots, X_{n+1}]/(F, F_0,\ldots, F_{n+1})
\]
with grading induced by the usual grading by degree on $k[X_0,\ldots, X_{n+1}]$. In this section, we describe the isomorphism of a suitable graded subspace of $J(F)$ with the primitive Hodge cohomology of $X$, for $X$ smooth and of degree prime to the characteristic. As mentioned in the introduction, this result is already in the literature, but we include a discussion here for the reader's convenience and to present constructions necessary for our applications. 

Suppose that $X$ is smooth and $\dg$ is prime to the characteristic. The Euler equation
\[
\dg\cdot F=\sum_{i=0}^{n+1}F_i\cdot X_i
\]
shows that $F$ is in $(F_0,\ldots, F_{n+1})$. The ideal
$(F,  F_0,\ldots, F_n)$ is a $(X_0,\ldots, X_{n+1})$-primary ideal since $X$ is smooth,  and  thus $J(F)$ is a complete intersection. One has the following facts:
\begin{enumerate}
\item The non-zero graded summand of $J(F)$ of top degree has dimension one over $k$ and is in degree $(\dg-2)(n+2)$.\\
\item There is a canonical choice of generator $e_F$ for $J(F)_{(\dg-2)(n+2)}$, we call this the {\em Scheja-Storch generator}.
\end{enumerate}

The generator $e_F$ is defined as follows: We have $(F_0,\ldots, F_{n+1})\subset (X_0,\ldots, X_{n+1})$ since  $\dg>1$,  so there are (non-unique!) elements $a_{ij}\in k[X_0,\ldots, X_{n+1}]$ with
\[
F_i=\sum_j a_{ij}X_j;\quad i=0,\ldots, n+1
\]
Let $\tilde{e}_F:=\det(a_{ij})$ and let $e_F\in J(F)_{(\dg-2)(n+2)}$ be the image of $\tilde{e}_F$. One shows (see for instance \cite[Lemma 1.2($\alpha$)]{SchejaStorch}) that $e_F$ is independent of the choice of the $a_{ij}$. One can also show that $(\dim_{k} J(F))\cdot e_{F}$ is the Hessian determinant $\det\left(\frac{\partial^{2}f}{\partial X_{i}\partial X_{j}}\right)$ \cite[Korollar 4.7]{SchejaStorch}.

Let $N=(\dg-2)(n+2)$. We have the symmetric bilinear form 
\[
B_\Jac:J(F)\otimes J(F)\to k
\]
defined on $J(F)_m\otimes J(F)_{N-m}$ by $B_\Jac(x,y)=\lambda\Leftrightarrow x\cdot y=\lambda\cdot e_F$; $B_\Jac(x,y)=0$ if $\deg x+\deg y\neq N$. This bilinear form is non-degenerate.

We let $\sI_X$ denote the ideal sheaf of $X$; the choice of  degree $\dg$ defining equation $F$ gives an isomorphism $\sI_X\cong \sO_{\P^{n+1}}(-\dg)$. For $q\ge0$,  let $\sO_{\P^{n+1}}(qX)$ denote the invertible sheaf of rational functions $g$ on $\P^{n+1}$, regular on $\P^{n+1}\setminus X$ and with pole order at most $q$ along $X$. For a coherent sheaf $\sE$ on $\P^{n+1}$ and $q\ge0$, we write $\sE(qX)$ for $\sE\otimes_{\sO_{\P^{n+1}}}\sO_{\P^n}(qX)$, and for $q<0$, set $\sE(qX)=\sE\otimes_{\sO_{\P^{n+1}}}\sI_X^{\otimes -q}$. Again, the choice of defining equation $F$ gives a canonical isomorphism $\sE(qX)\cong \sE(q\dg)$ for all $q\in \Z$.

We note that the logarithmic differential $dF/F$  gives a well-defined section of $i^*\Omega_{\P^{n+1}}^1(X)$ and that we have  the exact sequence
\begin{equation}\label{eqn:CanDifSeq}
0\to  \sO_X(-X)\xrightarrow{dF/F\cdot-}i^*\Omega_{\P^{n+1}}^1\xrightarrow{\pi} \Omega^1_X\to 0
\end{equation}

\begin{lemma}\label{lem:ExactSeq}  Suppose that $\dg$ is prime to the characteristic of $k$. \\[5pt]
1. We have an exact sequence of sheaves on $X$
\begin{multline*}
0\to \Omega^p_X\xrightarrow{dF/F\wedge-} i^*\Omega_{\P^{n+1}}^{p+1}(X)\xrightarrow{dF/F\wedge-} i^*\Omega_{\P^{n+1}}^{p+2}(2X)\xrightarrow{dF/F\wedge-}\\\ldots\xrightarrow{dF/F\wedge-}i^*\Omega_{\P^{n+1}}^n(qX)\xrightarrow{\pi_q} \Omega_X^n(qX)\to 0
\end{multline*}
for  integers $p,q\ge0$ with $p+q=n$. \\[2pt]
2.  Let $\sC(p)^*$ denote the complex
\begin{multline*}
0\to  i^*\Omega_{\P^{n+1}}^{p+1}(X)\xrightarrow{dF/F\wedge-} i^*\Omega_{\P^{n+1}}^{p+2}(2X)\xrightarrow{dF/F\wedge-}\\\ldots\xrightarrow{dF/F\wedge-}i^*\Omega_{\P^{n+1}}^n(qX)\xrightarrow{\pi_q} \Omega_X^n(qX)\to 0
\end{multline*}
with $ i^*\Omega_{\P^{n+1}}^{p+1}(X)$ in degree 0. 
The map $\delta:H^0(X, \Omega_X^n(qX))\to \H^q(X, \sC(p)^*)\cong H^q(X, \Omega^p_X)$ induced by the exact sequence (1) gives an exact sequence
\[
H^0(X, i^*\Omega_{\P^{n+1}}^n(qX))\xrightarrow{\pi_q} H^0(X, \Omega_X^n(qX))\xrightarrow{\delta} H^q(X, \Omega^p_X)
\]
for all $q\ge0$, $p+q=n$. Moreover, $\delta$ is surjective if $p\neq q$ and has image 
$H^p(X, \Omega_X^p)_{prim}$ in case $p=q$, $n=2p$. 
\end{lemma}

\begin{proof} (1) follows  by patching together the exact sequences
\begin{equation}\label{eqn:SSS}
0\to \Omega_{X}^{n-j-1}((q-j-1)X)\xrightarrow{dF/F\wedge-} i^*\Omega_{\P^{n+1}}^{n-j}((q-j)X)\to \Omega_X^{n-j}((q-j)X)\to 0.
\end{equation}
derived from \eqref{eqn:CanDifSeq}
for $j=0, \ldots, q-1$.

For (2), we have the spectral sequence
\[
E_1^{a,b}=H^b(X, \sC(p )^a)\Rightarrow \H^{a+b}(X, \sC(p )^*)\cong H^{a+b}(X, \Omega^p_X).
\]
The exact sequence
\begin{equation}\label{eqn:ses1}
0\to \Omega_{\P^{n+1}}^{n-j}((q-j-1)X) \to  \Omega_{\P^{n+1}}^{n-j}((q-j)X) \to i_*i^*\Omega_{\P^{n+1}}^{n-j}((q-j)X)\to0
\end{equation}
together with Bott's theorem  show that 
\[
H^j(X, i^*\Omega_{\P^{n+1}}^{n-j}((q-j)X))=0
\]
for  $j\ge1$, $q-j\ge 2$, and also for $j\ge1$,  $q-j=1$, except possibly for the case $n=2j+1=2q-1$, $n-j=p+1$. 

In this case,  we are considering the exact sequence
\[
0\to \Omega_\P^{p+1} \to  \Omega_\P^{p+1}(X) \to i_*i^*\Omega_{\P}^{p+1}(X)\to0
\]
whose  long exact cohomology sequence gives
\[
H^{q-1}(X, i^*\Omega_{\P}^{p+1}(X))\cong
H^q(\P, \Omega_{\P}^{p+1})
\]
which is non-zero ($\cong k$) exactly if $q=p+1$, forcing $n=2p+1$. 

Feeding this into the exact sequences \eqref{eqn:SSS}, we see that 
$H^1(X, \Omega_X^n(qX))\cong H^p(X,\Omega^{p+1}_X(X))$ and we have the exact sequence
\begin{equation}\label{eqn:ses2}
H^p(X, \Omega^p_X)\to H^p(X, i^*\Omega_\P^{p+1}(X))\xrightarrow{\alpha} H^p(X,\Omega^{p+1}_X(X))\to
H^{p+1}(X,\Omega^p_X)\to 0
\end{equation}

Thus $E_1^{a,b}=0$ except for the cases
\[
E_1^{a,b}=\begin{cases} 
H^{p+1}(\P, \Omega^{p+1}_\P)&\text{ for } a=0, b=p\\
H^0(X, i^*\Omega^{p+1-a}_\P((a+1)X)&\text{ for }b=0, 0\le a\le q-1\\
H^b(X, \sO_X(qX))&\text{ for }a=q \end{cases}
\]
and $E_1^{q,1}\cong H^p(X,\Omega^{p+1}_X(X))$.
Thus, the only non-zero differentials are $d_1^{a,0}:E_1^{a,0}\to E_1^{a+1,0}$ and, in case $p=q$, the differential
$d_p^{0,p}=\pm\alpha:H^{p+1}(\P, \Omega^{p+1}_\P)\to H^1(X,\Omega^n(qX))\cong H^p(X,\Omega^{p+1}_X(X))$.

We are only interested in computing $H^q(X,\Omega^p_X)$ for $p+q=n$, and for this case the spectral sequence  shows that, for $p\neq q$, the edge homomorphism $H^0(X, \sO_X(qX))=E_1^{q,0}\to \H^q(X, \sC(p )^*)=H^q(X, \Omega^p_X)$ gives the exact sequence
\[
H^0(X, i^*\sO_\P(qX))\xrightarrow{\pi_q}H^0(X, \sO_X(qX))\to H^q(X, \Omega^p_X)\to0
\]
For $p=q=n/2$,  we have the additional edge homomorphism $\H^p(X, \sC(p )^*)=H^p(X, \Omega^p_X)\to E_1^{0,p}\cong H^{p+1}(\P, \Omega^{p+1}_\P)$, giving the exact sequence
\begin{multline*}
H^0(X, i^*\Omega^n_\P(pX))\xrightarrow{\pi_q}H^0(X, \Omega^n_X(pX))\to H^p(X, \Omega^p_X)\\\xrightarrow{\beta}
H^{p+1}(\P, \Omega^{p+1}_\P) \xrightarrow{\alpha}H^p(X,\Omega^{p+1}_X(X))
\end{multline*}
with the map $\beta$ the composition of $H^p(dF/F\wedge-):H^p(X, \Omega^p_X)\to H^p(X, i^*\Omega^{p+1}_\P(X))$ with the isomorphism $H^p(X, i^*\Omega_{\P}^{p+1}(X))\cong
H^{p+1}(\P, \Omega_{\P}^{p+1})$ we have constructed above. 

We claim that the map $\beta$ is (up to sign) the map $i_*:H^p(X, \Omega^p_X)\to H^{p+1}(\P, \Omega^{p+1}_\P)$ given (by Lemma~\ref{lem:ResiduePushforward}) by the coboundary in the residue sequence
\[
0\to \Omega^{p+1}_\P\to \Omega^{p+1}_\P(\Log X)\xrightarrow{\res_X}i_*\Omega^p_X\to 0
\]
This follows from the commutative diagram
\[
\xymatrix{
0\ar[r]&\Omega_\P^{p+1}\ar@{=}[d]\ar[r] &\Omega_\P^{p+1}(\Log X) \ar[r]^{res} \ar@{_(->}[d]&i_*\Omega_X^{p}\ar[d]^{dF/F\wedge -}\ar[r]&0\\
0\ar[r]&\Omega_\P^{p+1} \ar[r]&  \Omega_\P^{p+1}(X) \ar[r]& i_*i^*\Omega_{\P}^{p+1}(X)\ar[r]&0
}
\]

Now suppose that the degree $\dg$ is prime to the characteristic of $k$.  The projection formula gives
\[
i_*((c_1(\sO_X(1)))^p)=\dg\cdot c_1(\sO_\P(1)))^{p+1},
\]
hence $i_*((c_1(\sO_X(1)))^p)\neq0$, as $c_1(\sO_\P(1)))^{p+1}$ is a generator of the $k$-vector space $H^{p+1}(\P, \Omega^{p+1})$. Using  the projection formula again, we see that $\ker\beta=H^p(X,\Omega^p_X)_{prim}$.   This proves (2).

\end{proof}

Let $\Omega$ be the $(n+1)$-form
\begin{equation}\label{eqn:Omega}
\Omega:=\sum_{i=0}^{n+1}(-1)^iX_idX_0\wedge\ldots\wedge\widehat{dX}_i\wedge\ldots\wedge dX_{n+1}
\end{equation}
on $\A^{n+2}$. Since $\Omega$ is killed by interior multiplication with the Euler vector field $\sum_iX_i\del/\del X_i$, we may view $\Omega$ as a generating section of $\Omega^{n+1}_{\P^{n+1}}(n+2)\cong \sO_{\P^{n+1}}$. Thus, for each integer $q\ge0$, sending $A\in k[X_0,\ldots, X_{n+1}]_{\dg(q+1)-n-2}$ to $A\Omega/F^{q+1}$ gives an isomorphism
\[
\hat\psi_q: k[X_0,\ldots, X_{n+1}]_{\dg(q+1)-n-2}\xrightarrow{\sim} H^0(\P^{n+1}, \Omega^{n+1}_{\P^{n+1}/k}((q+1)X)),
\]

\begin{proposition} \label{prop: Main} Suppose that $\dg$ is prime to the characteristic of $k$. Then for each $q\ge0$, the map $\delta\circ \res_X\circ\hat\psi_q$ descends to an isomorphism 
\[
\psi_q:J(F)_{(q+1)\dg-n-2}\to H^q(X, \Omega^{n-q}_X)_{prim}
\]
 \end{proposition}
 
\begin{remark} This result, as well as its extension to the case of a hypersurface in a weighted projective space (Proposition~\ref{prop:JacobianQuotientHypersurf}), has been proven by Dolgachev \cite[\S 4.2]{Dolgachev}; the characteristic zero case goes back to Griffiths \cite{Griffiths}   and in the weighted case to Steenbrink \cite{Steenbrink}. We include the proof here and also the extension to the weighted case to give the proper background for our discussion of the multiplicative structure.
\end{remark}
 
\begin{proof} 
 
First assume $n\ge1$. For all $q\ge0$, Bott's theorem gives us the exact sequence
\[
0\to H^0(\P^{n+1}, \Omega^{n+1}_{\P^{n+1}}(qX))\to H^0(\P^{n+1}, \Omega^{n+1}_{\P^{n+1}}((q+1)X))\xrightarrow{i^*} H^0(X, i^*\Omega^{n+1}_{\P^{n+1}}((q+1)X))\to0
\]
so $\res_X$ descends to the map
\[
\bar\res_X:H^0(X, i^*\Omega^{n+1}_{\P^{n+1}}((q+1)X))\to H^0(X, \Omega_X^n(qX))
\]
with the same image as $\res_X$. This also shows that the isomorphism $\hat\psi_q$ descends to an isomorphism
\[
\tilde\psi_q:\left(k[X_0,\ldots, X_{n+1}]/(F)\right)_{\dg(q+1)-n-2}\to H^0(X, i^*\Omega^{n+1}_{\P^{n+1}}((q+1)X)).
\]

From the description of the residue map, we have the commutative diagram
\[
\xymatrix{
 H^0(X, i^*\Omega^n_{\P^{n+1}}(qX))\ar[r]^-{dF/F\wedge-}\ar[dr]_{\pi_q}& H^0(X, i^*\Omega^{n+1}_{\P^{n+1}}((q+1)X))\ar[d]^{\bar\res_X}\\
&H^0(X, \Omega^n_{X}(qX))
}
\]

Let $R=k[X_0,\ldots, X_{n+1}]$. We have the graded piece of the Koszul complex 
\[
\ldots\to R_{m-i}\otimes \bigwedge^iR_1\to R_{m-i+1}\otimes \bigwedge^{i-1}R_1\to\ldots,
\]
which gives us the resolution of $H^0(\P^{n+1},\Omega^n_{\P^{n+1}}(m))$ as
\[
0\to R_{m-n-2}\to R_{m-n-1}\otimes \bigwedge^{n+1} R_1\xrightarrow{\phi} H^0(\P^{n+1},\Omega^n_{\P^{n+1}}(m))\to 0
\]
Letting $\ell$ be an integer, $0\le \ell\le n+1$, one checks that for $m=n+1$, $\phi(X_0\wedge\ldots\wedge\widehat{X}_\ell\wedge\ldots\wedge X_{n+1})$ is the $n$-form 
\[
\Omega^{(\ell)}:=\sum_{i=0}^{\ell-1}(-1)^{i}X_idX^{\widehat{i,\ell}}-\sum_{i=\ell}^{n+1}(-1)^{i}X_idX^{\widehat{\ell,i}},
\]
where for $0\le a<b\le n+1$, 
\[
dX^{\widehat{a,b}}:=dX_0\wedge\ldots\wedge\widehat{dX}_a\wedge\ldots\wedge\widehat{dX}_b\wedge\ldots\wedge dX_{n+1}.
\]
Similarly, taking $m=\dg\cdot q$, we see that the map 
\[
(k[X_0,\ldots, X_{n+1}]_{\dg\cdot q-n-1})^{n+2}\to 
H^0(\P^{n+1}, \Omega^n_{\P^{n+1}}(qX))
\]
sending $(A_0,\ldots, A_{n+1})$ to  $\sum_{i=0}^{n+1}A_\ell\Omega^{(\ell)}/F^q$  is surjective. Moreover, a direct computation shows that (as forms on $\A^{n+2}$)
\[
dF\wedge \Omega^{(\ell)}=(-1)^{\ell-1}F_\ell\cdot \Omega+ \dg\cdot F\cdot dX^{\widehat{\ell}}
\]
Together with the surjectivity of $H^0(\P^{n+1}, \Omega^n_{\P^{n+1}}(qX))\to H^0(X, i^*\Omega^n_{\P^{n+1}}(qX))$, this shows that the image of 
\[
dF/F\wedge-:H^0(X, i^*\Omega^n_{\P^{n+1}}(qX))\to H^0(X, i^*\Omega^{n+1}_{\P^{n+1}}((q+1)X))
\]
 is   the image under $\hat\psi_q$ of the degree $\dg(q+1)-n-2$ part of  $(F, F_0,\ldots, F_{n+1})/(F)$. As the image of $\pi_q$ is equal to the kernel of $\delta:H^0(X, \Omega^n_{X}(qX))\to H^q(X, \Omega^p_X)$, this completes the proof.  
 
For the case $n=0$,  we have the exact residue sequence
 \[
 0\to H^0(\P^1, \Omega^1_{\P^1}(X))\xrightarrow{\res_X}H^0(X,\sO_X) \xrightarrow{i_*}
 H^1(\P^1, \Omega^1_{\P^1})\to 0
 \]
with $H^0(X,\sO_X)_{prim}:=\ker i_*$ by definition. This replaces the exact sequence at the beginning of the argument for $n\ge1$; after this the proof is the same, except much simpler.
\end{proof}

Let  $\Tr_{X/k}:H^n(X, \Omega^n_X)\to k$ be the canonical trace map. We have the perfect pairing
\[
B_{\Hodge}:=\Tr_{X/k}\circ\cup:\oplus_{p+q=n}H^q(X, \Omega^p_X)\otimes_k \oplus_{p+q=n}H^q(X, \Omega^p_X)\to k
\]
defined on $H^q(X, \Omega^p_X)\otimes_kH^{n-q}(X, \Omega^{n-p}_X)$ as 
\[
B_\Hodge(\alpha,\beta):=\Tr_{X/k}(\alpha\cup\beta)
\]
and is zero on the other pairs of summands.

In a number of articles on the comparison of the Jacobian ring with the primitive cohomology, it is mentioned that  $B_\Jac$ and $B_\Hodge$ agree up to a non-zero rational number. This is proven  by Carlson-Griffiths \cite{CG} assuming $k$ has characteristic zero.  Here we give a proof in arbitrary characteristic prime to $\dg$, using the isomorphism given above rather  than the isomorphism defined by Carlson-Griffiths, and determine the factor with respect to this isomorphism.

 Carlson-Griffiths express the cohomology groups as \v{C}ech cohomology to compute the product; we follow their method. In an Appendix, we give another approach.

To fix notation, let $\sV=V_0,\ldots, V_m$ be an open cover of some scheme $Y$. For an index $J=(j_0,\ldots, j_q)$, $0\le j_i\le m$,  we let $V_J=\cap_{i=0}^q V_{j_i}$. We use the convention that for a sheaf $\sF$ on $Y$, $\sC^q(\sV, \sF)$ is the set of families $\{s_J\in \sF(V_J)\}_J$ with $J=(j_0,\ldots, j_q)$ and $0\le j_0<\ldots <j_q\le m$. We define $|J|=q$ for $J=(j_0,\ldots, j_q)$, and for $0\le i\le q$ let $J^{\hat{i}}=(j_0,\ldots,j_{i-1}, j_{i+1}, \ldots, j_q)$. We use  the usual \v{C}ech coboundary map
\[
\delta:\sC^q(\sV, \sF)\to \sC^{q+1}(\sV, \sF),\ \delta(\{s_I\in \sF(V_I)\})_J=
\sum_{i=0}^{q+1}(-1)^is_{J^{\hat{i}}}|_{V_J}
\]
to define the \v{C}ech complex $(\sC^*(\sV, \sF), \delta)$.

Returning to the situation at hand, let  $\sV=\{V_0,\ldots, V_{n+1}\}$ be the open cover of $X$, with $V_j$  the complement of $F_j=0$ in $X$.

Over $V_j$, the exact sequence \eqref{eqn:CanDifSeq} splits, with splitting $H_j$ of $dF/F\wedge-:\sO_X(-X)\to i^*\Omega^1_{\P^{n+1}}$ defined by
\[
H_j(\sum_{l=0}^{n+1}a_l dX_l)=(F/F_j)\cdot a_j
\]
This extends to the map
\[
H_j:i^*\Omega^a_{\P^{n+1}}(bX)_{|V_j}\to i^*\Omega^{a-1}_{\P^{n+1}}((b-1)X)_{|V_j}
\]
Let $K_j:\Omega_{\A^{n+2}}^a\to \Omega_{\A^{n+2}}^{a-1}$ be (sign-graded) interior multiplication with $\del/\del X_j$,
\[
K_j(g dX_{i_1}\wedge\ldots\wedge dX_{i_a})=\begin{cases} 0&\text{ if }j\not\in \{i_1,\ldots, i_a\}\\
(-1)^{b-1}g dX_{i_1}\wedge\ldots\wedge dX_{i_{b-1}}\wedge dX_{i_{b+1}}\wedge\ldots\wedge dX_{i_a}\kern-80pt\\
&\text{ if } i_b=j.
\end{cases}
\]
We can express $H_j$ as
\[
H_j:=(F/F_j)\cdot K_j
\]
We define the map $H$ on \v{C}ech cochains 
\[
H:\sC^q(\sV, i^*\Omega^a_{\P^{n+1}}(bX))\to \sC^q(\sV, i^*\Omega^{a-1}_{\P^{n+1}}((b-1)X))
\]
by 
\[
H(\{\omega_{j_0,\ldots, j_q}\})=\{H_{j_0}(\omega_{j_0,\ldots, j_q})\}.
\]

Recall our differential form $\Omega\in H^0(\A^{n+2}, \Omega^{n+1}_{\A^{n+2}})$,
\[
\Omega:=\sum_{i=0}^{n+1} (-1)^iX_i dX_0\cdots \widehat{dX_i}\cdots dX_{n+1}.
\]
For $J=(j_0,\ldots, j_q)$, let $K_J=K_{j_q}\circ\ldots\circ K_{j_0}$, $\Omega_J=K_J(\Omega)$ and $F_J=F_{j_0}\cdots F_{j_q}$. Let $\pi:i^*\Omega^a_{\P^{n+1}}(bX)\to \Omega^a_{X}(bX)$ be the canonical surjection.

For $A\in k[X_0,\ldots, X_{n+1}]$ of degree $\dg\cdot(b+1)-n-2$, the differential form $A\Omega/F^{b+1}$ gives a well-defined global section of $\Omega^{n+1}_{\P^{n+1}}((b+1)X)$. 

\begin{lemma}\label{lem:Main} Take $A\in k[X_0,\ldots, X_{n+1}]$ of degree $\dg\cdot(b+1)-n-2$. Then:\\[5pt]
1. $\res_X(A\Omega/F^{b+1})\in H^0(X, \Omega_X^n(bX))$ is represented by 
\[
\{\pi(A\Omega_i)/F_iF^b\}_i\in \sC^0(\sV, \Omega_X^n(bX)).
\]
2. For $\{(A\Omega_I/F_IF^b)\}_I\in \sC^q(\sV,i^*\Omega^{a-1}_{\P^{n+1}}((b-1)X))$, 
\[
H(\{(dF/F)\wedge (A\Omega_I/F_IF^b)\}_I)=\{A\Omega_I/F_IF^b\}_I
\]
3. Let $\delta$ be the \v{C}ech differential. Then for $\{A\Omega_I/F_IF^b\}_I\in \sC^q(\sV,\Omega^{a}_{\P^{n+1}}(bX))$,  
\[
\delta(\{A\Omega_I/F_IF^b\}_I)=
\{(-1)^{q+1}(dF/F)\wedge( A\Omega_J/F_JF^{b-1}\}_J\in \sC^{q+1}(\sV,\Omega^{a-1}_{\P^{n+1}}((b-1)X)).
\]
4. $\pi\circ H:\sC^q(\sV, i^*\Omega^a_{\P^{n+1}}(bX))\to \sC^q(\sV,\Omega^{a-1}_{X}((b-1)X))$ is a splitting to the injective map $dF/F\wedge-: \sC^q(\sV,\Omega^{a-1}_{X}((b-1)X))\to \sC^q(\sV,i^*\Omega^a_{\P^{n+1}}(bX))$.
\end{lemma}

\begin{proof}[Proof, taken from the proof of \hbox{\cite[Lemma, pg. 10]{CG}}]  For all these identities, we may work in differential forms on suitable open subsets of  $\A^{n+2}$ or of the affine cone $\tilde{X}$ over $X$.

For (1),  let $E=\sum_{i=0}^{n+1}X_i\del/\del X_i$ be the Euler vector field and  $dV=dX_0 \ldots dX_n$ the volume form. Then  interior multiplication with $E$ satisfies
$\iota(E)(dV)=\Omega$ and thus
\[
0=\iota(E)(dF\wedge dV)=\sum_{i=0}^{n+1}X_iF_i\wedge dV-dF\wedge\Omega,
\]
so the Euler equation yields the identity $dF\wedge\Omega=\dg\cdot F\cdot dV$.  Applying $K_i$ gives 
\[
F_i\cdot\Omega-dF\wedge\Omega_i=\dg\cdot F\cdot K_idV
\]
or
\[
A\Omega/F^{b+1}-(dF/F)\wedge (A\Omega_i/F_iF^b)=\dg\cdot F\cdot A\cdot K_idV/(F^{b+1}F_i)
\]
Restricting to the affine cone over $X$, $\tilde{i}:\tilde{X}\to \A^{n+2}$, gives
\[
\tilde{i}^*(A\Omega/F^{b+1})=\tilde{i}^*((dF/F)\wedge (A\Omega_i/F_iF^b))
\]
as sections of $\tilde{i}^*(\Omega^{n+1}_{\A^{n+2}}((b+1)\tilde{X})$ over the inverse image of $V_i$. This proves (1).

For (2),  $K_idF=F_i$ and for $I=(i_0,\ldots, i_q)$,  $K_{i_0}\Omega_I=K_{i_0}K_I(\Omega)=0$.
Thus
\[
H_{i_0}(dF/F\wedge(A\Omega_I/F_IF^b))= (F/F_{i_0})\cdot (F_{i_0}/F)\cdot (A\Omega_I/F_IF^b)=A\Omega_I/F_IF^b.
\]

For (3), we have
\begin{align*}
\delta(\{A\Omega_I/F_IF^b\}_I)_J&=\sum_{l=0}^{q+1}(-1)^l
A\Omega_{J^{\widehat{j_l}}} /F_{J^{\widehat{j_l}}}F^b\\
&=\sum_{l=0}^{q+1}(-1)^l
A\Omega_{J^{\widehat{j_l}}}F_{j_l} /F_JF^b
\end{align*}
Applying $K_J$ to the identity $dF\wedge\Omega=\dg\cdot F\cdot dV$ gives as above
\[
(-1)^{q+1}\tilde{i}^*(dF\wedge\Omega_J)=\tilde{i}^*(\sum_{l=0}^{q+1}(-1)^lF_{j_l}\cdot \Omega_{J^{\widehat{j_l}}})
\]
so 
\[
\tilde{i}^*(\sum_{l=0}^{q+1}(-1)^lA\cdot \Omega_{J^{\widehat{j_l}}}/F_{J^{\widehat{j_l}}}F^b)
=\tilde{i}^*(-1)^{q+1}(dF/F)\wedge (A\Omega_J/F_JF^{b-1}).
\]
Thus
\[
\delta(\{A\Omega_I/F_IF^b\}_I)_J=(-1)^{q+1}(dF/F)\wedge (A\Omega_J/F_JF^{b-1}),
\]
which proves (3).

(4) follows from  the fact that $H_j$ gives a splitting to $\sO_X(-X)\to i^*\Omega^1_{\P^{n+1}}$ over $V_j$.
\end{proof}

\begin{proposition} \label{prop:Rep} For $\bar{A}\in J(F)_{(q+1)\dg-n-2}$, represented by a homogeneous form $A\in H^0(\P^{n+1}, \sO_{\P^{n+1}}((q+1)\dg-n-2))$, the cohomology class $\psi_q(\bar{A})\in H^q(X, \Omega^{n-q}_X)$ is represented by the cochain $\{ \pi(A\Omega_I/F_I)\}\in \sC^q(\sV, \Omega^{n-q}_X)$.
\end{proposition}

\begin{proof} The short exact sequence \eqref{eqn:SSS} gives us the coboundary map
\[
\delta_j:H^j(X, \Omega_X^{n-j}((q-j)X))\to H^{j+1}(X, \Omega_{X}^{n-j-1}((q-j-1)X))
\]
We first prove by induction on $j$ that $\delta_{j-1}\circ\ldots\circ \delta_0(\res_X(A\Omega/F^{q+1}))$ is  represented by the a \v{C}ech cocycle $\{(-1)^{j(j-1)/2}\pi(A\Omega_I/F_IF^{q-j})\}\in \sC^j(\sV, \Omega_X^{n-j}((q-j)X))$. 

The case $j=0$ is Lemma~\ref{lem:Main}(1). Supposing the induction hypothesis holds for $j$, we follow our sign convention and represent $\delta_j\circ \delta_{j-1}\circ\ldots \circ\delta_0(\res_X(A\Omega/F^{q+1}))$ by  the  coboundary of the cocycle $\{(-1)^{j(j-1)/2}\pi(A\Omega_I/F_IF^{q-j})\}\in \sC^j(\sV, \Omega_X^{n-j}((q-j)X))$. Following Remark~\ref{rem:Sign}, this is defined by lifting to the cochain 
\[
\{(-1)^{j(j-1)/2}A\Omega_I/F_IF^{q-j})\}\in \sC^j(\sV, i^*\Omega_{\P^{n+1}}^{n-j}((q-j)X))
\]
and applying the negative of the \v{C}ech coboundary operator 
\[
\delta:\sC^j(\sV, \Omega_{\P^{n+1}}^{n-j}((q-j)X))\to \sC^{j+1}(\sV, i^*\Omega_{\P^{n+1}}^{n-j}((q-j)X)), 
\]
which has image a cocycle in 
$ \sC^{j+1}(\sV, \Omega_{X}^{n-j-1}((q-j-1)X))$. By Lemma~\ref{lem:Main}(2,3,4), we have
\begin{multline*}
(-1)\delta(\{(-1)^{j(j-1)/2}A\Omega_I/F_IF^{q-j}\}_{|I|=j})\\=\{(-1)^{j}(-1)^{j(j-1)/2}\pi(A\Omega_I/F_IF^{q-j-1})\}_{|I|=j+1}
\end{multline*}
and since $j(j-1)/2+j=(j+1)j/2$,  the induction goes through. 

By Remark~\ref{rem:Sign}, the cohomology class $\psi_q(\bar{A})\in H^q(X, \Omega^{n-q}_X)$ is represented by $(-1)^{q(q-1)/2}\cdot \delta_{q-1}\circ\ldots\circ \delta_0(\res_X(A\Omega/F^{q+1}))$, that is, $\psi_q(\bar{A})\in H^q(X, \Omega^{n-q}_X)$ is represented by $\{ \pi(A\Omega_I/F_I)\}\in \sC^q(\sV, \Omega^{n-q}_X)$, as desired.
\end{proof}

\begin{proposition}\label{prop:Product}  Take non-negative integers $p, q$ with $p+q=n$, and take  $A\in J(F)_{(q+1)\dg-n-2}$, $B\in J(F)_{(p+1)\dg-n-2}$. Let $\omega_A= \psi_q(A)\in H^q(X, \Omega^p_X)$, $\omega_B= \psi_p(B)\in H^p(X, \Omega^q_X)$.\\[5pt]
1. The cup product $\omega_A\cup \omega_B\in H^n(X, \Omega^n_X)$ is  represented by the cochain 
\[
\{\pi(AB\Omega_{Rs}\Omega_{sT}/F_RF_s^2F_T), |Rs|=q, |sT|=p, RsT=I\}_I\in \sC^n(\sV, \Omega^n_X).
\] 
2. Let $\sU=\{U_0,\ldots, U_{n+1}\}$ be the open cover of $\P^{n+1}$, where $U_j$ is the complement of $F_j=0$. The coboundary $\delta_X(\omega_A\cup \omega_B)\in H^{n+1}(\P^{n+1}, \Omega_{\P^{n+1}}^{n+1})$ with respect to the residue sequence
\[
0\to \Omega_{\P^{n+1}}^{n+1}\to \Omega_{\P^{n+1}}^{n+1}(\Log X)\to i_*\Omega^n_X\to 0
\]
is represented by the cochain
\[
(-\dg)\cdot \frac{AB\Omega}{F_0\cdots F_{n+1}}\in \sC^{n+1}(\sU, \Omega_{\P^{n+1}}^{n+1}).
\]
\end{proposition}

\begin{proof} The proof is exactly the same as in \cite{CG}. We lift $A$ and $B$  to $k[X_0,\ldots, X_{n+1}]$, which we also denote as  $A$, $B$. (1) is a direct consequence of Proposition~\ref{prop:Rep} and the definition of cup product in terms of \v{C}ech cochains.  

For (2), we represent $\delta_X(\omega_A\cup \omega_B)$ as a \v{C}ech cochain by lifting 
\[
\pi(AB\Omega_{Rs}\Omega_{sT}/F_RF_s^2F_T)
\]
 to the section $dF/F\wedge AB\Omega_{Rs}\Omega_{sT}/F_RF_s^2F_T$ of  $\Omega_{\P^{n+1}}^{n+1}(\Log X)$ over $U_{RsT}$ and then taking the negative of the \v{C}ech coboundary. Following \cite[Lemma, pg. 14]{CG}, we have the identity in $\tilde{i}^*\Omega^{n+1}_{\A^{n+2}}$
\[
dF\wedge\Omega_{Rs}\Omega_{sT}=(-1)^\nu X_\nu F_s\Omega
\]
where $\nu$ is the complement of $RsT$ in $\{0,\ldots, n+1\}$ and $\tilde{i}:\tilde{X}\to \A^{n+2}$ is the inclusion of the affine cone over $X$. Noting that there is a single index $(01\ldots n+1)$ for $\sC^{n+1}(\sU,-)$,  the \v{C}ech coboundary is given by
\begin{align*}
\delta(\{\frac{dF}{F}\wedge \frac{AB\Omega_{Rs}\Omega_{sT}}{F_RF_s^2F_T}\})&=\sum_{\nu=0}^{n+1} \frac{AB\cdot X_\nu\Omega}{F\cdot\prod_{j\neq \nu}F_j}\\
&=\frac{\sum_{\nu=0}^{n+1}AB\cdot X_\nu F_\nu\Omega}{F\cdot F_0\cdots F_{n+1}}\\
&=\dg\cdot \frac{AB\cdot \Omega}{F_0\cdots F_{n+1}}\in \Omega^{n+1}_{\P^{n+1}}(U_{0\ldots n+1}).
\end{align*}
\end{proof}
 
It remains to compute $\Tr_{X/k}(\omega_A\cup\omega_B)$. Since the trace map is compatible with push-forward, we have
\[
\Tr_{X/k}(\omega_A\cup\omega_B)=\Tr_{\P^{n+1}/k}(i_{X*}(\omega_A\cup\omega_B))
\]
where $i_{X*}:H^n(X, \Omega_X^n)\to H^{n+1}(\P^{n+1}, \Omega^{n+1}_{\P^{n+1}})$ is the Gysin map. 

Choosing  $a_{ij}\in k[X_0,\ldots, X_{n+1}]$ with $F_i=\sum_{j=0}^{n+1}a_{ij}X_j$, the Scheja-Storch element $e_F\in J(F)_{(\dg-2)(n+2)}$ is the image in $J(F)$ of $\det(a_{ij})$. Setting $F_{ij}:=\del^2F/\del X_i\del X_j$, let $\Hess(F)$  denote the Hessian matrix $(F_{ij})$ of $F$.  

\begin{lemma}\label{lem:TrSS} Let $U_i\subset \P^{n+1}$ be the complement of $F_i=0$, giving the cover $\sU=\{U_0,\ldots, U_{n+1}\}$ of $\P^{n+1}$. Let $[e_F\Omega/\prod_{i=0}^{n+1}F_i]\in H^{n+1}(\P^{n+1}, \Omega^{n+1}_{\P^{n+1}})$ denote the class represented by the cocycle $e_F\Omega/\prod_{i=0}^{n+1}F_i\in \sC^{n+1}(\sU,   \Omega^{n+1}_{\P^{n+1}})$.  Then 
\[
\Tr_{\P^{n+1}/k}([e_F\Omega/\prod_{i=0}^{n+1}F_i])=1.
\]
\end{lemma}

\begin{proof} We first reduce to the case of characteristic zero. Let $\sO$ be a dvr  with residue field $k$ and quotient field $K$. Lift $F\in k[X_0,\ldots, X_{n+1}]$ to a homogeneous $F_\sO\in \sO[X_0,\ldots, X_{n+1}]_\dg$,  and let $\sX\subset \P^{n+1}_\sO$ be the hypersurface defined by $F_\sO$. Then $\sX$ is smooth over $\Spec \sO$ and thus there are homogeneous elements $a^{\sO}_{ij}\in \sO[X_0,\ldots, X_{n+1}]_{\dg-2}$ with $F_{\sO i}:=\del F_\sO/\del X_i=\sum_ja_{ij}^{\sO}X_j$. The image of $\det(a_{ij}^\sO)$ in $J(F_\sO)$ gives us the Scheja-Storch element $e_{F_\sO}\in J(F_\sO)$. Letting $\sU_\sO=\{U_{\sO 0}, \ldots, U_{\sO n+1}\}$ be the open cover of $\P^{n+1}_\sO$ with $U_{\sO i}$ the complement of $F_{\sO i}=0$,  we have the \v{C}ech cocycle  $e_{F_\sO}\Omega/\prod_{i=0}^{n+1}F_{\sO i}\in \sC^{n+1}(\sU_\sO, \Omega^{n+1}_{\P^{n+1}_\sO/\sO})$ representing a cohomology class $\rho_\sO\in H^{n+1}(\P^{n+1}_\sO, \Omega^{n+1}_{\P^{n+1}_\sO/\sO})$. Since the canonical trace map is well-defined for an arbitrary noetherian base-scheme and is compatible with base-change, we have the element $\Tr_{\P^{n+1}/\sO}(\rho_\sO)\in \sO$ which reduces to $\Tr_{\P^{n+1}/k}([e_F\Omega/\prod_{i=0}^{n+1}F_i])$ in $k$ and yields $\Tr_{\P^{n+1}_K}([e_{F_K}\Omega/\prod_{i=0}^{n+1}F_{Ki}])$ in $K$ upon localization. Thus, if we show that $\Tr_{\P^{n+1}/K}([e_{F_K}\Omega/\prod_{i=0}^{n+1}F_{Ki}])=1$, it will follow that 
$\Tr_{\P^{n+1}/k}([e_F\Omega/\prod_{i=0}^{n+1}F_{i}])=1$. Taking $\sO$ such that $K$ has characteristic zero and changing notation, we reduce to the case of a characteristic zero base-field $k$.

The element
\[
d\ln F_1/F_0\wedge\cdots\wedge d\ln F_{n+1}/F_n\in \sC^{n+1}(\sU, \Omega^{n+1}_{\P^{n+1}})
\]
represents $c_1(\sO_{\P^{n+1}}(\dg-1))^{n+1}\in H^{n+1}(\P^{n+1}, \Omega^{n+1}_{\P^{n+1}})$, and thus
\[
\Tr_{\P^{n+1}/k}([d\ln F_1/F_0\wedge\cdots\wedge d\ln F_{n+1}/F_n])= (\dg-1)^{n+1}. 
\]
Explicitly,  we have
\[
d\ln F_1/F_0\wedge\cdots\wedge d\ln F_{n+1}/F_n=\frac{\sum_{i=0}^{n+1}(-1)^iF_idF_0\cdots\widehat{dF_i}\cdots dF_{n+1}}{F_0\cdots F_{n+1}}.
\]
Using the Euler relation $(\dg-1)F_i=\sum_{j=0}^{n+1}F_{ij}\cdot X_j$ and some elementary properties of the determinant gives the identity
\[
(\dg-1)\cdot \frac{\sum_{i=0}^{n+1}(-1)^iF_idF_0\cdots\widehat{dF_i}\cdots dF_{n+1}}{F_0\cdots F_{n+1}}=\frac{\det\Hess(F)\cdot\Omega}{F_0\cdots F_{n+1}}, 
\]
Indeed,  we have
\begin{align*}
(\dg-1)\cdot \sum_{i=0}^{n+1}(-1)^iF_idF_0\cdots&\widehat{dF_i}\cdots dF_{n+1}\\
&=
\sum_{i=0}^{n+1}(-1)^i \sum_{j=0}^{n+1}F_{ij}X_jdF_0\cdots\widehat{dF_i}\cdots dF_{n+1}\\
&=\sum_{j=0}^{n+1}(-1)^jX_j\cdot  \sum_{i=0}^{n+1}(-1)^{i+j} F_{ij} dF_0\cdots\widehat{dF_i}\cdots dF_{n+1}.
\end{align*}
The coefficient of $dX_0\cdots\widehat{dX_l}\cdots dX_{n+1}$ in $\sum_{i=0}^{n+1}(-1)^{i+j} F_{ij} dF_0\cdots\widehat{dF_i}\cdots dF_{n+1}$ is $\det\Hess F$ for $l=j$, and for $l\neq j$ this is   $(-1)^{l-j}$ times the determinant of the matrix formed from $\Hess(F)$ by removing the $l$th column  and replacing it with the $j$th column, hence zero. 

Thus 
\[
\Tr_{\P^{n+1}/k}([\frac{\det\Hess(F)\cdot\Omega}{F_0\cdots F_{n+1}}])=(\dg-1)^{n+2}. 
\]
On the other hand, we have 
\[
F_i=\sum_{j=0}^{n+1}\frac{F_{ij}}{\dg-1}\cdot X_j
\]
so taking $a_{ij}=\frac{1}{\dg-1}F_{ij}$ gives the Scheja-Storch element 
\[
e_F=\frac{1}{(\dg-1)^{n+2}}\cdot \det\Hess(F) 
\]
and thus
\[
\Tr_{\P^{n+1}/k}([\frac{e_F\cdot \Omega}{F_0\cdots F_{n+1}}])=1.
\]
\end{proof}

\begin{theorem}\label{thm:Main} Let $k$ be a field and let $X\subset \P^{n+1}_k$ be a smooth hypersurface of degree $\dg$ prime to $\Char k$ and defined by a homogeneous $F\in k[X_0,\ldots, X_{n+1}]$. Take non-negative integers $p, q$ with $p+q=n$, take $A\in J(F)_{(q+1)\dg-n-2}$, $B\in J(F)_{(p+1)\dg-n-2}$ and let $\omega_A=\psi_q(A)\in H^q(X, \Omega_X^p)$, $\omega_B=\psi_p(B)\in H^p(X, \Omega_X^q)$. Let $e_F\in J(F)_{(\dg-2)(n+2)}$ be the Scheja-Storch element and suppose that $AB=\lambda\cdot e_F$ in $J(F)$, $\lambda\in k$. Then 
\[
\Tr_{X/k}(\omega_A\cup\omega_B)=(-\dg)\cdot\lambda\in k.
\]
\end{theorem}

\begin{proof} This follows from Proposition~\ref{prop:Product}, Lemma~\ref{lem:ResiduePushforward} applied to the push-forward $i_*:
H^n(X, \Omega^n_X)\to H^{n+1}(\P^{n+1},\Omega^{n+1}_{\P^{n+1}})$, and Lemma~\ref{lem:TrSS}.
\end{proof}

\begin{corollary}\label{cor:HdgJacQ} Let $k$, $X$ and $F$ be as above. Then  the isomorphism 
\[
\oplus_q\psi_q:\oplus_{q=0}^n J(F)_{\dg(q+1)-n-2}\to \oplus_{p+q=n}H^q(X,\Omega_X^p)_{prim}
\]
induces an isometry
\[
B_\Hodge\cong \begin{cases} -\dg\cdot B_\Jac\perp\<\dg\>&\text{ for $n$ even,}\\
-\dg\cdot B_\Jac&\text{ for $n$ odd,}
\end{cases}
\]
where $\<\dg\>$ is the rank one form $(x,y)\mapsto \dg\cdot xy$. 
\end{corollary}

\begin{proof} We take $p+q=n$. Recall that $H^p(X, \Omega^q_X)=H^p(X, \Omega^q_X)_{prim}$ if $p\neq q$ and if $p=q=n/2$, then $H^p(X, \Omega^p_X)_{prim}$ is the orthogonal complement of $k\cdot c_1(\sO_X(1))^p$ in $H^p(X, \Omega^p_X)$, with respect to $B_\Hodge$. Thus, if $n$ is odd, it follows from 
Theorem~\ref{thm:Main} that 
\[
B_\Hodge\cong -\dg\cdot B_\Jac.
\]

If $n=2m$ is even, let $B_\Hodge^{m,m}$ denote the restriction of $B_\Hodge$ to $H^m(X, \Omega^m_X)$ and $B_{\Hodge\ prim}^{m,m}$  the restriction of $B_\Hodge$ to $H^m(X, \Omega^m_X)_{prim}$. As $\Tr(c_1(\sO_X(1))^n)=\dg$, we see that
\[
B_\Hodge^{m,m}=\<\dg\>+B_{\Hodge\ prim}^{m,m}
\]
Using Theorem~\ref{thm:Main} again, we find
\[
B_\Hodge\cong \<\dg\>+ -\dg\cdot B_\Jac.
\]
This completes the proof.
 \end{proof}
 
 \begin{remark} Referring to  Remark~\ref{rem:Sign}, the reader may prefer the convention for the boundary map given by $\delta^m([z])=+[x]$, rather than the one we have been using, $\delta^m([z])=-[x]$. This yields the following changes:\\[5pt]
(i) In Lemma~\ref{lem:ResiduePushforward}, we have $\delta^{p,q}_X=-i_*$\\[2pt]
(ii) In Proposition~\ref{prop:Rep}, the cohomology class $\psi_q(\bar{A})\in H^q(X, \Omega^{n-q}_X)$ is represented by the cochain $\{(-1)^{q}\pi(A\Omega_I/F_I)\}\in \sC^q(\sV, \Omega^{n-q}_X)$.\\[2pt]
(iii) In Proposition~\ref{prop:Product}(2), $\delta_X(\omega_A\cup \omega_B)$ is represented by $(-1)^n\cdot \dg\cdot \frac{AB\Omega}{F_0\cdots F_{n+1}}$; this follows from the change (ii) and the change of sign in our coboundary convention.\\[2pt]
(iv) The conclusion of Theorem~\ref{thm:Main}  is changed by replacing $(-\dg)$ with $(-1)^{n+1}\dg$, since (i) changes the trace of the term in (iii)  by a factor of $(-1)$. \\[2pt]
(v) Corollary~\ref{cor:HdgJacQ} remains the same. If $n$ is even, this is clear. If $n$ is odd, then since $B_\Hodge$ pairs $H^p(X,\Omega^q_X)$ with $H^q(X,\Omega^p_X)$, and is zero when restricted to each of these summands if $p\neq q$, we see that $B_\Hodge$ is hyperbolic, so $(-1)\cdot B_\Hodge\cong B_\Hodge$. 
\end{remark}

 \section{Weighted-homogeneous hypersurfaces}
 
 We extend the results of the previous section to hypersurfaces in a weighted projective space, with an eye toward conductor formulas for  quasi-homogenous singularities. 
 
As mentioned in the introduction, much of what is covered here was done by Steenbrink \cite{Steenbrink} and also by Dolgachev \cite{Dolgachev}. We include this material to establish the notation and basic constructions and also to give details on the multiplicative structure, which we did not find in the literature. We refer the reader to  Delorme \cite{Delorme} for  foundational material on weighted projective spaces. 
 
Fix a sequence of non-negative integers $a_*:=(a_0,\ldots, a_{n+1})$ and a base-scheme $B$, and let $\rho_{a_*}:\G_{m,B}\to \Aut \A_B^{n+2}\setminus\{0\}$ be the representation 
\[
\rho_{a_*}(t)(y_0,\ldots, y_{n+1})=(t^{a_0}y_0,\ldots, t^{a_{n+1}}y_{n+1})
\]
We define $\P_B(a_*)$ to be the quotient scheme $\G_m\backslash (\A_B^{n+2}\setminus\{0\})$ with quotient map $p_{a_*}:\A_B^{n+2}\setminus\{0\}\to\P^n_B(a_*)$.  For $a_*=(1^{n+2})$,  we recover the usual projective space $\P^{n+1}_B$. 

Alternatively, we give the sheaf of $\sO_B$-algebra $\sO_B[Y_0,\ldots, Y_{n+1}]$ the grading with $\deg(Y_i)=a_i$, then $\P_B(a_*)=\Proj_{\sO_B}\sO_B[Y_0,\ldots, Y_{n+1}]$.

We can also describe $\P_B(a_*)$ as a quotient of $\P^{n+1}_B$.  Let $\mu(a_*)=\prod_{i=0}^n\mu_{a_i}$, with $\mu_\ell$ the group-scheme of $\ell$th roots of 1. For simplicity, we will henceforth assume that $a:=\prod_ia_i$ is prime to all residue characteristics of $B$. We have the map of sheaves of polynomial rings $\pi_{a_*}^*:\sO_B[Y_0,\ldots, Y_{n+1}]\to 
\sO_B[X_0,\ldots, X_{n+1}]$ with $\pi_{a_*}(Y_i)=X_i^{a_i}$, which intertwines the action $\rho_{a_*}$ and $\rho_{1^{n+2}}$ and gives the map of quotient schemes $\pi_{a_*}:\P^{n+1}_B\to \P_B(a_*)$, equivalently, the induced map $\Proj_{\sO_B}\sO_B[X_0,\ldots, X_{n+1}]\to \Proj_{\sO_B}\sO_B[Y_0,\ldots, Y_{n+1}]$.  We have the action $\bar{\rho}_{a_*}:\mu(a_*)\to \Aut\P^{n+1}_B$ by 
\[
\bar{\rho}_{a_*}(\lambda_0,\ldots, \lambda_{n+1})[x_0,\ldots, x_{n+1}]=
[\lambda_0x_0,\ldots,\lambda_{n+1} x_{n+1}]
\]
which identifies $\pi_{a_*}$ with the quotient map $\P^{n+1}_B\to \mu(a_*)\backslash \P^{n+1}_B$. 

We define the sheaf $\tilde\Omega^p_{\P(a_*)}$ via a modified Euler sequence. Given a sheaf $\sF$ on $\A^{n+1}\setminus\{0\}$ with a $\G_m$-action, and an integer $m$, we let $\sF(m)$ be the sheaf with $\G_m$-action twisted by the characters $\chi_m(t)=t^m$. The vector field $\Eul(a_*):=\sum_{i=0}^{n+1}a_iY_i\del/\del Y_i$ is invariant with respect to $\rho_{a_*}$.  We define $\hat\Omega^p_{\P(a_*)}$ as the kernel of
\[
\bigwedge^p_{\sO_{\A^{n+2}}}(\oplus_{i=0}^{n+1}\sO_{\A^{n+2}}(-a_i)e_i)\xrightarrow{\iota^p_{\Eul(a_*)}}
\bigwedge^{p-1}_{\sO_{\A^{n+2}}}(\oplus_{i=0}^{n+1}\sO_{\A^{n+2}}(-a_i)e_i)
\]
where $\iota^1_{\Eul(a_*)}(e_i)=a_iY_i$, extended as usual to the map $\iota^p_{\Eul(a_*)}$. The induced    $\G_m$-action on $\hat\Omega^p_{\P(a_*)}$ defines the sheaf $\tilde\Omega^p_{\P(a_*)}$ on $\P_B(a_*)$ as the $\G_m$-invariants of $\hat\Omega^p_{\P(a_*)}$; the sheaf $\hat\Omega^p_{\P(a_*)}(m)$ defines  the twisted version $\tilde\Omega^p_{\P(a_*)}(m)$. 

As the $\mu(a_*)$-action on $\P^{n+1}_B$ arises from the $\GL_{n+2}$-action on $\A^{n+2}$, all $\GL_{n+2}$-linearized sheaves on $\P^{n+1}_B$, such as $\Omega^p_{\P^{n+1}_B/B}(m)$, have a natural 
$\mu(a_*)$-action. Via the description of $\P_B(a_*)$ as $\mu(a_*)\backslash\P^{n+1}_B$ we have another description of  $\tilde\Omega^p_{\P(a_*)}(m)$ as the $\mu(a_*)$-invariants of $\Omega^p_{\P^n_B}(m)$.  The fact that $\sO_B[X_0^{a_0},\ldots, X_n^{a_n}]\subset \sO_B[X_0,\ldots, X_n]$ is the subring of $\mu(a_*)$-invariants in $\sO_B[X_0,\ldots, X_n]$ shows that that our two descriptions of $\tilde\Omega^p_{\P(a_*)}(m)$ agree.

We note that the sheaf $\sO_{\P(a_*)}(m)$ on $\P(a_*)$ is an invertible sheaf if $m$ is divisible by the least common multiple    $\lcm(a_*)$   of $a_0,\ldots, a_{n+1}$: the section $X_i^{m/a_i}$ is a generator over the open subscheme $U_i(a_*)\subset \P(a_*)$ defined by $X_i\neq0$.

\begin{theorem}[Bott's theorem for weighted projective space \hbox{\cite[Theorem (2.3.2)]{Dolgachev}}] Suppose all the residue characteristics of  $B$ are  prime to $\prod_ia_i$.  The cohomology of $\tilde\Omega^q_{\P_B(a_*)}(m)$ satisfies
\[
H^p(\P_B(a_*), \tilde\Omega^q_{\P_B(a_*)/B}(m))=\begin{cases} 0&\text{ for }p>0\text{ and }m\ge q-n, m\neq0,\\0& \text{ for }p>0, m=0\text{ and }p\neq q\\
0&\text{ for }p=0,   m\le q, \text{ except for }m=p=q=0.
\end{cases}
\]
\end{theorem}

\begin{proof} We have the isomorphism
\[
H^p(\P_B(a_*), \tilde\Omega^q_{\P_B(a_*)/B}(m))\cong H^0(\mu(a_*), H^p(\P_B^{n+1}, \Omega^q_{\P^n/B}(m)))
\]
from which the result follows from Bott's theorem for $\P^n$.
\end{proof}

Suppose for simplicity that $B=\Spec A$ for some noetherian ring $A$. We now consider a sequence $a_*:=(a_0,\ldots, a_{n+1})$ and a positive integer $\dg$ divisible by $\lcm(a_*)$.  A homogeneous polynomial $F(Y_0,\ldots, Y_{n+1})\in A[Y_0,\ldots, Y_{n+1}]$ of weighted-degree $\dg$ defines a hypersurface  $Y_F\subset \P_A(a_*)$. Let $G:=\pi_{a_*}^*(F)\in A[X_0,\ldots, X_{n+1}]$, giving the degree $\dg$ hypersurface $X:=X_G\subset \P^{n+1}_A$.  We define
$\tilde\Omega^p_{\P(a_*)}(qY)$ as the $\G_m$-invariants in the sheaf $\hat\Omega^p_{\P(a_*)}(q\cdot p_{a_*}^{-1}(Y))$ of $p$-forms with poles along $p_{a_*}^{-1}(Y)\subset \A^{n+2}$. As in the classical case, we have the canonical isomorphism  
$\tilde\Omega^p_{\P(a_*)}(qY)\cong \tilde\Omega^p_{\P(a_*)}(q\cdot \dg)$. Alternatively, $\tilde\Omega^p_{\P(a_*)}(qY)$ is the sheaf of $\mu(a_*)$-invariants in $\Omega^p_{\P^{n+1}_A/A}(q\cdot X)$. 

\begin{definition}  Let $F\in A[Y_0,\ldots, Y_{n+1}]$ be a homogeneous polynomial of weighted-degree $\dg$ with $\dg$ divisible by $\lcm(a_*)$. 
We call $Y_F\subset \P_A(a_*)$ a {\em smooth quotient hypersurface} if the hypersurface $X_G\subset \P^{n+1}_A$ defined by $G$ is smooth over $A$. We call $Y_F\subset \P_A(a_*)$ a smooth hypersurface if in addition $Y_F$ is a smooth $A$-scheme.  
\end{definition}
In particular, if $Y_F\subset \P_A(a_*)$ is a  smooth quotient hypersurface, then $Y_F$ has at worst (tame) quotient singularities.

For $Y_F\subset \P_A(a_*)$ a smooth quotient hypersurface, the map $\pi_{a_*}:\P^{n+1}_A\to \P_A(a_*)$ induces a map $\pi_{a_*, F}:X_G\to Y_F$ identifying $Y_F$ with $\mu_{a_*}\backslash X_G$. We define the sheaf $\tilde\Omega^p_{Y_F/A}$ on $Y_F$ as the $\mu_{a_*}$-invariants in $\Omega^p_{X_G/A}$. The map $\pi_{a_*, F}$ gives   a corresponding map on Hodge cohomology
\[
\pi_{a_*,F}^*:H^q(Y_F, \tilde\Omega^p_{Y_F/A})\to H^q(X_G, \Omega^p_{X_G/A})
\]
giving an isomorphism
\[
H^q(Y_F, \tilde\Omega^p_{Y_F/A})\cong H^q(X_F, \Omega^p_{X_F/A})^{\mu_{a_*}}.
\]
If $Y_F$ is a smooth hypersurface over $A$, then the canonical map $\Omega^p_{Y_F/A}\to \tilde\Omega^p_{Y_F/A}$ is an isomorphism and $H^q(Y_F, \tilde\Omega^p_{Y_F/A})=H^q(Y_F, \Omega^p_{Y_F/A})$ is the usual Hodge cohomology of $Y_F$. 

We now take $A=k$ a field. Let $i:Y:=Y_F\to \P_A(a_*)$ be the inclusion of a smooth quotient hypersurface. The logarithmic differential $dF/F$ gives a well-defined section of $i^*\tilde\Omega_{\P(a_*)}^1(Y)$ and we have the residue map
\[
\res_Y:H^0(\P(a_*), \tilde\Omega^{n+1}_{\P(a_*)}((q+1)Y)\to H^0(Y, \tilde\Omega_Y^n(qY))
\]
both of these being induced by the corresponding constructions for $X$ by taking $\mu_{a_*}$-invariants. Similarly, Lemma~\ref{lem:ExactSeq}  extends to 
\begin{lemma}\label{lem:ExactSeqWeighted}  Suppose that $\dg\cdot\prod_ia_i$ is prime to the characteristic of $k$. \\[5pt]
1. We have an exact sequence 
\begin{multline*}
0\to \tilde\Omega^p_Y\xrightarrow{dF/F\wedge-} i^*\tilde\Omega_{\P(a_*)}^{p+1}(Y)\xrightarrow{dF/F\wedge-} i^*\tilde\Omega_{\P(a_*)}^{p+2}(2Y)\xrightarrow{dF/F\wedge-}\\\ldots\xrightarrow{dF/F\wedge -}i^*\Omega_{\P(a_*)}^n(qY)\to \tilde\Omega_Y^n(qY)\to 0
\end{multline*}
for  integers $p,q\ge0$ with $p+q=n$. \\[2pt]
2. The map $\delta:H^0(Y, \tilde\Omega_Y^n(qY))\to H^q(Y, \tilde\Omega^p_Y)$ induced by the exact sequence (1) gives an exact sequence
\[
H^0(Y, i^*\tilde\Omega_{\P^{n+1}}^n(qY))\xrightarrow{\pi_q} H^0(Y, \tilde\Omega_Y^n(qY))\xrightarrow{\delta} H^q(Y, \tilde\Omega^p_Y)
\]
for all $q\ge0$, $p+q=n$. Moreover $\delta$ is surjective if $p\neq q$ and has image 
$H^p(Y, \tilde\Omega^p_Y)_{prim}$.
\end{lemma}
Here $H^p(Y, \tilde\Omega_Y^p)_{prim}$ is defined as the  $\mu_{a_*}$-invariants in 
$H^p(X, \Omega_X^p)_{prim}$. In case $Y$ is smooth over $k$, this is the same as the usual notion of primitive, where we may use the image of  $c_1(O_{\P_A(a_*)}(m))$ in $H^1(Y, \Omega_Y^1)$ as the ample class for any $m$ prime to the characteristic with $\lcm(a_*)\ |\ m$.

We let $\Omega_{a_*}$ be the $n+1$-form on $\A^{n+2}_k$
\[
\Omega_{a_*}:=\sum_{i=0}^{n+1}(-1)^ia_iY_idY_0\wedge\ldots\wedge\widehat{dY}_i\wedge\ldots\wedge dY_{n+1}
\]
We note that 
\begin{equation}\label{eqn:OmegaPullback}
\pi_{a_*}^*(\Omega_{a_*})=\prod_{i=0}^{n+1}a_i\cdot \prod_{i=0}^{n+1}X_i^{a_i-1}\cdot \Omega
\end{equation}
where 
\[
 \Omega=\sum_{i=0}^{n+1}(-1)^iX_idX_0\wedge\ldots\wedge\widehat{dX}_i\wedge\ldots\wedge dX_{n+1}
 \]
 is the $n+1$-form we considered before. Letting $p:\A^{n+2}\setminus\{0\}\to \P^{n+1}$ be the quotient map, $\Omega$ is a local generator of $p^*\Omega^{n+1}_{\P^{n+1}/k}\subset 
 \Omega^{n+1}_{\A^{n+1}/k}$, and $\pi_{a_*}^*(\Omega(a_*))$ is a local generator of the $\mu(a_*)$-invariants in $p^*\Omega^{n+1}_{\P^{n+1}/k}$.  With respect to the $\G_m$ action on $\A^{n+2}\setminus\{0\}$ defining $\P(a_*)$,  $\Omega(a_*)$ has weight $|a_*|:=\sum_ia_i$. Thus, sending $A\in k[Y_0,\ldots, Y_{n+1}]_{(q+1)\cdot\dg-|a_*|}$, to the form $A\Omega(a_*)/F^{q+1}$ defines an isomorphism
 \[
\hat\psi_{Y,q}: k[Y_0,\ldots, Y_{n+1}]_{(q+1)\cdot\dg-|a_*|}\xrightarrow{\sim} H^0(\P(a_*), \tilde\Omega^{n+1}_{\P(a_*)}((q+1)Y).
 \]
 
As above, we take a weighted homogeneous polynomial $F\in k[Y_0,\ldots, Y_{n+1}]$ of degree $\dg$ and let $G=\pi_{a_*}^*F\in k[X_0,\ldots, X_{n+1}]$. We assume that $Y_F$ is a smooth quotient hypersurface in $\P(a_*)$.

Let  $F_i=\del F/\del Y_i$. The condition that $Y_F$ is a smooth quotient hypersurface in $\P(a_*)$ implies that the ideal $(F, F_0,\ldots, F_{n+1})$ is $(Y_0,\ldots, Y_{n+1})$-primary, so the Jacobian ring $J(F)$ is an Artinian $k$-algebra. Indeed, 
\[
G_i=\del G/\del X_i=\del F/\del Y_i\cdot (a_iX_i^{a_i-1}), 
\]
and $(G, G_0,\ldots, G_{n+1})$ is by assumption $(X_0,\ldots, X_{n+1})$-primary.   We have the modified Euler equation
\[
\dg\cdot F=\sum_i a_iY_i F_i
\]
so $(F, F_0,\ldots, F_n)=(F_0,\ldots, F_n)$ is a complete intersection ideal and $J(F)$ has a 1-dimensional socle generated by the Scheja-Storch element $e_F$. 

Define the map of $k$-vector spaces
\[
\tilde\theta_{a_*}^*: k[Y_0,\ldots, Y_{n+1}]\to k[X_0,\ldots, X_{n+1}]
\]
by 
\[
\tilde\theta_{a_*}^*(H):=\pi_{a_*}^*(H)\cdot \prod_{i=0}^na_i\cdot X_i^{a_i-1}
\]
One checks that $\tilde\theta_{a_*}^*$ descends to a well-defined map of $k$-vector spaces
\[
\theta_{a_*}^*: J(F) \to J(G)
\]
Indeed, $G_i=a_iX_i^{a_i-1}\pi_{a_*}^*(F_i)$, so $\theta_{a_*}^*(F_i)$ is in $(G_0,\ldots, G_{n+1})$. 

For a version of the following result in characteristic zero,  see \cite[Theorem 4.3.2]{Dolgachev}, where this result is attributed to Steenbrink.

\begin{proposition}\label{prop:JacobianQuotientHypersurf} Suppose $Y_F$ is a smooth quotient hypersurface, and $\dg\cdot\prod_ia_i$ is prime to $\Char k$. Let $|a_*|:=\sum_ia_i$.  Sending $A\in k[Y_0,\ldots, Y_{n+1}]_{(q+1)\dg-|a_*|}$ to 
$\delta\res_Y(A\Omega_{a_*}/F^{q+1})$ descends to an isomorphism
\[
\psi_{Y,q}: J(F)_{(q+1)\dg-|a_*|}\to H^q(Y, \tilde\Omega^p_Y)_{prim}
\]
Moreover, the diagram
\[
\xymatrix{
J(F)_{(q+1)\dg-|a_*|}\ar[r]^{\theta_{a_*}^*}\ar[d]_{\psi_Y,q}&J(G)_{(q+1)\dg-n-2}\ar[d]^{\psi_{X,q}}\\
H^q(Y, \tilde\Omega^p_Y)_{prim}\ar[r]^{\pi_{a_*}^*}&H^q(X, \Omega^p_X)_{prim}
}
\]
commutes.
\end{proposition}

\begin{proof} We follow the proof of Proposition~\ref{prop: Main}, taking $\mu(a_*)$-invariants. This shows that the map $\hat\psi_{Y,q}$ defines an isomorphism
\[
\bar\psi_{Y,q}:\left(k[Y_0,\ldots, Y_{n+1}]/(F)\right)_{(q+1)\cdot\dg-|a_*|}\xrightarrow{\sim} H^0(Y, i^*\tilde\Omega^{n+1}_{\P(a_*)}((q+1)Y)),
\]
we have the  commutative diagram
\[
\xymatrix{
 H^0(Y, i^*\tilde\Omega^n_{\P(a_*)}(qY))\ar[r]^{dF/F\wedge-}\ar[dr]_{\pi_q}& H^0(Y, i^*\tilde\Omega^{n+1}_{\P(a_*)}((q+1)X))\ar[d]^{\bar\res_Y}\\
&H^0(Y, \tilde\Omega^n_{Y}(qY))
}
\]
and the map $H^0(\P(a_*), \tilde\Omega^n_{\P(a_*)}(qY))\to H^0(Y, i^*\tilde\Omega^n_{\P(a_*)}(qY))$ is surjective.

We replace the Koszul resolution used in the proof of Proposition~\ref{prop: Main} with the modified version, where we use generators $a_0Y_0,\ldots, a_{n+1}Y_{n+1}$ for the ideal $(Y_0,\ldots, Y_{n+1})\subset R=k[Y_0,\ldots, Y_{n+1}]$. The corresponding Koszul complex, with $e_i$ having degree $a_i$,  is
\[
\ldots \to \bigwedge^{j+1}\oplus_{i=0}^{n+1} R(-a_i)e_i\xrightarrow{\Eul_{a_*}}
 \bigwedge^{j}\oplus_{i=0}^{n+1} R(-a_i)e_i\to \ldots, 
 \]
which breaks up into graded pieces. This   presents 
$H^0(\P(a_*), \tilde\Omega^n_{\P(a_*)}(m))$ as
\[
0\to R_{m-|a_*|}\to [\bigwedge^{n+1}\oplus_{i=0}^{n+1} R(-a_i)e_i]_m\to H^0(\P(a_*), \tilde\Omega^n_{\P(a_*)}(m))\to0
\]
with $\deg(e_i)=a_i$.
Taking $m=|a_*|-a_\ell$   gives us the element in $H^0(\P(a_*), \tilde\Omega^n_{\P(a_*)}(m))$, 
 \[
\Omega^{(\ell)}_{a_*}:=\Eul_{a_*}(e^{\hat{\ell}})=\sum_{i=0}^{\ell-1}(-1)^ia_iY_idY^{\widehat{i,\ell}}-\sum_{i=\ell+1}^{n+1}(-1)^ia_iY_idY^{\widehat{\ell,i}}
\]
Taking $m=q\cdot\dg$, we see that 
$H^0(\P(a_*), \tilde\Omega^n_{\P(a_*)}(qY))$ is generated by elements of the form
\[
\omega:=\frac{\sum_{\ell=0}^{n+1}A_\ell\cdot \Omega_{a_*}^{(\ell)}}{F^q};\quad A_\ell\in k[Y_0,\ldots, Y_{n+1}]_{q\cdot\dg-|a_*|+a_\ell}.
\]
One computes $dF\cdot \Omega^{(\ell)}_{a_*}=(-1)^{\ell-1}F_\ell\cdot \Omega_{a_*}+\dg\cdot F\cdot dY^{\hat{\ell}}$. For $\omega$ as above, we thus have
\[
\frac{dF}{F}\wedge i^*\omega=\frac{\sum_\ell (-1)^{\ell-1}A_\ell\cdot F_\ell}{F^{q+1}}\cdot \Omega_{a_*}.
\]
This shows that the kernel of $\delta\bar\res_Y\circ\hat\psi_{Y,q}$ is exactly the degree $(q+1)\dg-|a_*|$ part of the ideal $(F, F_0,\ldots, F_{n+1})/(F)$, that is, $\psi_{Y,q}:J(F)_{(q+1)\dg-|a_*|}\to H^q(Y, \tilde\Omega^p_Y)_{prim}$ is an isomorphism.

The commutativity of  the diagram in the statement of the Proposition follows from the identity \eqref{eqn:OmegaPullback}.
\end{proof}

We can now complete our comparison result for smooth hypersurfaces $Y_F\subset \P(a_*)$. As before, we have the non-degenerate bilinear forms $B_{Y,\Hodge}$ on $\oplus_{p+q=n}H^q(Y,\Omega^p_{Y/k})$ and $B_{Y,\Jac}$  on $\oplus_{q=0}^n J(F)_{(q+1)\cdot\dg-|a_*|}$.

\begin{theorem}\label{thm:HdgJacQH} Let $a_*=(a_0,\ldots, a_{n+1})$ be a sequence of non-negative integers with $\gcd=1$ and let  $F\in k[Y_0,\ldots, Y_{n+1}]$ be a weighted homogeneous polynomial of degree $\dg$, where we give $Y_i$ degree $a_i$. Suppose that $\dg\cdot\prod_{i=0}^{n+1}a_i$ is  prime to the characteristic of $k$ and that the hypersurface $Y\subset\P_k(a_*)$ defined by $F$ is smooth over $k$. Then the isomorphisms $\psi_{Y,q}: J(F)_{(q+1)\cdot\dg-|a_*|}\to H^q(Y, \Omega^{n-q}_{Y/k})_{prim}$ induce an isometry 
\[
B_{Y,\Hodge}\cong \<\dg\cdot \prod_{i=0}^{n+1}a_i\>+\<-\dg\>B_{\Jac(F)} 
\]
in case $n$ is even and 
\[
B_{Y,\Hodge}\cong  \<-\dg\>B_{\Jac(F)}  
\]
in case $n$ is odd.
\end{theorem}

\begin{proof} Let $G=\pi_{a_*}^*(F)\in k[X_0,\ldots, X_{n+1}]$ and let $X\subset \P^{n+1}_k$ be the smooth hypersurface defined by $G$.  We have the finite morphism $\pi_{a_*}:X\to Y$. We  begin by comparing the Scheja-Storch elements $e_F\in J(F)_{\dg\cdot (n+2)-2|a_*|}$ and $e_G\in J(G)_{(\dg-2)(n+2)}$.  

Write $\del F/\del Y_i=\sum_ja_{ij}\cdot Y_j$. Then
\[
\del G/\del X_i=\pi_{a_*}^*(\del F/\del Y_i)\cdot a_iX_i^{a_i-1}=\sum_j\pi_{a_*}^*(a_{ij})\cdot X_j^{a_j-1}a_iX_i^{a_i-1}\cdot X_j
\]
Letting $\text{diag}(\alpha_0,\ldots, \alpha_{n+1})$ denote the diagonal matrix with entries $\alpha_i$, we have
\begin{align*}
e_G&=\det(\pi_{a_*}^*(a_{ij})\cdot X_j^{a_j-1}a_iX_i^{a_i-1})\\&=
\det\left(\text{diag}(\ldots, a_iX_i^{a_i-1}, \ldots)\cdot
(\pi_{a_*}^*(a_{ij}))\cdot \text{diag}(\ldots, X_i^{a_i-1}, \ldots)\right)\\&=
(\prod_ia_i)\cdot(\prod_iX_i^{a_i-1})^2\pi_{a_*}^*(e_F).
\end{align*}

Suppose we have $x\in J(F)_{(q+1)\dg-|a_*|}$, $y\in J(F)_{(p+1)\dg-|a_*|}$ with $p+q=n$ and with
\[
xy=\lambda\cdot e_F,
\]
in other words,  $B_{\Jac(F)}(x,y)=\lambda$. Then
\[
\theta_{a_*}^*(x)\cdot \theta_{a_*}^*(y)=(\prod_ia_iX_i^{a_i-1})^2\pi_{a_*}(xy)=\lambda\cdot\prod_ia_i\cdot  e_G
\]
that is, 
\[
B_{\Jac(G)}(\theta_{a_*}^*(x), \theta_{a_*}^*(y))=(\prod_ia_i)\cdot   B_{\Jac(F)}(x, y)
\]

On the other hand, for $\omega_x\in H^q(Y, \Omega^p_Y)_{prim}$, $\omega_y\in H^p(Y, \Omega^q_Y)_{prim}$ we have
\[
B_{Y,\Hodge}(\omega_x, \omega_y)=\Tr_{Y/k}(\omega_x\cup\omega_y)
\]
Since $\Tr_{X/k}=\Tr_{Y/k}\circ \pi_{a_*, *}$, the projection formula gives us
\begin{multline*}
B_{X,\Hodge}(\pi_{a_*}^*(\omega_x), \pi_{a_*}^*(\omega_y))=\Tr_X(\pi_{a_*}^*(\omega_x)\cup\pi_{a_*}^*(\omega_y))\\=\deg(\pi_{a_*})\cdot \Tr_{Y/k}(\omega_x\cup\omega_y)=
(\prod_ia_i)\cdot B_{Y,\Hodge}(\omega_x, \omega_y).
\end{multline*}

Since the restriction of $B_{X,\Hodge}$ to $\oplus_{p+q=n}H^q(X, \Omega^p_X)_{prim}$ is isometric via $\psi_X$ to $\<-\dg\>\cdot B_{\Jac(G)}$ (Theorem~\ref{thm:Main}), we see with the help of Proposition~\ref{prop:JacobianQuotientHypersurf} that the restriction of $B_{Y,\Hodge}$ to 
$\oplus_{p+q=n}H^q(Y, \Omega^p_Y)_{prim}$ is isometric via $\psi_Y$ to 
$\<-\dg\>\cdot B_{\Jac(F)}$.

In case $n$ is odd, this completes the computation. For $n$ even, we need to compare the top self-intersection numbers of $c_1(\sO_Y(m))$ on $Y$ and $c_1(\sO_X(m))$ on $X$; up to a square, this is independent of the choice  of $m$ (prime to the characteristic). Since $\pi_{a_*}^*(\sO_Y(m))=
\sO_X(m)$, using the projection formula we see that 
\[
\deg_k(c_1(\sO_X(m))^n)=(\prod_ia_i)\deg_k(c_1(\sO_Y(m))^n).
\]
This completes the proof.
\end{proof}

\begin{remark} As above, consider a hypersurface $Y_F\subset \P_k(a_*)$,  $a_*=(a_0,\ldots, a_{n+1})$, with $\gcd(a_0,\ldots, a_{n+1})=1$, $\dg\cdot\prod_ia_i$ prime to the characteristic and $\lcm(a_*)\ |\ \dg$.  For $n\ge2$, the condition that $Y_F\subset \P_k(a_*)$  is smooth over $k$ implies that $\P_k(a_*)$  has only isolated singularities; indeed, the assumption that $\lcm(a_*)\ |\ \dg$ implies that $Y_F$ is a Cartier divisor on $\P_k(a_*)$, so if $Y_F$ is smooth over $k$, then $Y_F$ must be contained in the smooth locus of $\P_k(a_*)$ and since $Y_F$ is ample, the singular locus must have dimension zero. In fact, the singular locus must then be contained in the set of ``vertices'' $(0,\ldots, 0,1,0\ldots, 0)$ of $\P_k(a_*)$.

For $n\ge3$, this is equivalent with the $a_i$ being pairwise relatively prime and implies that $F=\sum_{i=0}^{n+1}\lambda_i\cdot Y_i^{\dg/a_i}+\ldots$ with the $\lambda_i\neq0$. In particular $\dg=m\cdot \prod_ia_i$ for some integer $m$, and one can rewrite the ``constant'' term $\<\dg\cdot \prod_ia_i\>$ in the case of even $n>2$ as $\<\dg/\prod_ia_i\>$.  
\end{remark}

\begin{ex} A common example of a smooth weighted homogenous hypersurface is given by the $d$-fold cyclic cover $Y$ of $\P^n$ branched along a smooth hypersurface of degree $md$:
\[
F(Y_0,\ldots, Y_{n+1})=Y_{n+1}^d-f_{md}(Y_0,\ldots, Y_n)
\]
In this case $a_*=(1,\ldots, 1, m)$, $\dg=md$, $e_F=(-d)Y_{n+1}^{d-2}\cdot e_f$, $J(F)=J(f)\otimes J(Y_{n+1}^d)$. The bilinear form $B_{\Jac(F)}$ does not however decompose as $\<-d\>\cdot B_{\Jac(f)}$ (even modulo hyperbolic factors): for $d=2l+2$, modulo hyperbolic terms, $B_{\Jac(F)}$ arises from the pairings, with $p+q=n$
\[
 J(f)_{(q+1)\dg-lm-(n+m+1)}\times J(f)_{(p+1)\dg-lm-(n+m+1)}\to J(f)_{(n+1)(\dg-2)}=k\cdot e_f
\]
Calling this $B_{\Jac(f)'}$, we have $B_{\Jac(F)}=\<-d\>B_{\Jac(f)'}+r\cdot H$ for some $r$. However, since $f=0$ has dimension $n-1$, these are not the graded pieces involved in the pairing $B_{\Jac(f)}$. 

\end{ex}

 \section{Quadratic conductor formulas} Before attempting a general conjecture, we would like to work out a few basic examples. The first one is the case of a finite flat morphism $f:X\to S$ in $\Sch_k$, with $X=\Spec A$, $S=\Spec \sO$, and both $A$ and $\sO$ discrete valuation rings, with $f^*:\sO\to A$ a local homomorphism, in other words, an exercise in elementary ramification theory. We first determine which data are canonical and where we need to make choices. We assume that $\sO$ has residue field $k$, that the ramification is tame and that $k$ has characteristic different from two.
  
 Let $K$ be the quotient field of $\sO$. Then $K\to A_K:=A\otimes_\sO K$ is a finite   field extension , which we will assume is separable, giving the  non-degenerate trace form $q_{A_K/K}(x)=\Tr_{A_K/K}(x^2)$.  The class of $q_{A_K/K}$ in $\GW(K)$ is the motivic Euler characteristic $\chi(\Spec A_K/K)$. Letting $k'$ be the residue field of $A$, we similarly have the trace form $q_{k'/k}(x)=\Tr_{k'/k}(x^2)$ representing $\chi(\Spec k'/k)\in \GW(k)$. Both $\chi(\Spec A_K/K)\in \GW(K)$ and $\chi(\Spec k'/k)\in \GW(k)$ are canonical, but in contrast to the geometric situation, they live in different rings, so cannot be compared directly.

\begin{remark}\label{rem:specialization}
We recall the definition of the specialization map
\[
\spc_t:\GW(K)\to \GW(k).
 \]
 
Let $F$ be a field and let $K^{MW}_{*}(F)$ denote the graded Milnor-Witt K-theory algebra of $F$ \cite[Definition 3.1]{Morel-book}. For $u\in F^\times$, we have the element $[u]\in K^{MW}_1(F)$; the elements $[u]$ together with the element $\eta\in K^{MW}_{-1}(F)$ generate $K^{MW}_{*}(F)$ as a ring.  Morel \cite[Lemma 6.3.8]{MorelICTP} defines an isomorphism of rings
\[
\phi_{0,F}:\GW(F)\xrightarrow{\sim} K^{MW}_{0}(F).
\]
with $\phi_{0,F}(\<u\>)=1+\eta[u]$; since the rank one forms $\<u\>$, $u\in F^\times$, generate 
$\GW(F)$ as an abelian group, this formula determines $\phi_{0,F}$.

To the dvr  $\mathcal{O}$ with parameter $t$, Morel associates a specialization morphism on Milnor-Witt K-theory \cite[Lemma 3.16, discussion on pg. 57]{Morel-book}:
\[
s^{t}_{v}:K^{MW}_{*}(K)\to K^{MW}_{*}(k),
\]
which is the unique homomorphism of graded rings such that $s^{t}_{v}(\eta)=\eta$, $s^{t}_{v}([t^nu])=\bar{u}$ for $u\in \mathcal{O}^{\times}$, $n\in \Z$, where $\bar{u}\in k^\times$ is the image of $u$ under the quotient map $\pi:\sO\to k$. Here $v$ is the induced valuation on $K$.

Via the isomorphisms $\phi_{0,K}$, $\phi_{0,k}$, the specialization morphism $s^{t}_{v}:K^{MW}_0(K)\to K^{MW}_0(k)$ defines the ring homomorphism $\spc_t:\GW(K)\to \GW(k)$, and characterizes $\spc_t$ as the unique ring homomorphism with $\spc_t(\<t^nu\>)=\<\bar{u}\>$ for 
all $n\in \Z$, $u\in \sO^\times$.

Specialization maps on the Witt ring are classical although the setting of the Grothendieck-Witt ring does not seem to have been considered; the two are essentially equivalent constructions. The case of a dvr  can be derived from Springer's theorem  \cite[Propositions 3,   4 and 5]{Springer} on the structure of the Witt ring of a henselian, discretely valued field.  Knebusch  \cite[Theorem 3.1]{Knebusch} constructs a specialization map for a general valuation. See \cite[Lemma 19.10]{EKM} for a recent treatment of the case of a dvr .   These all have the property that the specialization of $\<t^nu\>\in W(K)$ is $\<\bar{u}\>\in W(k)$ for all $n\in \Z$, $u\in \sO^\times$, so they agree with the map $W(K)\to W(k)$ induced by $\spc_t$.
\end{remark}

 Thus, we can consider
\[
\Delta_t(A/\sO):=\spc_t(\chi(A_K/K))-\chi(k'/k)\in \GW(k).
\]

If we set $n=[A_K:K]$, $f:=[k':k]$ and let $e$ be the usual ramification index: $t=u\cdot s^e$ where $s\in A$ is a parameter and $u\in A^\times$, we have
\[
\rnk(\Delta_t(A/\sO))=n-f=(e-1)\cdot f.
\]
Since we are in the tamely ramified case, $(e-1)\cdot f$ is just the degree over $k$ of the zero-cycle given by the vanishing of the section $dt\in \Omega_{A/k}$. For the quadratic version, writing $t=us^e$ gives us $dt=e\cdot u\cdot s^{e-1}ds+ s^e\cdot du$. As $ds$ is a generator for $\Omega_{A/k}$ as $A$-module, we have
\[
dt=v\cdot s^{e-1}ds
\]
for some $v\in A^\times$ with $v\equiv e\cdot u\mod s$. Letting $0_A=\Spec k'\subset \Spec A$, the section $dt$ of $\Omega^1_{A/k}$ gives us the corresponding Euler class with support $e_{0_A}(\Omega^1_{A/k}, dt)\in H^1_{0_A}(\Spec A, \sK^{MW}_1(\omega^{-1}_{A/k}))$. Explicitly, this is just
\[
\<\bar{v}\>\cdot\sum_{i=0}^{e-2}\<-1\>^i\otimes ds^{-1}\otimes \del/\del s
\]
where $\bar{v}$ is the image of $v$ in $k'$ (see \cite[Example 4.5]{LevineEuler}). The purity isomorphism $H^1_{0_A}(\Spec A, \sK^{MW}_1(\omega^{-1}_{A/k}))\cong \GW(k')$ transforms this element to $\<\bar{v}\>\cdot\sum_{i=0}^{e-2}\<-1\>^i\in \GW(k')$. To sum up, this element of $\GW(k')$ does depend on the choice of parameter $t$ for $\sO$, but not on the choice of parameter $s$ for $A$. Note that 
$\rnk(e_{0_A}(\Omega^1_{A/k}, dt))=e-1$, so
\[
\rnk(\Tr_{k'/k}(e_{0_A}(\Omega^1_{A/k}, dt)))=(e-1)\cdot f=\rnk(\Delta_t(A/\sO)).
\]
This is just the classical, albeit trivial, conductor formula in relative dimension zero: $[A_K:K]-[k':k]=\deg_k(dt=0)$. 

To compute the quadratic version, note that the specialization $\spc_t$ factors through passing to the completion $\hat{\sO}\cong k[[t]]$, so we may assume from the beginning that $\sO=k[[t]]$, and thus $A=k'[[s]]$. Moreover, we may insert the intermediate  complete dvr  $k'[[t]]$; by the multiplicativity of the trace form, this reduces us to the case $k'=k$, $n=e$. Writing $t=us^e$ for $u\in k[[s]]^\times$, we have $u=\sum_{i\ge0}u_is^i$ with $u_0\neq0$, so $t=u_0s^n\cdot\sum_{i\ge0}(u_i/u_0)s^i$. Replacing $s$ with $(\sum_{i\ge0}(u_i/u_0)s^i)^{1/e}\cdot s$, we may assume that $s^e=at$ with $a\in k^\times$. 

The trace form $q_{A_K/K}$ is easily computed: Let $H=\<1\>+\<-1\>$ denote the hyperbolic form. Then
\[
q_{A_K/K}=\begin{cases} \<e\>+\frac{e-1}{2}\cdot H&\text{ if $e$ is odd,}\\
 \<e\>+\<eat\>+\frac{e-2}{2}\cdot H&\text{ if $e$ is even.}
 \end{cases}
 \]
 so
 \[
\Delta_t(A/\sO)=\begin{cases} \<e\>-\<1\>+\frac{e-1}{2}\cdot H
&\text{ if $e$ is odd,}\\
 \<e\>-\<1\>+\<ea\>+\frac{e-2}{2}\cdot H&\text{ if $e$ is even.}
\end{cases}
\]
while the class $e_{0_A}(\Omega^1_{A/k}, dt)\in \GW(k)$ is
\[
e_{0_A}(\Omega^1_{A/k}, dt)=\begin{cases} \frac{e-1}{2}\cdot H
&\text{ if $e$ is odd,}\\
\<ea\>+\frac{e-2}{2}\cdot H&\text{ if $e$ is even.}
\end{cases}
\]

We observe that there is an ``error term'' $\<e\>-\<1\>$ (of rank zero!) in both cases:
\[
\Delta_t(A/\sO)=\<e\>-\<1\>+e_{0_A}(\Omega^1_{A/k}, dt)\in \GW(k).
\]
On the positive side, we do get a formula valid for both even and odd $e$, and although both 
$\Delta_t(A/\sO)$ and $e_{0_A}(\Omega^1_{A/k}, dt)$ depend {\it a priori} on $t$, their difference is independent of the choice of $t$. Going back to the general case $\sO\subset A$, we have
\[
\Delta_t(A/\sO)=\Tr_{k'/k}(\<e\>-\<1\>+e_{0_A}(\Omega^1_{A/k}, dt))\in \GW(k).
\]

Our next case is the specialization of a smooth projective hypersurface to the projective cone over a hyperplane section. Specifically, we take a dvr  $\sO$  with residue field $k$, quotient field $K$ and parameter $t$,  let $\sX\subset \P^{n+1}_{\sO}$ be the hypersurface defined by an equation of the form $F_t=0$, where
\[
F_t=F(X_0,\ldots, X_n)-t\cdot X_{n+1}^e
\]
and $F\in \sO[X_0,\ldots, X_n]$ is homogeneous of degree $e$, with $e$ prime to the characteristic if this is positive. We assume that the generic fiber $X_t\subset \P^{n+1}_K$ is smooth over $K$ and that $\bar{X}_0\subset \P^n_k$ defined by $F=0$ is also smooth (to simplify the notation we will write $F$ for the reduction of $F$ in $k[X_0,\ldots, X_n]$ and also for the image of $F$ in $K[X_0,\ldots, X_n]$). Clearly the special fiber $X_0$ is the projective cone $C_{\bar{X}_0}$ over $\bar{X}_0$. 

At this point we need to decide what motivic Euler characteristic we will be using for such singular varieties. After some consideration, we choose the ``Euler characteristic with compact support'': For $\pi_Y:Y\to \Spec k$ a separated finite type $k$-scheme, we have the object $\pi_{Y!}(1_Y)\in \SH(k)$. We give a quick review of the properties of these objects and their associated Euler characteristics; for details we refer the reader to the discussion in \cite[\S 1]{LevineEuler}.

If $k$ has characteristic zero, the fact that separated finite type $k$-schemes admit resolution of singularities implies that that 
$\pi_{Y!}(1_Y)$ is a strongly dualizable object of $\SH(k)$, so the Euler characteristic 
\[
\chi(\pi_{Y!}(1_Y)):=\Tr(\id_{\pi_{Y!}(1_Y)})\in  \End_{\SH(k)}(1_{\Spec k})\cong \GW(k)
\]
is well-defined. We denote $\chi(\pi_{Y!}(1_Y))$ by $\chi_c(Y/k)$. This Euler characteristic has the following properties:\\[5pt]
1. For $Y$ smooth and integral of dimension $d_Y$ over $k$, we have $\chi_c(Y/k)=\<-1\>^{d_Y}\cdot \chi(Y/k)$.\\[2pt]
2. For $Y$ smooth and proper over $k$, we have  $\chi_c(Y/k)= \chi(Y/k)$.\\[2pt]
3. Let $Y$ be a separated finite type $k$-scheme with closed subscheme $i:Z\subset Y$ having open complement $U\subset Y$. Then $\chi_c(Y/k)=\chi_c(U/k)+\chi_c(Z/k)$.\\[2pt]
4. Let 
\[
\xymatrix{
E\ar[r]\ar[d]&\tilde{Y}\ar[d]\\
F\ar[r]& Y
}
\]
be an abstract blow-up square of separated finite type $k$-schemes, that is, the diagram is cartesian, the maps $E\to \tilde{Y}$, $F\to Y$ are closed immersions, the map $\tilde{Y}\to Y$ is proper, and the induced map $\tilde{Y}\setminus E\to Y\setminus F$ is an isomorphism. Then
\[
\chi_c(Y/k)=\chi_c(\tilde{Y}/k)-\chi_c(E/k)+\chi_c(F/k).
\]

Indeed, (1) follows from the identities $\pi_{Y!}=\pi_\#\circ\Sigma^{-T_{Y/k}}$ (since $\pi_Y$ is smooth) and $\chi(Y\wedge(\P^1)^{\wedge m}/k)=\<-1\>^m\cdot\chi(Y/k)$. (2) follows from (1) and the fact that for $Y$ smooth and proper of odd dimension over $k$, $\chi(Y/k)$ is a multiple of the hyperbolic form $H$, and $\<-1\>\cdot [H]=[H]$. For (3), this is a consequence of the localization distinguished triangle for the closed immersion $i:Z\to Y$ and open complement $j:U\to Y$,
\[
j_!j^!\to \id_{\SH(Y)}\to i_*i^*
\]
applied to $1_Y$, together with the identities $j^!(1_Y)=1_U$, $i^*(1_Y)=1_Z$, $i_*=i_!$. One then applies $\pi_{Y!}$ and uses the additivity of trace in $\SH(k)$. (4) follows from (3). 

For our situation $\sX\to \Spec \sO$, we have $\chi_c(X_K/K)=\chi(X_K/K)$. For computing 
$\chi_c(X_0)=\chi_c(C_{\bar{X}_0})$, we have the blow-up square
\[
\xymatrix{
\bar{X}_0\ar[r]\ar[d]&\tilde{C}_{\bar{X}_0}\ar[d]\\
0\ar[r] &C_{\bar{X}_0}\ar@{=}[r]&X_0
}
\]
where $0=(0:\ldots:0:1)$ is the vertex of the cone $C_{\bar{X}_0}$, and $\tilde{C}_{\bar{X}_0}$ is the $\P^1$-bundle over $\bar{X}_0$, $\tilde{C}_{\bar{X}_0}=\Proj_{\bar{X}_0}(\sO_{\bar{X}_0}\oplus \sO_{\bar{X}_0}(1))$. This gives
\begin{align*}
\chi_c(X_0/k)&=\chi(\tilde{C}_{\bar{X}_0}/k)+\<1\>-\chi({\bar{X}_0}/k)\\
&=(\<1\>+\<-1\>)\cdot \chi({\bar{X}_0}/k)+\<1\>-\chi({\bar{X}_0}/k)\\
&=\<-1\>\cdot \chi({\bar{X}_0}/k)+\<1\>
\end{align*}

In the positive characteristic case $p>2$, Riou \cite[Appendix]{EllipticFlop} shows that $\pi_{Y!}(1_Y)$ is a strongly dualizable object of $\SH(k)[1/p]$. Moreover, the map $\GW(k)\to \GW(k)[1/p]$ is injective and the corresponding Euler characteristic $\chi_c(Y/k)\in \GW(k)[1/p]$ lands in $\GW(k)\subset \GW(k)[1/p]$ and has all the properties listed above. Thus, we can proceed just as in characteristic zero. A similar argument shows that for $k$ not perfect, we can safely pass to perfection and still have a well-defined Euler characteristic $\chi_c(Y/k)\in \GW(k)$.

We have the Jacobian rings ($F_i:=\del F/\del X_i$)
\begin{align*}
J_{F_t}&:=K[X_0,\ldots, X_{n+1}]/(F_0,\ldots, F_n, tX_{n+1}^{e-1})\\
J_F&:=k[X_0,\ldots, X_n]/(F_0,\ldots, F_n)
\end{align*}
and the version of $J_F$ over $K$,
\[
J_{F,K}:=K[X_0,\ldots, X_n]/(F_0,\ldots, F_n)
\]
As above, define
\[
\Delta_t(\sX/\sO)=\spc_t\chi_c(X_t/K)-\chi_c(X_0)\in \GW(k)
\]
We have the classical conductor identity
\[
\rnk \Delta_t(\sX/\sO)=(-1)^n\dim_kJ_F.
\]
In equal characteristic, this follows from the results of Milnor (characteristic 0) and Deligne (positive characteristic, noting that there are no Swan conductor terms in this simple case). In general, one can use the classical formulas for topological/\'etale Euler characteristics of the hypersurfaces $X_0$ and $X_K$ to establish this identity directly.

Suppose first that $n$ is odd. Then $\chi_c(X_t/K)=\chi(X_t/K)$ is hyperbolic, $\chi(X_t/K)=N_{n,e}\cdot H$, where $2N_{n,e}$ is the topological Euler characteristic. Also $\bar{X}_0\subset \P^n_k$ has even dimension $2m$ and $H^{m,m}(\bar{X}_0)_{prim}\cong J_{F, (m+1)\cdot e-n-1}$. Let $q_{J_F}$ be the quadratic form associated to the $k$-bilinear form $B_{\Jac, F}$ and let $q_{J_F,0}$ be the restriction of $q_{J_F}$ to the ``primitive Hodge part''  $\oplus_{q=0}^{2m}J_{F, (q+1)\cdot e-n-1}$. Then (Corollary~\ref{cor:HdgJacQ})
\[
\chi_c(\bar{X}_0/k)=\chi(\bar{X}_0/k)=\<e\>+\<-e\>\cdot [q_{J_{F,0}}]+m\cdot H.
\]
and 
\[
q_{J_F}=q_{J_F,0}+r\cdot H
\]
for some integer $r\ge0$. 
Thus
\begin{align*}
 \Delta_t(\sX/\sO)&=N_{n,e}\cdot H-(\<1\>+\<-1\>\cdot \chi_c(\bar{X}_0))\\
 &=N_{n,e}\cdot H-(\<1\>+\<-e\>+\<e\>\cdot [q_{J_{F,0}}]+m\cdot H)\\
 &=-\<1\>-\<-e\>-\<e\>\cdot [q_{J_F}]+(-1)^nM\cdot H+1\cdot H\\
 \end{align*}
with $M=1-N_{n,e}-m+r$.
Noting that 
\begin{align*}
-\<1\>-\<-e\>+1\cdot H&=-\<1\>-\<-e\>+\<e\>+\<-e\>\\
&=\<e\>-\<1\>
\end{align*}
we have
\[
 \Delta_t(\sX/\sO)=\<e\>-\<1\>+(-\<e\>)^n([q_{J_F}]+M\cdot H)
\]
Comparing ranks we have $M=0$ and thus
\begin{equation}\label{eqn:Odd}
\Delta_t(\sX/\sO)=\<e\>-\<1\>+(-\<e\>)^n[q_{J_F}]\in\GW(k)
\end{equation}

Now take $n=2m$ even. Then $\bar{X}_0$ has odd dimension, so (Corollary~\ref{cor:HdgJacQ} again)
\[
\chi_c(\bar{X}_0/k) = N_{e,n}\cdot H
\]
\[
\chi_c(X_t/K)=\chi(X_t/K)=\<e\>+\<-e\>\cdot q_{J_{F_t},0}+ m\cdot H
\]
where as above $q_{J_{F_t},0}$ is the restriction of $q_{J_{F_t}}$ to the primitive Hodge part of $J_{F_t}$ and
\[
q_{J_{F_t}}=q_{J_{F_t}, 0}+ r\cdot H
\]
for some $r\ge0$. Moreover, we have
\[
J_{F_t}=J_{F,K}\otimes_KK[X_{n+1}]/(X_{n+1}^{e-1})
\]
and 
\[
e_{F_t}=e_{F,K}\cdot (-et)X_{n+1}^{e-2}
\]
The graded ring $K[X_{n+1}]/(X_{n+1}^{e-1})$  with socle generator $(-et)X_{n+1}^{e-2}$ represents the quadratic form $q_e$ with 
\[
q_e=\begin{cases} \frac{e-1}{2}\cdot H&\text{ for $e$ odd,}\\
\<-et\>+ \frac{e-2}{2}\cdot H&\text{ for $e$ even.}
\end{cases}
\]
Thus there is an integer $L_{e,n}$ such that 
\[
q_{J_{F_t}}=\begin{cases} L_{e,n}\cdot H&\text{ for $e$ odd,}\\
\<-et\>\cdot q_{J_{F,K}}+L_{e,n}\cdot H&\text{ for $e$ even.}
\end{cases}
\]
and therefore there is an integer $R_{e,n}$ such that
\[
\chi(X_t/K)=\begin{cases} \<e\>+R_{e,n}\cdot H&\text{ for $e$ odd,}\\
\<e\>+\<t\>\cdot q_{J_{F,K}}+R_{e,n}\cdot H&\text{ for $e$ even.}
\end{cases}
\]

Note that $q_{J_F}=\spc_t q_{J_{F,K}}=\spc_t(\<t\>\cdot q_{J_{F,K}})$, and that for $e$ odd and $n$ even, the socle of $J_{F,K}$ is in odd degree $(e-2)(n+1)$, so $q_{J_{F,K}}$ is itself hyperbolic, say $q_{J_{F,K}}=s_{e,n}\cdot H$.  

Putting this together gives, for $n$  even and $e$ odd
\begin{align*}
\Delta_t(\sX/\sO)&=\<e\>+R_{e,n}\cdot H-(\<1\>+\<-1\>\cdot N_{e,n}\cdot H)\\
&=\<e\>-\<1\>+M_{e,n}\cdot H;\quad M_{e,n}=R_{e,n}-N_{e,n}.
\end{align*}
Comparing ranks again, we see that $M_{e,n}=s_{e,n}$ and
\[
\Delta_t(\sX/\sO)=\<e\>-\<1\>+(-\<e\>)^n\cdot[q_{J_F}].
\]
For $n$ even and $e$ even, we have
\begin{align*}
\Delta_t(\sX/\sO)&=\<e\>+[q_{J_F}]+R_{e,n}\cdot H-(\<1\>+\<-1\>\cdot N_{e,n}\cdot H)\\
&=\<e\>-\<1\>+ (-\<e\>)^n\cdot [q_{J_F}]+M_{e,n}H;\quad M_{e,n}=R_{e,n}-N_{e,n}.
\end{align*}
Comparing ranks once more, we have $M_{e,n}=0$ and
\[
\Delta_t(\sX/\sO)=\<e\>-\<1\>+(-\<e\>)^n\cdot[q_{J_F}].
\]

Finally, we note that the class $[q_{J_F}]\in \GW(k)$ is exactly the local Euler class with supports at $0\in \sX$ 
\[
e_0(\Omega_{\sX/k}, df)\in H^{n+1}_0(\sX, \sK^{MW}_{n+1}(\omega_{\sX/k}))
\cong \GW(k).
\]
This is a result of \cite{BachmannWickelgren} (see Theorem 2.18, Proposition 2.32 and Theorem 7.6), namely, that for a rank $m$ vector bundle $V$ on smooth $k$-scheme $Y$ of dimension $m$ and a section $s$  of $V$ with an isolated zero at some closed point $p\in Y$, the local Euler class 
\[
e_p(V, s)\in 
H^m_p(Y, \sK^{MW}(\det^{-1}(V)))\cong \GW(k(p ); \det^{-1}(V)\otimes\omega^{-1}_{Y/k}\otimes k(p ))
\]
is given as follows: suppose we have a local frame $(\lambda_1,\ldots, \lambda_m)$ for $V$ near $p$ and parameters $t_1,\ldots, t_m$ on $Y$ at $p$. Expressing $s=\sum_{i=1}^ms_i\lambda_i$ and using the $t_1,\ldots, t_m$ to define an isomorphism $\hat{\sO}_{Y,p}\cong k(p )[[t_1,\ldots, t_m]]$, we have the complete intersection ring $J(s):=k[[t_1,\ldots, t_m]]/(s_1,\ldots, s_m)$ and the Scheja-Storch element $e=\det(a_{ij})$, where $s_i=\sum_{j=1}^ma_{ij}t_j$. This gives us the quadratic form $q(s)$ on  $J(s)$, $q(s)(x)=\lambda\Leftrightarrow x^2=\lambda\cdot e$ and one has
\[
e_p(V, s)=[q(s)]\otimes (\lambda_m\wedge\ldots\wedge\lambda_1)^{-1}\otimes \del/\del t_1\wedge\ldots\wedge \del/\del t_m \in \GW(k(p ), \det^{-1}V\otimes \omega_{Y/k}^{-1}\otimes k(p )).
\]
In case  $Y=\sX$, $V=\Omega_{\sX/k}$, we have $\det^{-1}V\otimes \omega_{Y/k}^{-1} =\omega_{\sX/k}^{\otimes -2}=(\omega_{\sX/k}^{\otimes -1})^{\otimes 2}$ and as this is a square, we can safely ignore the twist. See the discussion of local indices in \cite[\S 4]{LevineEuler} for more details.

In summary, we have
\begin{theorem}\label{thm:ConductorFormula1} Let $k$ be a field of characteristic different from 2, let $\sO$ be a dvr    with quotient field $K$, residue field $k$ and parameter $t$, and let $f:\sX\to \Spec \sO$ be a flat projective regular $\sO$-scheme.  In addition, we suppose that   the generic fiber $X_t$ is smooth over $K$. Let $X_0$ denote the special fiber.  Finally, we suppose that either\\[5pt]
1. $\dim_\sO\sX=0$ and $\sX$ has  closed points $p_1,\ldots, p_s$ and $\sX\to \Spec \sO$  is tamely ramified. \\[2pt]
or\\[2pt]
2. $\sX\subset \P^{n+1}_\sO$ is a hypersurface defined by  a homogeneous equation of the form $F(X_0,\ldots, X_n)-tX_{n+1}^e$ with $e$ prime to the characteristic of $k$ if this is positive,  and the hypersurface $\bar{X}_0\subset\P^n_k$ defined by $F$ is smooth over $k$.\\[5pt]
Let $p_1,\ldots, p_s\in \sX$ be the closed points of $X_0$ in case (1) and let $p_1$ be the vertex of $X_0$ in case (2). Let $k_i$ be the residue field of $\sX$ at $p_i$ (so $k_1=k$ in case (2) and is a separable extension in case (1)). Let $e_i$ be the multiplicity of $X_0$ at $p_i$ and define 
\[
\Delta_t(\sX/\sO):=\spc_t\chi_c(X_t/K)-\chi_c(X_0/k).
\]
Then
 \[
 \Delta_t(\sX/\sO)= \sum_{i=1}^s\Tr_{k_i/k}[\<e_i\>-\<1\>+(-\<e_i\>)^{\dim_\sO\sX}\cdot e_{p_i}(\Omega_{\sX/k}, dt)]\in \GW(k) 
 \]
 \end{theorem}
 Note the two main differences from the classical case: the sign $(-1)^{\dim_\sO\sX}$ is replaced with the ``sign'' $(-\<e_i\>)^{\dim_\sO\sX}$ and one needs to introduce the error term $\<e_i\>-\<1\>$.  We also note that in case $k=\R$, we have $\<e_i\>=\<1\>$ since $e_i>0$, so these subtleties disappear.

The case of a quasi-homogeneous singularity defined by $F\in k[Y_0,\ldots, Y_n]$, such that the corresponding hypersurface $Y_F\subset \P_k(a_*)$ is smooth, is essentially the same. Suppose that $F$ has weighted degree $\dg$ and assume that $\dg\cdot\prod_ia_i$ is prime to the characteristic if this is positive. We consider the family $\sX\subset \P_\sO((a_*,1))$ defined by $F_t=F-tY_{n+1}^\dg$ in $\P_\sO((a_*,1))$. For $t=0$, the affine open subset $U_{n+1}\subset \P_k((a_*,1))$ is just $\A^{n+1}$ and the hypersurface $X_0\subset \P_k((a_*,1))$ has intersection with $U_{n+1}$ equal to the affine hypersurface defined by $F$. $\sX$ is smooth over $\sO$ away from $0\in U_{n+2}(k)$ and $X_0\setminus\{0\}$ is a (stratified) $\A^1$-bundle over $(F=0)\subset \P_k(a_*)$. The only significant changes from the homogeneous case are seen in the difference between Corollary~\ref{cor:HdgJacQ} and Theorem~\ref{thm:HdgJacQH}, giving the following result

\begin{theorem}\label{thm:ConductorFormula2}  Let $k$ be a field of characteristic different from 2, let $\sO$ be a dvr   with quotient field $K$, residue field $k$ and parameter $t$, and let $f:\sX\subset\P_\sO(a_*,1)$ be defined by $F_t=F(Y_0,\ldots, Y_n)-tY_{n+1}^\dg$, where $F$ is weighted homogeneous of degree $\dg$ and weights $a_*$ and $\{F=0\}\subset \P_k(a_*)$ is a smooth hypersurface. Let $X_0$ denote the special fiber and let $p$ denote the vertex of the ``cone'' $X_0$. Suppose $\lcm(a_*)$ divides $\dg$ and that $\dg\cdot\prod_ia_i$ is prime to the characteristic of $k$ if this is positive.   Let
\[
\Delta_t(\sX/\sO):=\spc_t\chi_c(X_t/K)-\chi_c(X_0/k).
\]
Then
 \[
 \Delta_t(\sX/\sO)= \<\dg\cdot\prod_ia_i\>-\<1\>+(-\<\dg\>)^n\cdot e_p(\Omega_{\sX/k}, dt)\in \GW(k) 
 \]
 \end{theorem}

 Based on these result, we make the following conjecture.\footnote{The formula for $\Delta_t(\sX/\sO)$ in Conjecture~\ref{conj:Main}, under the  assumption that $\sO$ has characteristic zero and that the singularity of $X_0$ at $p_i$ can be resolved by a single weighted blow-up with smooth   exceptional divisor, has been verified in a   paper by Ran Azouri \cite{Azouri}.}

 \begin{conjecture}\label{conj:Main} Let $\sO$ be a dvr  with  residue field $k$ of characteristic different from two, quotient field $K$ and parameter $t$. Let $f:\sX\to \Spec \sO$ be a flat and proper morphism with $\sX$ smooth over $k$ and with generic fiber $X_K$ smooth over $K$. Suppose that the special fiber $X_0$ has only finitely many singular points $p_1,\ldots, p_s$. Suppose in addition that each   $(X_0, p_i)$  is (over the algebraic closure $\overline{k(p_i)}$) a quasi-homogeneous hypersurface singularity  of degree $\dg_i$, with weights $a_*^{(i)}$, and with $\lcm(a_*^{(i)})=1$. Let $a^{(i)}=\prod_ja_j^{(i)}$. Finally, suppose that  $\dg_i$, $a^{(i)}$ and $[k(p_i):k]$ are prime to $\Char k$ if this is positive. Define
 \[
\Delta_t(\sX/\sO):=\spc_t(\chi_c(X_K/K))-\chi_c(X_0/k)\in \GW(k).
\]
Then
\[
\Delta_t(\sX/\sO)=\sum_{i=1}^s \Tr_{k(p_i)/k}[\<\dg_i\cdot a^{(i)}\>-\<1\>+(-\<\dg_i\>)^{\dim_\sO\sX}e_{p_i}(\Omega_{\sX/k}, dt)].
\]
 \end{conjecture}
 
 \begin{remark}  We are not able to formulate a reasonable general conjecture without the assumption of quasi-homogeneity.
 \end{remark}

For $k=\C$, the conjecture reduces to the classical case, that is,  Milnor's theorem computing the Milnor number as the dimension of the Jacobian ring. 
For $k=\R$, as the pair $(\rnk, \text{sig}):\GW(\R)\to \Z\times\Z$ is injective, the conjecture follows from the case $k=\C$ together with results of Eisenbud-Levine \cite{EL} and Khimshiashvili \cite{Kh}. Indeed, let $\sX\to \Spec \sO$ be a proper morphism of $\R$-schemes with  $\sO$ the henselization of $\R[t]_{(t)}$, satisfying the hypotheses of the conjecture. Using $t$ as parameter, it makes sense to speak of the fiber $X_t$ over $t\in \R$ for $|t|<<1$.   Khimshiashvili  computes the Milnor number of an isolated hypersurface singularity $f=0$ at $0\in \R^n$, for $f:\R^n\to \R$ a real analytic map, in terms of the topological degree of the normalized  gradient
\[
\frac{\nabla f}{||\nabla f||}: S_\epsilon^{n-1}\to S^{n-1},\ 0<\epsilon<<1. 
\]
Combining this with  theorem of Eisenbud-Levine and Khimshiashvili computes the difference $\chi^\topo(X_t(\R))-\chi^\topo(X_0(\R))$ (for $0<t<<1$) by applying the signature to the expression in our conjecture, noting that  $\chi^\topo(X_t(\R))=\text{sig}(\chi_c(X_t/\R))$ and $\chi^\topo(X_0(\R))=\text{sig}(\chi_c(X_0/\R))$. Since we are orienting $\R$ by our choice of parameter $t$, the specialization $\spc_t(\chi_c(X_K/K))\in \GW(\R)$ has signature equal to $\lim_{t\to 0^+}\chi^\topo(X_t(\R))$ (the term $\chi^\topo(X_t(\R))$ is constant in $t$ for $0<t<<1$): this follows by diagonalizing the form $\chi_c(X_K/K)$ over $K=\R[t]_{(t)}^h[1/t]$ and then noting that for $ut^n\in \R[t]_{(t)}^h[1/t]$, $u=u(t)$ a unit in  $\R[t]_{(t)}^h$, $n\in \Z$, we have
$\<ut^n\>=\<u(0)\>\in \GW(\R)$ for  $0<t<<1$ and $\spc_t\<ut^n\>=\<u(0)\>$ as well.

In both cases $k=\C, \R$, the conjecture holds without the assumption that the singularities are quasi-homogeneous, and just reduces to the results of Milnor, Eisenbud-Levine and Khimshiashvili as explained  above. 

Of course, it would be nice to have an intrinsic explanation of the error term $\<\dg\prod_ia_i\>-\<1\>$ and the reason why  the factor $(-\<\dg\>)^{\dim_\sO\sX}$ replaces the classical factor $(-1)^{\dim_\sO\sX}$, as well as a formulation valid for arbitrary isolated singularities. As a first step in this direction, we reformulate the problem in terms of Ayoub's motivic nearby cycles and vanishing cycles functors.  

\section{Motivic nearby cycles}
\label{sec:nearby-cycles}

Let us recall some facts about the motivic nearby cycle functors introduced by Ayoub in \cite[Chapter 3]{ayoub-thesis-2}. 

As in the rest of the paper, we write $\mathcal{O}$ for a dvr with fraction field $K$ and residue field $k$ with uniformizer $t$. We assume that $K$ and $k$ are of characteristic $0$. Let $f:\mathcal{X}\to \Spec(\mathcal{O})$ be a finite type separated morphism\footnote{The theory as developed in \cite{ayoub-thesis-2} is restricted to quasi-projective morphisms. This restriction is imposed because the six operation formalism in \cite{ayoub-thesis-1} is also restricted to quasi-projective morphisms; however, standard arguments based on Zariski descent and Chow's lemma can be used to lift these restrictions.}, with generic fiber $f_{K}:X_{K}\to \Spec(K)$. The motivic nearby cycle functor is a functor
\[
\Psi^{t}_{f}:\SH(X_{K})\to \SH(X_{k}).
\]
This is often simply denoted by $\Psi_{f}$; but the functor really depends on the uniformizer $t$, and this dependency is important in our case, so we keep $t$ in the notation. More precisely, the functor $\Psi^{t}_{f}$ defined in \cite[Chapter 3]{ayoub-thesis-2} is the \emph{tame} motivic nearby cycle functor (i.e. the motivic analogue of the tame nearby cycle in the \'etale theory): since we are working in equal characteristic $0$, we omit the adjective ``tame'' everywhere.

The special case where $f=\id_{\Spec(\mathcal{O})}$ is particularly important:
\[
\Psi^{t}:=\Psi^t_{\id}:\SH(K)\to \SH(k).
\]
Let $g:\mathcal{Y}\to \mathcal{X}$ be another finite type separated morphism. Then we have a canonical comparison natural transformation
\[
\alpha_{g}:g_{k}^{*}\Psi^t_{f}\to \Psi^t_{f\circ g}g_{K}^{*}
\]
which is a natural isomorphism if $g$ is smooth. We also have a natural transformation obtained from $\alpha_{g}$ by passing to right adjoints
\[
\beta_{g}:\Psi^t_{f}(g_{K})_{*}\to (g_{k})_{*}\Psi^t_{f\circ g}
\]
which is a natural isomorphism if $g$ is proper.

In particular, we have 
\[
\alpha_{f}:f_{k}^{*}\Psi^t\to \Psi^t_{f}f_{K}^{*}
\]
which is an isomorphism for $f$ smooth and
\[
\beta_{f}:\Psi^t(f_{K})_{*}\to (f_{k})_{*}\Psi^t_{f}
\]
which is an isomorphism for $f$ proper.

All the properties of $\Psi^t$ so far are part of the fact that $\Psi^t$ is a ``specialization system'' in the sense of \cite[Definition 3.1.1]{ayoub-thesis-2}; see \cite[Definition 3.5.6]{ayoub-thesis-2}. From this, one can deduce further compatibilities with purity isomorphisms and exceptional operations, as explained in \cite[\S 3.1.2]{ayoub-thesis-2}. Let us recall some of those results.

The functor $\Psi^t$ commutes with Thom equivalences: for any vector bundle $\mathcal{V}$ on $\mathcal{X}$, there are canonical isomorphisms \cite[Proposition 3.1.7]{ayoub-thesis-2}:
\[
\Th^{\pm 1}(\mathcal{V}_{k})\circ\Psi^t_{f}\simeq \Psi^t_{f}\Th^{\pm 1}(\mathcal{V}_{K}).
\]
Let $g:\mathcal{Y}\to \mathcal{X}$ be a finite type separated morphism. Then we have canonical natural transformations
\[
\mu_{g}:(g_{k})_{!}\Psi^t_{f\circ g}\to \Psi^t_{f}(g_{K})_{!}
\]
and
\[
\nu_{g}:\Psi^t_{f\circ g}(g_{K})^{!}\to (g_{k})^{!}\Psi^t_{f}.
\]
If $g$ is proper (resp. smooth), then $\mu_{g}$ is invertible, and via the natural isomorphism $g_{!}\simeq g_{*}$ identifies with $\beta_{g}^{-1}$ (resp. $\nu_{g}$ is invertible, and via purity isomorphisms identifies with $\alpha_{g}^{-1}$ up to Thom isomorphisms).

The properties so far actually hold without any assumptions on the residual characteristics. However, we need to be in equal characteristic $0$ for constructibility and the interaction of the (tame) motivic nearby cycle with the symmetric monoidal structure of $\SH$, to which we now turn.

The functor $\Psi^{t}_{f}$ sends constructible objects of $\SH(X_{K})$ to constructible objects of $\SH(X_{k})$ by \cite[Theorem 3.5.14]{ayoub-thesis-2}. Motivic nearby cycle form a lax-monoidal functor (as an application of \cite[Proposition 3.2.12]{ayoub-thesis-2}): we have a natural transformation
\[
\Psi^t_{f}(-)\otimes \Psi^t_{f}(-)\to \Psi^t_{f}(-\otimes -).
\]
This lax-monoidal structure commute with exterior products in the following sense: given another morphism $f':\mathcal{Y}\to \Spec(\mathcal{O})$, the morphism
\[
\Psi^t_{f}(-)\boxtimes \Psi^t_{f'}(-)\to \Psi^t_{ (f\times_{\mathcal{O}}f')}(-)
\]
defined using $\alpha_{f},\alpha_{f'}$ and the lax-monoidal structure is an isomorphism \cite[Theorem 3.5.17]{ayoub-thesis-2}.

In particular, $\Psi^t$ is a symmetric monoidal functor: we have $\Psi^t(\one_{K})=\one_{k}$ and $\Psi^t(M\otimes N)=\Psi^t(M)\otimes \Psi^t(N)$ for all $M,N\in \SH(K)$. This implies that $\Psi^{t}$ sends dualizable objects to dualizable objects (with our assumptions on the characteristics, this also follows from the preservation of constructibility).

Moreover, $\Psi^t$ commutes with Verdier duality in $\SH$ \cite[Theorem 3.5.20]{ayoub-thesis-2}: we have a canonical isomorphism
\[
\Psi^t_{f}\underline{\Hom}_{\SH(K)}(-,f_{K}^{!}\one)\simeq \underline{\Hom}_{\SH(k)}(\Psi^t_{f},f_{k}^{!}\one).
\]

\section{Motivic vanishing cycles}
\label{sec:vanish-cycles}

Write $j:X_{K}\to \mathcal{X}$ (respectively $i:X_{k}\to \mathcal{X}$) for the open immersion of the generic fiber (resp. the closed immersion of the special fiber) of $f$ in the total space.

To define motivic vanishing cycles, we have to look closer at the construction of $\Psi^{t}_{f}$ (for varying $f$) in terms of specialization systems. With the notations of \cite[Definition 3.2.3, Definition 3.5.3]{ayoub-thesis-2}, we have $\Psi^{t}=\mathcal{R}\bullet \chi$ where $\mathcal{R}$ is a certain diagram of schemes indexed by $\Delta\times\mathbb{N}^{\times}$ and $\chi$ is the ``canonical specialization system'' $i^{*}j_{*}$ (in particular, for the rest of this paragraph, $\chi$ does not mean Euler characteristic!) The diagram of schemes $(\mathcal{R},\Delta\times\mathbb{N}^{\times})$ satisfies $\mathcal{R}([0],1)=\mathbb{G}_{m}$ and this provides a morphism of diagrams of schemes $(\mathbb{G}_{m},e)\to (\mathcal{R},\Delta\times\mathbb{N}^{\times})$. By \cite[Definition 3.2.13]{ayoub-thesis-2}, this induces a morphism of specialization systems $\chi\simeq (\mathbb{G}_{m},e)\bullet \chi\to \mathcal{R}\bullet \chi=\Psi^{t}$  which, applied to $\mathcal{X}$
yield a natural transformation
\[
i^{*}j_{*}\to \Psi^{t}_{f}.
\]
Combined with the unit of the adjunction $j^{*}\dashv j_{*}$, this provides a natural transformation
\[
\SP_{f}:i^{*}\to \Psi^{t}_{f}j^{*}.
\]
We use the same convention as for $\Psi$ and just write $\SP$ for $\SP_{\id}$.

This transformation $\SP_{f}$ is usually called ``specialization'' (because it induces the classical specialization morphism on cohomology), hence the notation. The word ``specialization'' is already overloaded is this paper, and we will avoid this terminology.

The natural transformation $\SP_{f}$ comes from a morphism of specialization systems, and this implies that it is compatible with the natural transformations $\alpha,\beta,\nu,\mu$. We spell out this compatibility for the more involved cases of $\mu$ and $\nu$ which we need later.

Let $g:\mathcal{Y}\to \mathcal{X}$ be a finite type separated morphism. Then for any $M\in \SH(\mathcal{Y})$, we have a commutative diagram
\begin{equation}\label{SP-mu}
\begin{tikzcd}[column sep=large]
  (g_{k})_{!}i^{*}M \arrow[r,"\SP_{f\circ g}(M)"] \arrow[d,"\sim"] & (g_{k})_{!} \Psi^{t}_{f\circ g}j^{*}M \arrow[r, "\mu_{g}"] & \Psi^{t}_{f}(g_{K})_{!} j^{*}M \arrow[d,"\sim"] \\
  i^{*}g_{!}M \arrow[rr,"\SP_{f}(g_{!}M)"] & & \Psi^{t}_{f}j^{*}g_{!} M
\end{tikzcd}
\end{equation}
where the vertical arrows are base change isomorphisms. Similarly, we have for any $N\in \SH(\mathcal{X})$ a commutative diagram
\begin{equation}\label{SP-nu}
\begin{tikzcd}[column sep=large]
i^{*}g^{!} \arrow[r,"\SP_{f\circ g}(g^{!}N)"] \arrow[d,"\mathrm{Ex}^{!,*}"] & \Psi^{t}_{f\circ g}j^{*}g^{!}N \arrow[r, "\mathrm{Ex}^{!,*}", "\sim"'] & \Psi^{t}_{f\circ g}(g_{K})^{!} j^{*}N \arrow[d,"\nu_{g}"] \\
  (g_{k})^{!}i^{*}N \arrow[rr,"\SP_{f}(N)"] & & (g_{k})^{!}\Psi^{t}_{f}j^{*}N
\end{tikzcd}
\end{equation}
where $\mathrm{Ex}^{!,*}$ denotes exchange morphisms (and the vertical one is not necessarily an isomorphism).

\begin{definition}
The motivic vanishing cycle $\Phi^{t}_{f}$ is defined to be the cofiber of $\SP_{f}$: for $M\in \SH(\mathcal{X})$, we put
  \[
i^{*}M\stackrel{\SP_{f}}{\to}\Psi^{t}_{f}j^{*}M\to\Phi^{t}_{f}M\stackrel{+}{\to}.
\]
\end{definition}  
Given the lack of functoriality of cofibers in triangulated categories, this definition does not specify a vanishing cycle functor but rather the isomorphism class of individual motivic vanishing cycle objects. Since we are primarily interested in their Euler characteristic, this is not a problem for us. Nevertheless this lack of functoriality can be remedied, either by doing a more refined construction with diagrams of schemes as in \cite{ayoub-thesis-2} or by translating the construction of loc.cit. in the language of $\infty$-categories.

From the classical analytic and \'etale theory, we know that vanishing cycles of the unit measure singularities of $f$. In particular, we expect them to disappear when $f$ is smooth. This is indeed the case.

\begin{proposition}\label{prop:smooth-vanish}
Assume that $f$ is a smooth morphism. Then $\Phi^{t}_{f}\one=0$ and $\Psi^{t}_{f}\one\simeq \one$.
\end{proposition}
\begin{proof}
The natural transformation $\SP_{f}$ is constructed from a morphism of specialization systems, and this implies that it is compatible with the natural transformations $\alpha_{f}$, which are isomorphisms because $f$ is smooth. This shows that $\Phi^{t}_{f}\one\simeq f^{*}\Phi^{t}\one$. But $\Psi^{t}$ is monoidal and $\SP(\one):\one\simeq i^{*}\one\to \Psi^{t}\one\simeq \one$ is the identity, hence $\Phi^{t}\one=0$.
\end{proof}

% In this situation, we can also describe the functoriality of $\Psi^{t}$ more precisely, and this is information that we will use later to identify specialization morphisms. Assume still that $f$ is a smooth morphism. We can pass to left adjoints in the natural isomorphism $\alpha_{f}^{-1}$ to define a natural transformation
% \[
% \gamma_{f}:(f_{k})_{\sharp}\Psi^{t}_{f}\to \Psi^{t}(f_{K})_{\sharp}.
% \]
% When evaluated at the unit and combined with the previous proposition, this gives a morphism
% \[
% \gamma_{f}(\one):(f_{k})_{\sharp}\one\to \Psi^{t}(f_{K})_{\sharp}(\one).
% \]

\begin{lemma}\label{lem:SP-iso}
  Assume that $f$ is smooth and proper, or more generally that $f$ admits a smooth compactification $f=\bar{f}\kappa$ over $\Spec(\mathcal{O})$ with $\kappa:\mathcal{X}\to\overline{\mathcal{X}}$ an open immersion and $\bar{f}:\overline{\mathcal{X}}\to\Spec(\mathcal{O})$ a smooth proper morphism. Assume as well that the reduced complement $\overline{\mathcal{X}}\setminus \mathcal{X}$ is smooth over $\mathcal{O}$. Then
  \[
\SP(f_{!}f^{!}\one):i^{*}f_{!}f^{!}\one \to \Psi^{t}j^{*}f_{!}f^{!}\one
\]
is an isomorphism.
\end{lemma}  
\begin{proof}
Let $f=\bar{f}\kappa$ be the smooth compactification of $f$. Write $\iota:\mathcal{Z}\to \overline{\mathcal{X}}$ for the complementary reduced closed immersion and $g=\bar{f}\iota$ for the structure morphism of the $\mathcal{O}$-scheme $\mathcal{Z}$. We have a localization triangle
  \[
f_{!}f^{!}\one \to \bar{f}_{!}\bar{f}^{!}\one\to g_{!}\iota^{*}\bar{f}^{!}\one\stackrel{+}\to.
  \]
Note that $\iota^{*}\bar{f}^{!}$ and $g^{!}$ differ only by a Thom isomorphism by purity since $\bar{f}$ and $g$ are both smooth. So by naturality of $\SP$, we only need to prove that both $\SP_{\id}(\bar{f}_{!}\bar{f}^{!}\one)$ and $\SP_{\id}(g_{!}g^{!}\one)$ are isomorphisms, i.e., we can assume that $f$ is smooth and proper.
  
We start with the diagram \eqref{SP-mu} (for $f=\mathrm{id},g=f$ and $M=f^{!}\one$)
  \[
\begin{tikzcd}[column sep=large]
  (f_{k})_{!}i^{*}f^!\one \arrow[r,"\SP_{f}(f^!\one)"] \arrow[d,"\sim"] & (f_{k})_{!} \Psi^{t}_{f}j^{*}f^!\one \arrow[r, "\mu_{f}"] & \Psi^{t}(f_{K})_{!} j^{*}f^!\one \arrow[d,"\sim"] \\
  i^{*}f_{!}f^!\one \arrow[rr,"\SP(f_{!}f^!\one)"] & & \Psi^{t}j^{*} f^!\one.
\end{tikzcd}
  \]
Since $f$ is proper, then $\mu_{f}$ is an isomorphism. So we are reduced to proving that
  \[
\SP_{f}(f^!\one):i^{*}f^!\one \to \Psi^{t}_{f}j^{*}f^{!}\one
\]
is an isomorphism. We also have the diagram \eqref{SP-nu} (for $f=\mathrm{id},g=f$ and $N=\one$)
  \[
\begin{tikzcd}[column sep=large]
i^{*}f^{!}\one \arrow[r,"\SP_{f}(f^{!}\one)"] \arrow[d,"\mathrm{Ex}^{!,*}"] & \Psi^{t}_{f}j^{*}f^{!}\one \arrow[r, "\mathrm{Ex}^{!,*}", "\sim"'] & \Psi^{t}_{f}(f_{K})^{!} j^{*}\one \arrow[d,"\nu_{f}"] \\
  (f_{k})^{!}i^{*}\one \arrow[rr,"\SP(\one)"] & & (f_{k})^{!}\Psi^{t}j^{*}\one.
\end{tikzcd}
\]
Because $f$ is smooth, the morphism $\nu_{f}$ is an isomorphism, and so is the vertical exchange transformation $\mathrm{Ex}^{!,*}$. So we just need to know that $\SP(\one)$ is an isomorphism, which is a special case of Proposition \ref{prop:smooth-vanish}.
\end{proof}

\section{The quadratic Euler characteristic of vanishing cycles}
\label{sec:euler-char-motiv}

Once again, we assume that we are in equal characteristic $0$. In that case, $\Psi^t$ is monoidal, so there is an induced morphism $(\Psi^t)_{*}:\GW(K)\simeq \End_{\SH(K)}(\one)\to \End_{\SH(k)}(\one)\simeq \GW(k)$, and since categorical traces commute with monoidal functors, the following holds. For a field $F$,  let $\chi(-/F)$ denote the categorical Euler characteristic, defined on strongly dualizable objects of $\SH(F)$. 

\begin{lemma}
  Let $M\in \SH(K)$ be dualizable. Then $\Psi^t(M)\in \SH(k)$ is dualizable and we have
\[
\chi(\Psi^t(M)/k)=(\Psi^t)_{*}(\chi(M)/K)
\]
\end{lemma}  

Our first goal is to compute $(\Psi^t)_{*}$ explicitly. 

Because $\Psi^t$ is a monoidal functor and commutes with Thom equivalences, it induces a morphism of graded algebras
\[
(\Psi^t)_{*}:\Hom_{\SH(K)}(\one,\Sigma^{*,*}\one)\to \Hom_{\SH(k)}(\one,\Sigma^{*,*}\one)
\]
whose $0$-th graded part is the morphism $(\Psi^t)_{*}$ from the previous paragraph. Let $F$ be a field and $K^{MW}_{*}(F)$ denote the graded Milnor-Witt K-theory algebra of $F$. Recall from \cite[Theorem 1.23]{Morel-book} or \cite[Theorem 6.4.1]{MorelICTP} that if $F$ is perfect, we have an isomorphism $\End_{\SH(F)}(\one)\simeq K^{MW}_{0}(F)$, and that more generally, we have an isomorphism of graded algebras
\[
\Hom_{\SH(F)}(\one,\Sigma^{*,*}\one)\simeq K^{MW}_{*}(F).
\]
Via the isomorphism $\phi_{0,F}:\GW(F)\to K^{MW}_0(F)$ mentioned in Remark~\ref{rem:specialization}, we have the isomorphism $\GW(F)\cong 
\End_{\SH(F)}(\one)$

\begin{proposition}\label{prop:psi-spe}
The induced morphism
  \[
(\Psi^t)_{*}:K^{MW}_{*}(K)\to K^{MW}_{*}(k) 
\]
is the specialization morphism $s^{t}_{v}$.

In particular, the  morphism
\[
(\Psi^t)_{*}:\GW(K)\to \GW(k)
\]
induced by Morel's isomorphisms $\phi_{0,K}$, $\phi_{0,k}$
is the specialization morphism $\spe_{t}$.
\end{proposition}

\begin{proof}
The last statement of the proposition follows from the first part, together with the definition of 
$\spe_{t}$ described in Remark~\ref{rem:specialization}. To prove the first claim, since $(\Psi^{t})_{*}$ is a homomorphism of graded rings, it suffices to check that $(\Psi^{t})_{*}(\eta)=\eta$, $(\Psi^{t})_{*}([u])=[\bar{u}]$ for $u\in \mathcal{O}^{\times}$ and $(\Psi^{t})_{*}([t])=1$.

Since we are in equal characteristic, $\sO$ contains a subfield $k_0$ with $k_0\to k$ a finite extension. Let $p:\Spec \sO\to \Spec k_0$ be the corresponding morphism. 
By  \cite[Lemma 3.5.10]{ayoub-thesis-2}, the natural transformation $\SP\circ p^*:i^*p^*\to \Psi\circ p^*$ is an isomorphism. Applying this to $\Sigma^{a,b}\one_{k_0}$, we find that $\SP_{\Sigma^{a,b}\one_\sO}:i^*(\Sigma^{a,b}\one_\sO)\to \Psi(\Sigma^{a,b}\one_\sO)$ is an isomorphism. Given a map $\alpha:\one_\sO\to \Sigma^{a,b}\one_\sO$, the naturality of  $\SP$ yields the commutative diagram
\[
\xymatrixcolsep{50pt}
\xymatrix{
\one_k\ar[d]_{i^*\alpha}\ar[r]_-\sim^-{\SP_{\one_\sO}}&\Psi(\one_\sO)\ar[d]^{\Psi(\alpha)}\\
\Sigma^{a,b}\one_k\ar[r]_-\sim^-{\SP_{\Sigma^{a,b}\one_\sO}}&\Psi(\Sigma^{a,b}\one_\sO
}
\]
Applying this to $\alpha=\eta_\sO:\one_\sO\to \Sigma^{-1,-1}\one_\sO$ and $\alpha=[u]:\one_\sO\to \Sigma^{1,1}\one_\sO$ gives the result we want for $\eta$ and $[u]$.

Finally, we compute $(\Psi_{t})_{*}[t]$. The idea is to use a ``partial N\'eron model'' $\mathcal{G}$ of $\mathbb{G}_{m}$, which is the smooth $\mathcal{O}$-group scheme obtained by gluing two copies of $\mathbb{G}_{m,\mathcal{O}}$ (which we suggestively denote by $\mathbb{G}_{m,\mathcal{O}}$ and $t\cdot \mathbb{G}_{m,\mathcal{O}}$) along their generic fibers, with the transition isomorphism given by multiplication by $t$. To sum up, we have $\mathcal{G}=\mathbb{G}_{m,\mathcal{O}}\cup_{\mathbb{G}_{m,K}}t\cdot \mathbb{G}_{m,\mathcal{O}}$, with generic fiber $\mathbb{G}_{m,K}$ and special fiber $\mathbb{G}_{m,k}\coprod t\cdot\mathbb{G}_{m,k}$. By construction, the point $t\in \mathbb{G}_{m}(K)\simeq \mathcal{G}(K)$ extends to a section $\tilde{t}\in \mathcal{G}(\mathcal{O})$ whose reduction $\tilde{t}_{k}$ is $1\in (t\cdot \mathbb{G}_{m,k})(k)$. We can alternatively glue two copies of $\mathbb{P}^{1}$ in the same way, and get a compactification $\overline{\mathcal{G}}$ of $\mathcal{G}$ which has generic fiber $\mathbb{P}^{1}_{K}$ and special fiber $\mathbb{P}^{1}_{k}\cup t\cdot\mathbb{P}^{1}_{m,k}$ glued along one point. We write $f$ (resp. $\bar{f}$) for the structure morphism of $\mathcal{G}$ (resp. $\overline{\mathcal{G}}$).
 
By naturality of $\SP$, we have a commutative diagram
\[
\begin{tikzcd}
  \one \arrow[d,"\tilde{t}_{k}"] & i^{*}\one \arrow[r, "\SP(\one)", "\sim"'] \arrow[l, "\sim"] \arrow[d,"\tilde{t}_{\mathcal{O}}"] & \Psi^{t}j^{*} \one \arrow[d,"\tilde{t}_{\mathcal{O}}"] \arrow[r, "\sim"] &  \Psi^{t}\one \arrow[d,"\tilde{t}_{K}"]  \\
  (f_{k})_{!}(f_{k})^{!}\one  & i^{*}f_{!}f^{!}\one \arrow[r, "\SP(f_{!}f^{!}\one)"] \arrow[l, "\sim"] & \Psi^{t}j^{*} f_{!}f^{!}\one \arrow[r, "\sim"]  &  \Psi^{t}(f_{K})_{!}(f_{K})^{!}\one 
\end{tikzcd}  
\]

The morphism $f$ does not satisfy the assumptions of Lemma \ref{lem:SP-iso}; and indeed, the bottom middle morphism is not an isomorphism in this case. We have $\mathcal{G}_{k}\simeq \mathbb{G}_{m}^{\coprod 2}$, so
\[
i^{*}f_{!}f^{!}\one\simeq (f_{k})_{!}(f_{k})^{!}\one\simeq ((r_{k})_{!}(r_{k})^{!}\one)^{\oplus 2}
\]
while we have 
\[
\Psi^{t}(f_{K})_{!}(f_{K})^{!}\one\simeq (r_{k})_{!}(r_{k})^{!}\one
\]
by Lemma \ref{lem:SP-iso}. We claim that the induced map $((r_{k})_{!}(r_{k})^{!}\one)^{\oplus 2}
\to (r_{k})_{!}(r_{k})^{!}\one$ is simply the fold map $\id+\id$. This can be checked separately on each factor, and then it follows from the computation of motivic nearby cycles on a branch of the semi-stable morphism with reduced special fiber $\bar{f}$ (case $m=1$ in \cite[Theorem 3.3.43]{ayoub-thesis-2}). This fact, together with the commutative diagram above, implies that $(\Psi_{t})_{*}[t]=1$ and finishes the proof.
\end{proof}

\begin{proposition}
Let $f:\mathcal{X}\to \mathrm{Spec}(\mathcal{O})$ be a finite type separated morphism.
  Then 
\[
\Delta_{t}(\mathcal{X}/\mathcal{O})=\chi(\Phi^t f_{!}\one/k).
\]
If $f$ is moreover proper, we have
\[
\Delta_{t}(\mathcal{X}/\mathcal{O})=\chi((f_{k})_{!}\Phi^t_{f}\one/k).
\]
Finally, if $f$ is proper, $X_{K}$ is smooth over $K$, the field $k$ is perfect and $f$ only has isolated singularities $x_{1},\ldots x_{s}\in X_{k}$, then $\Phi^t_{f}$ is supported on $x_1,\ldots,x_{s}$ and we have
\[
\Delta_{t}(\mathcal{X}/\mathcal{O})=\sum_{i=1}^{s}\Tr_{k(x_{i})/k}\chi(\Phi^t_{ f}(\one)_{x_{i}}/k(x_i)).
\]
\end{proposition}
\begin{proof}
We start with the triangle defining motivic vanishing cycles, applied to the object $f_{!}\one$, and use base change:
  \[
(f_{k})_{!}\one\to\Psi^{t}(f_{K})_{!}\one\to\Phi^{t}f_{!}\one\stackrel{+}{\to}.
\]
Every term is constructible in $\SH(k)$ hence dualizable, and by additivity of Euler characteristics we find
\[
\chi(\Phi^{t}f_{!}\one/k)=\chi(\Psi^{t}(f_{K})_{!}\one/K)-\chi_{c}(X_{k}/k).
\]
By Proposition \ref{prop:psi-spe}, we have $\chi(\Psi^{t}(f_{K})_{!}\one/k)=\spe_{t}\chi_{c}(X_{K}/K)$ and we deduce
the first formula.

Let us now assume $f$ proper. The fact that $\SP$ is induced from a morphism of specialization systems implies that there is a commutative diagram of natural transformations
\[
\begin{tikzcd}
  i^{*}(f_{K})_{*}  \arrow[r,"\SP\circ (f_{K})_{*}"] \arrow[d] & \Psi^{t}(f_{K})_{*} \arrow[d,"\beta_{f}"] \\
  (f_{k})_{*}i^{*} \arrow[r,"(f_{k})_{*}\SP"] & (f_{k})_{*}\Psi^{t}_{f}
\end{tikzcd}
\]
where the vertical arrows are isomorphisms because $f$ is proper. So we get an induced isomorphism $\Phi^t f_{!}\one\simeq (f_{k})_{!}\Phi^t_{f}\one$ of cofibers. This shows the second formula.

Finally, we assume in addition that $f$ has smooth generic fiber and isolated singularities in the special fiber. Let $\iota:\mathcal{X}^{\circ}\to \mathcal{X}$ be the open immersion of the smooth locus of $f$ into the total space. Then $\iota^{*}\Phi^{t}_{f}\one$ is the cofiber of the morphism
\[
\one\simeq \iota^{*}i^{*}\one\stackrel{\iota^{*}\SP}{\to}\iota^{*}\Phi^{t}_{f}\one
\]
Since $\iota$ is smooth, we have an isomorphism $\alpha_{\iota}:\iota^{*}\Phi^{t}_{f}\one\simeq \Phi^{t}_{f\circ \iota}\one$, and because $\SP$ comes from a morphism of specialization systems, the composite morphism
\[
\one \to \Phi_{f\circ\iota}^{t}\one 
\]
is the morphism $\SP$ corresponding to $f\circ\iota$. By Proposition \ref{prop:smooth-vanish}, this morphism is an isomorphism, and we deduce that $\iota^{*}\Phi^{t}_{f}\one=0$. By localization, this precisely means that $\Phi^{t}_{f}$ is supported on the finitely many points $x_{1},\ldots,x_{s}$. It remains to identify the action of $f_{*}$ with the trace. Since $k$ is perfect, the extensions $k(x_i)/k$ are all separable. In this case, the identification of $f_*$ with the trace map on $\GW(-)$ is proven in \cite[Corollary 8.5]{BachmannWickelgren}.
\end{proof}

\appendix

\section{A coordinate-free approach to Theorem~\ref{thm:Main}} Here we describe another approach for proving Theorem~\ref{thm:Main}, suggested by the third author, that avoids most of the explicit computations used by Carlson-Griffiths, and is essentially coordinate-free.

 \subsection{A dual resolution}
Let $V$ be an $n+2$-dimensional $k$-vector space with symmetric algebra $S(V)$ over $k$, and let $i: X\hookrightarrow \P= \Proj_k(S(V))\cong \P^{n+1}_k$ be a smooth $n$-dimensional hypersurface of degree  $\dg$ in $\P$, defined by a homogeneous prime ideal $(F)$ for some $F\in S^\dg(V)$.

For $Y$ a smooth $k$-scheme, we let $\Omega^p_Y=\Omega^p_{Y/k}$ denote the sheaf of $p$-forms, and 
$\sT^p_Y=\sH{om}_{\sO_Y}(\Omega^p_Y,\sO_Y)$ the sheaf of tangent $p$-vectors, 
so that $\sT_Y=\sT^1_Y$ is the tangent sheaf of $Y$. 

The canonical pairing
\[
\<-,-\>:\sT_Y\times \Omega^1_Y\to \sO_Y
\]
extends to the interior multiplication of tangent vectors on $p+1$-forms 
\[
\sT_Y\times \Omega^{p+1}_Y\to \Omega^p_Y.
\]
Explicitly, for each local section $v$ of $\sT_Y$ on some open $j:U\to Y$ the interior multiplication operator 
\[
\iota_v:j^*\Omega^{p+1}_Y\to j^*\Omega^p_Y
\]
is defined for local sections $\eta_0,\ldots, \eta_p$ of $j^*\Omega^1_Y$ by
\[
\iota_v(\eta_0\wedge\ldots\wedge\eta_p):=\sum_{j=0}^p(-1)^j\<v,\eta_j\>\cdot \eta_0\wedge\ldots\widehat{\eta}_j\wedge\ldots\wedge \eta_p.
\]
This extends to the operator $\iota_v:j^*\Omega^{p+r}_Y\to j^*\Omega^p_Y$ for $v$ a local section of $\ext^r\sT_Y$, with 
\[
\iota_{v_1\wedge\ldots\wedge v_r}(\eta):=i_{v_r}\circ\ldots\circ i_{v_1}(\eta)
\]
for local sections $v_1,\ldots, v_r$ of $\sT_Y$ and $\eta$ of $\Omega_Y^{p+r}$. Sending a local section $v$ of $\ext^r\sT_Y$ to $\iota_v:j^*\Omega^r_Y\to j^*\sO_Y$ defines an isomorphism 
\begin{equation}\label{eqn:ExteriorIdent}
\ext^r\sT_Y\cong \sT_Y^r; 
\end{equation}
we will identify $\ext^r\sT_Y$ and $\sT_Y^r$ via this isomorphism.

Dually, we have the interior multiplication $\iota_\eta:j^*\sT^{p+1}_Y \to j^*\sT^p_Y$ for $\eta$ a local section of $\Omega^1_Y$, with
\[
\iota_\eta(v_0\wedge\ldots\wedge v_p)=\sum_{j=0}^p (-1)^j \<v_j,\eta\>\cdot v_0\wedge\ldots\wedge \widehat{v}_j\wedge\ldots\wedge v_p
\]
and this extends to the interior multiplication $\iota_\eta:j^*\sT^{p+r}_Y \to j^*\sT^p_Y$ for $\eta$ a local section of $\Omega^r_Y$, with
\[
\iota_{\eta_1\wedge\ldots\eta_r}(v)= \iota_{\eta_r}\circ\ldots\circ\iota_{\eta_1}(v)
\]
for local sections $\eta_1,\ldots, \eta_r$ of $\Omega_Y$. For $p=0$,  this agrees with the evaluation pairing $\Omega^r_Y\otimes \sT^r_Y\to \sO_Y$.

Letting $\sN_X$ be the normal sheaf $\sH{om}_{\sO_X}(\sI_X/\sI_X^2, \sO_X)$,  we have the canonical identifications  $\sI_X/\sI_X^2\cong \sO_X(-X)$,  $ \sO_X(X)\cong \sN_X$. We have the 
exact sequence 
\[
0\to \sI_X/\sI^2_X\xrightarrow{d} i^*\Omega_\P\to \Omega_X\to0
\]
whose $\sO_X$-dual is the  tangent-normal sheaf sequence for $i:X\to \P$. This gives us the exact sequence
\begin{equation}\label{basic}
0\to \sT_X\to i^*\sT_{\P}\xrightarrow{d^*} \sO_X(X)\to 0.
\end{equation}
Explicitly, for $v$ a local section of $ i^*\sT_{\P}$, $d^*(v)=i^*(\<v,dF/F\>)$, where $F$ is a local generating section for  for $\sI_X$. 

Apply $\ext^j$ and piecing together the exact sequences
\[
0\to \ext^j\sT_X\to i^*\ext^j\sT_{\P}\xrightarrow{d^*} \ext^{j-1}\sT_X\otimes \sO_X(X)\to 0
\]
gives rise to an exact sequence
\[0\to \sT^p_X\to i^*\sT^p_{\P}\xrightarrow{\iota_{dF/F}} i^*\sT^{p-1}_{\P}(X)\xrightarrow{\iota_{dF/F}}\cdots\xrightarrow{\iota_{dF/F}}
i^*\sT_{\P}((p-1)X)\xrightarrow{\iota_{dF/F}}  \sO_X(pX)\to 0,
\]
for each $0\leq p\leq n$ (for $p=0$ this is the trivial ``sequence'' 
$0\to \sO_X\to \sO_X\to 0$).   We view this as a resolution of  
$\sT^p_X[p]$ by the complex of sheaves
\[\sD(p):=  \left(0\to i^*\sT^p_{\P}\to 
i^*\sT^{p-1}_{\P}(X) \to\cdots\to i^*\sT_{\P}((p-1)X)\to \sO_X(pX)\to 
0\right).\]
Shifting by $[-p]$ gives the resolution of $\sT^p_X$, 
\[
\sT^p_X\xrightarrow{\epsilon_p} \sD(p)[-p].
\]

Via the identifications $\sT^p_X=\ext^p\sT_X$, we have  the exterior product  pairings
\[
\sT^p_X\tensor_{\sO_X}\sT^q_X\to \sT^{p+q}_X,\ i^*\sT^a_\P\tensor_{\sO_X}i^*\sT^b_\P\to i^*\sT^{a+b}_\P.
\]
Since ${\sD}(p)$ is viewed as a complex concentrated in degrees $[-p,0]$,  these pairings
satisfy the required Leibniz rule  for defining a natural pairing of complexes 
\begin{equation}\label{eqn:TruncKosProd}
\tilde\mu_{p,q}:\sD(p)\tensor \sD(q)\to \sD(p+q).
\end{equation}
We have the canonical isomorphism
\[
\sD(p)[-p]\tensor \sD(q)[-q]\to (\sD(p)\tensor \sD(q))[-p-q]
\]
which, for local sections $a$ of $\sD(p)[-p]^i=\sD(p)^{i-p}$, and $b$ of $\sD(q)[-q]^j=\sD(q)^{j-q}$ sends $a\otimes b$ to $(-1)^{iq}\cdot a\otimes b\in \sD(p)^{i-p}\otimes \sD(q)^{j-q}\subset 
(\sD(p)\tensor \sD(q))^{i+j-p-q}=(\sD(p)\tensor \sD(q))[-p-q]^{i+j}$. This defines the multiplication
\[
\mu_{p,q}:\sD(p)[-p]\tensor \sD(q)[-p]\to \sD(p+q)[-p-q]
\]
as the composition
\[
\sD(p)[-p]\tensor \sD(q)[-p]\cong (\sD(p)\tensor \sD(q))[-p-q]\xrightarrow{\tilde{\mu}_{p,q}[-p-q]}
 \sD(p+q)[-p-q].
 \]
 
The augmentations $\epsilon_p:\sT_X^p\to  \sD(p)[-p]$, $\epsilon_q:\sT_X^q\to  \sD(q)[-q]$
$\epsilon_{p+q}:\sT_X^{p+q}\to  \sD(p+q)[-p-q]$ fit together to give a commutative diagram of multiplication maps
\[
\xymatrix{
\sT^p_X\otimes \sT_X^q\ar[r]^{-\wedge-}\ar[d]^{\epsilon_p\otimes\epsilon_q}&\sT^{p+q}_X\ar[d]^{\epsilon_{p+q}}\\
\sD(p)[-p]\otimes \sD(q)[-q]\ar[r]^{\mu_{p,q}}&\sD(p+q)[-p-q]
}
\]

Twisting the surjection in \eqref{basic} to give the surjection
\begin{equation}\label{basic2}
d^*:i^*\sT_{\P}(-X)\to\sO_X
\end{equation}
yields an exact Koszul complex $\Kos(d^*)$ of locally free sheaves on $X$,
\[0\to \ext^{n+1}i^*\sT_{\P}(-(n+1)X)\to\cdots\to i^*\sT_{\P}(-X)\to\sO_X\to 
0,\]
with $\sO_X$ in degree 0. The complex $\Kos(d^*)$ comes equipped with an algebra structure, and is self-dual under 
\[
\sH{om}_{\sO_X}(-,\ext^{n+1}i^*\sT_{\P}(-(n+1)X)[n+1]). 
\]
Our identifications  \eqref{eqn:ExteriorIdent} for $0\le i\le p$ defines an isomorphism  
\begin{equation}\label{eqn:KoszulIdent}
\sD(p)\cong \tau_{\ge -p}\Kos(d^*) \otimes\sO_X(pX), 
\end{equation}
and
the product map \eqref{eqn:TruncKosProd} is similarly induced from
algebra structure on the Koszul complex.

\subsection{Identifying the Hodge pairing}
For $Y$ a smooth $k$ scheme of dimension $N$ over $k$,  the interior product
defines the isomorphism
\begin{equation}\label{eqn:IsoY}
\iota_Y(p):\sT^p_Y\tensor_{\sO_Y}\omega_Y\xrightarrow{\sim}\Omega^{N-p}_Y
\end{equation}

For our hypersurface $X$, we have in addition  the residue isomorphism $\res_X:i^*\omega_\P(X)\xrightarrow{\sim} \omega_X$, with inverse $dF/F\wedge-: \omega_X\to i^*\omega_\P(X)$ sending a local section $\eta$ of $\omega_X$ to $dF/F\wedge \tilde{\eta}$, where $\tilde\eta$ is a local lifting of $\eta$ to a section of $i^*\Omega^n_\P$. Combined with $\iota_\P(p-j)$ this gives us the isomorphism
\[
\iota'_\P(p-j):=\iota_\P(p-j)\circ dF/F\wedge-:i^*\sT^{p-j}_\P\otimes_{\sO_X}\omega_X(jX)\xrightarrow{\sim}i^*\Omega^{q+j+1}_\P((j+1)X).
\]
We use this to reinterpret the complexes $\sC(q)^*$ of Lemma 3.1 as follows.

\begin{lemma} \label{lem:ComparisonIsoComplex}  The maps     
\[
(-1)^{p+\frac{j(j+1)}{2}}\iota'_\P(p-j):\sD(p)\otimes \omega_X[-p]^j\to \sC(q)^j
\]
for $j=0,\ldots, p-1$, the map
\[
(-1)^{p(p-1)/2}\cdot\id:\sD(p)\otimes \omega_X[-p]^p=\omega_X(pX)\to \sC(q)^p=\omega_X(pX)
\]
and the map
\[
\iota_X(p):\sT^p_X\otimes\omega_X\to \Omega_X^q
\]
define an isomorphism of augmented complexes
\[
\iota(p):[\sT^p_X\otimes\omega_X\to\sD(p)\otimes \omega_X[-p]]\to [ \Omega_X^q\to\sC(q)].
\]
\end{lemma}

\begin{proof} Take a point $x\in X$ and a regular system of parameters  $t_0,...,t_{n+1}$  in the local ring $\sO_{\P, x}$  with $t_0$ a defining equation for $X$. One checks by a direct computation that for $j=0,\ldots, p-1$, the diagram
\[
\xymatrixcolsep{40pt}
\xymatrix{
i^*\sT^{p-j}_\P\otimes_{\sO_X}\omega_X(jX)\ar[d]^{\iota_\P'(p-j)}\ar[r]^-{(-1)^p\iota_{dF/F}}&i^*\sT^{p-j-1}_\P\otimes_{\sO_X}\omega_X((j+1)X)\ar[d]^{\iota_\P'(p-j-1)}\\
i^*\Omega^{q+j+1}_\P((j+1)X)\ar[r]^-{dF/F\wedge-}&i^*\Omega^{q+j+2}_\P((j+2)X)
}
\]
commutes up to the factor $(-1)^{j+1}$, and the diagrams
\[
\xymatrix{
\sT_X^p\otimes\omega_X\ar[r]\ar[d]^{\iota_X(p)}&i^*\sT^p_\P\otimes\omega_X\ar[d]^{\iota_\P'(p)}\\
\Omega_X^q\ar[r]^-{dF/F\wedge-}&i^*\Omega^{q+1}_\P(X)
}
\]
and 
\[
\xymatrixcolsep{40pt}
\xymatrix{
i^*\sT_\P\otimes_{\sO_X}\omega_X((p-1)X)\ar[d]^{\iota_\P'(1)}\ar[r]^-{(-1)^p\iota_{dF/F}}&\omega_X(pX)\ar@{=}[d]\\
i^*\Omega^{n}_\P(pX)\ar[r]^-{\pi}&\omega_X(pX)
}
\]
commute up to the factor $(-1)^p$; here $\pi$ is the restriction map $\pi_p$ considered in Lemma~\ref{lem:ExactSeq}. The result follows from this.
\end{proof}

\begin{lemma}\label{lem:ProductCommutes} For $p+q=n$, define the pairing
\[
-\wedge-:\sT_X^p\otimes\omega_X\times \sT_X^q\otimes\omega_X\to \sT_X^n\otimes\omega_X\otimes\omega_X
\]
by sending $(\alpha\otimes\eta, \beta\otimes\eta')$ to $\alpha\wedge\beta\otimes(\eta\otimes\eta')$. Then the diagram
\[
\xymatrix{
\sT_X^p\otimes\omega_X\times \sT_X^q\otimes\omega_X\ar[r]^-{-\wedge-} \ar[d]^{\iota_X(p)\otimes\iota_X(q)}&\sT_X^n\otimes\omega_X\otimes\omega_X\ar[d]^{\iota_X(n)\otimes\id}\\
\Omega_X^q\times\Omega_X^p\ar[r]^-{-\wedge-}&\omega_X
}
\]
commutes.
\end{lemma}
\begin{proof} This follows by an easy computation after choosing local coordinates near a point of $X$ as in the proof of the previous lemma.
\end{proof}

We have the natural map
\begin{multline*}
\coker\left(H^0(X,i^*\sT_{\P}((p-1)X)\tensor_{\sO_X}\omega_X)
\to H^0(X,\sO_X(pX)\tensor_{\sO_X}\omega_X)\right)\\
\to{\mathbb H}^p(X,\sD(p)[-p]\tensor_{\sO_X}\omega_X).
\end{multline*}
The augmentation $\sT^p_X\otimes\omega_X\to \sD(p)[-p]\tensor_{\sO_X}\omega_X$ gives us the natural isomorphism
\[
{\mathbb H}^p(X,\sD(p)[-p]\tensor_{\sO_X}\omega_X)
\cong H^p(X,\sT^p_X\tensor_{\sO_X}\omega_X)
\]
 for each $0\leq p\leq n$; putting these together gives us the natural map
 \begin{multline}
\coker\left(H^0(X,i^*\sT_{\P}((p-1)X)\tensor_{\sO_X}\omega_X)
\to H^0(X,\sO_X(pX)\tensor_{\sO_X}\omega_X)\right)\\
\xrightarrow{\Phi_p}H^p(X,\sT^p_X\tensor_{\sO_X}\omega_X)
\end{multline}

It follows from Lemma~\ref{lem:ExactSeq} and Lemma~\ref{lem:ComparisonIsoComplex} that 
$\Phi_p$ is an isomorphism for $p\neq n/2$. For $n=2m$ even, $\Phi_m$ is injective and the composition of $\Phi_m$ with the isomorphism 
\[
H^m(\iota_X(m)):H^m(X,\sT^m_X\tensor_{\sO_X}\omega_X)\to
H^m(X, \Omega^m_X)
\]
has image $H^m(X, \Omega^m_X)_{prim}$.

We note that there is a {\em canonical} isomorphism 
\begin{equation}\label{eqn:Canon}
\can:\omega_{\P}\xrightarrow{\sim} \ext^{n+2}V\tensor_k\sO_{\P}(-n-2),
\end{equation}
coming from the canonical resolution of the diagonal on $\P\times \P$. This is 
in fact an isomorphism of ${\rm PGL}_{n+2}(k)$-equivariant sheaves. If we give $V$ a basis $x_0,\ldots, x_{n+1}$, this identifies $\P$ with $\P^{n+1}=\Proj\, k[X_0,\ldots, X_{n+1}]$ and  $V$ with $H^0(\P^{n+1}, \sO_{\P^{n+1}}(1))$ by sending $x_i$ to $X_i$.  Then the global section of 
$\omega_{\P}(n+2)$ corresponding via $\can$ to $x_0\wedge\ldots x_{n+1}$  is the $n+1$-form $\Omega$  \eqref{eqn:Omega}.

The choice of generator $F$ for $H^0(\P, \sI_X(\dg))$ gives us the isomorphisms
\[
\can_F:\sO_\P(nX)\cong \sO_\P(n\dg), \ \can_F:\sO_X(nX)\cong \sO_X(n\dg)
\]
and together with $\can$ the isomorphism
\[
 H^0(X, \omega_X(pX))\cong  H^0(X,\sO_X((p+1)\dg -n-2))\otimes\ext^{n+2}V
 \]

The pairing $V^\vee\times V\to k$ extends to the $k$-linear map
\[
\del:V^\vee\to \Der_k(S(V), S(V))_{-1}
\]
of $V^\vee$ to the degree -1 derivations of $S(V)$. Evaluation at $F\in S^\dg(V)$ defines the map $\ev_F:V^\vee\to S^{\dg-1}(V)$. We let  $W(X,V)\subset S^{\dg-1}(V)$ be the image of $V^\vee$ under $\ev_F$, giving the isomorphism
\begin{equation}\label{eqn:evF}
\ev_F:V^\vee\xrightarrow{\sim} W(X,V).
\end{equation}
We define $J(X,V)$ to be the quotient ring $S(V)/(W(X,V))$. 

A basis $x_0,\ldots, x_{n+1}$ for $V$ gives us the isomorphism $S(V)\cong k[X_0,\ldots, X_{n+1}]$, with  $x_i\in V$ mapping to $X_i$, and the map $V^\vee\to \Der_k(S(V), S(V))_{-1}$ sends $x^i$ in the dual basis $x^0,\ldots, x^{n+1}$ to $\del/\del X_i$. This gives us the corresponding isomorphism of $J(X, V)$ with $J(F)$.

By  Proposition~\ref{prop: Main}  the resulting surjection
\[
S^{(p+1)\dg -n-2}(V)\otimes\ext^{n+2}V\to  H^0(X, \omega_X(pX))
\]
defines an isomorphism
\begin{multline*}
J(X,V)_{(p+1)\dg -n-2}\otimes\ext^{n+2}V\\\xrightarrow{\xi_p} \coker\left(H^0(X,i^*\sT_{\P}((p-1)X)\tensor_{\sO_X}\omega_X)
\xrightarrow{\iota_{dF/F}\otimes\id} H^0(X,\omega_X(pX))\right)
\end{multline*}
We note that $\xi_p$ depends on the choice of defining equation $F$: changing $F$ to $\lambda\cdot F$ changes $\xi_p$ to $\lambda^{-p-1}\cdot \xi_p$.

Choosing a basis for $V$ as above, for $A\in k[X_0,\ldots, X_{n+1}]_{(p+1)\dg -n-2}$ lifting $\bar{A}\in J(F)_{(p+1)\dg -n-2}$,  we have the section
\[
\frac{A}{F^p}\otimes\frac{\Omega}{F}\in H^0(X, \sO_X(pX)(\dg-n-2)\otimes i^*\omega_\P(X)(n+2-\dg)),
\]
and $\xi_p(\bar{A})$ is the image in $\coker \iota_{dF/F}\otimes\id$ of $(A/F^p)\res(\Omega/F)\in H^0(X, \omega_X(pX))$.

Let 
\[
\phi_p:J(X,V)_{(p+1)\dg -n-2}\otimes\ext^{n+2}V\to H^p(X,\sT^p_X\tensor_{\sO_X}\omega_X)
\]
be the composition $\Phi_p\circ\xi_p$. 

We have the isomorphisms
\begin{multline*}
J(X,V)_{(n+2)(\dg -2)}\otimes(\ext^{n+1}V)^{\otimes 2}\\\xrightarrow{\xi'_n} \coker\left(H^0(X,i^*\sT_{\P}((p-1)X)\tensor_{\sO_X}\omega_X)
\to H^0(X, \omega_X^{\otimes 2}(nX))\right)
\end{multline*}
and
\begin{multline}
\coker\left(H^0(X,i^*\sT_{\P}((n-1)X)\tensor_{\sO_X}\omega_X^{\otimes 2})
\to H^0(X, \omega_X^{\otimes 2}(nX))\right)\\
\xrightarrow{\Phi'_n} H^n(X,\sT^n_X\tensor_{\sO_X}\omega_X^{\otimes 2})
\end{multline}
arising similarly from the resolution $\sT^n_X\tensor_{\sO_X}\omega_X^{\otimes 2}\to \sD(n)[-n]\tensor_{\sO_X}\omega_X^{\otimes 2}$. Explicitly,  if we choose a basis for $V$ as above, then for $A\in k[X_0,\ldots, X_{n+1}]_{(n+2)(\dg-2)}$ lifting $\bar{A}\in J(F)_{(n+2)(\dg-2)}$, 
$\xi'_n(\bar{A})$ is the image in the cokernel of the element
\[
\frac{A}{F^n}\otimes\res(\frac{\Omega}{F})\otimes\res(\frac{\Omega}{F}).
\]

Let 
\[
\phi'_n:J(X,V)_{(n+2)(\dg -2)}\otimes(\ext^{n+1}V)^{\otimes 2}
\to H^n(X,\sT^n_X\tensor_{\sO_X}\omega_X^{\otimes 2})
\]
be the composition $\Phi'_n\circ\xi'_n$. 

\begin{lemma}\label{lem:Products}  1. We have
\[
H^p(\iota_X(p))\circ\phi_p=(-1)^{(p-1)p/2}\psi_p.
\]
2. Take $p+q=n$. Then for elements $A\in J(X,V)_{(q+1)\dg -n-2}\otimes\ext^{n+2}V$ and  $B\in J(X,V)_{(p+1)\dg -n-2}\otimes\ext^{n+2}V$, we have
\[
\psi_q(A)\wedge\psi_p(B)=(-1)^{(n-1)n/2}\cdot H^n(\iota_X(n)\otimes\id)(\phi'_n(AB))
\]
in $H^n(X, \omega_X)$.
\end{lemma}

\begin{proof} (1) follows from the fact that in degree $p$,  the comparison isomorphism
 $\iota(p):\sD(p)[-p]\otimes \omega_X\to \sC(q)$ is $(-1)^{(p-1)p/2}\id:\omega_X(pX)\to 
 \omega_X(pX)$. 
 
 For (2), we have
 \[
 \phi_q(A)\wedge  \phi_p(B)\in H^n(X, \sT^n_X\otimes\omega_X^{\otimes 2}).
 \]
The isomorphism $\sD(p)[-p]^p\otimes\sD(q)[-q]^q\to (\sD(p)^0\otimes\sD(q)^0)[-n]^n$ used in defining the product $\sD(p)[-p]^p\otimes\sD(q)[-q]^q\to\sD(n)[-n]$, sends $a\otimes b$ to $(-1)^{pq}a\otimes b$. Thus,  we have
\begin{align*}
\phi_q(A)\wedge  \phi_p(B)&=\Phi_q(\xi_q(A))\wedge\Phi_p(\xi_p(B))\\
&=\Phi'_n((-1)^{pq}(\xi_q(A)\wedge \xi_p(B))\\
&=(-1)^{pq}\phi'_n(AB)
\end{align*}
On the other hand, by  Lemma~\ref{lem:ProductCommutes} and (1), we have
\begin{align*}
 H^n(\iota_X(n)\otimes\id)(\phi_q(A)\wedge \phi_p(B))&=
(-1)^{(p-1)p/2}\cdot (-1)^{(q-1)q/2}\cdot \psi_q(A)\wedge\psi_p(B)
\end{align*} 
By this and  Lemma~\ref{lem:ProductCommutes}, we have
 \[
H^n(\iota_X(n)\otimes\id)(\phi'_n(AB))=
 (-1)^{(p-1)p/2+(q-1)q/2+pq}\psi_q(A)\wedge\psi_p(B)
 \]
 With $q=n-p$, we have
 \[
\frac{1}{2}[(p-1)p+(n-p-1)(n-p)+2p(n-p)]=
\frac{1}{2}(n-1)n,
\]
completing the proof.
\end{proof}

\subsection{Scheja-Storch class and the trace map}

We let  $W:=W(X,F)\subset S^{\dg-1}(V)$ be the $k$-subspace spanned by the partial derivatives  $v(F)$, $v\in V^\vee$, so the evaluation map 
\[
\ev_F:V^\vee\to S(V),\quad \ev_F(v)=v(F),
\]
 induces an isomorphism $\ev_F:V^\vee\to W$. 

Using the complex $\sD(n)[-n]\tensor\omega_X^{\tensor 2}$, the twist by $(\ext^{n+2}V)^{\otimes 2}$ of the 1-dimensional socle 
$J(X,V)_{(n+2)(e-2)}$ is
identified with
$H^n(X,\omega_X)$ via the map $H^n(\iota_X(n)\otimes\id)\circ\phi'_n$; we denote this isomorphism by
\[
\phi_X:J(X,V)_{(n+2)(e-2)}\otimes(\ext^{n+2}V)^{\otimes 2}\xrightarrow{\sim} H^n(X,\omega_X).
\]
To complete our discussion, we need to compute  the composition 
\[
Tr_X\circ 
\phi_X:J(X,V)_{(n+2)(e-2)}\otimes(\ext^{n+2}V)^{\otimes 2}\to k. 
\]

We first discuss a more general situation.

Suppose we are given a subspace $W\subset H^0(\P,\sO_{\P}(e-1))$ with $\dim_k 
W=n+2$, such that the linear system $|W|$ has no base points; equivalently, the 
ideal in the polynomial algebra $S(V)$ generated by $W\subset S^{e-1}(V)$ has 
radical equal to $S(V)_+$, the graded maximal ideal. We have an induced 
surjection $W\tensor_k\sO_{\P}\to \sO_{\P}(e-1)$ which we may equivalently 
regard as a surjection
\begin{equation}\label{eqW}
 W\tensor_k\sO_{\P}(-(e-1))\to \sO_{\P}.
 \end{equation}

The quotient $R(W)=S(V)/W\cdot S(V)$ is an Artinian graded complete intersection 
ring, whose 1-dimensional socle is in degree $(n+2)(e-2)$; this follows from the 
computation, for such graded complete intersections, of their Hilbert series.  
 
Note that $V\subset S^1(V)$ generates the homogeneous maximal ideal $S(V)_+$, 
and we have that $W\subset H^0(\P,\sO_{\P}(e-1))=S^{e-1}(V)$. Hence we may find 
a linear map (not unique, of course) 
\[
\varphi:W\to V\tensor_kS^{e-2}(V)
\]
so that on composition with the multiplication map 
\[
V\tensor_kS^{e-2}(V)\to S^{e-1}(V)
\]
we obtain the  inclusion $W\subset S^{e-1}(V)$. We may equivalently view $\varphi$ as 
 giving a linear inclusion 
\[
\varphi:W\to H^0(\P,V\tensor_k\sO_{\P}(e-2)),
\]
or an $\sO_\P$ linear map of sheaves $W\otimes \sO_\P\to V\tensor_k\sO_{\P}(e-2)$, which after twisting is also the same as a map
\[
\tilde{\varphi}:W\tensor_k\sO_{\P}(-(e-1))\to V\tensor_k\sO_{\P}(-1).
\]

Choose a basis $x_0,\ldots, x_{n+1}$ for $V$, giving as above the identification $S(V)\cong k[X_0,\ldots, X_n]$ with $x_i$ mapping to $X_i$. If $W$ is the subspace $W(X,V)$ of $S^{\dg-1}(V)$ with basis the partial derivatives $\del F/\del X_i$,  we may view $\varphi$ as an 
$(n+2)\times(n+2)$ matrix $(a_{ij})\in M_{n-2, n-2}(k[X_0,\ldots, X_n]_{\dg-2})$, with
\[
\frac{\del F}{\del X_j}=\sum_{i=0}^{n+1}a_{ij}X_i,
\]
just as used in defining the Scheja-Storch element for $F$. 

The canonical surjective map 
\[
V\tensor_k\sO_{\P}(-1)\to \sO_{\P}
\]
corresponding to the identification  $V=H^0(\P,\sO_{\P}(1))$, yields 
an exact Koszul complex of locally free sheaves on $\P$:
\begin{multline}\label{canonical}
0\to  \ext^{n+2}V\tensor_k\sO_{\P}(-n-2)\to 
\ext^{n+1}V\tensor_k\sO_{\P}(-n-1)\to\\\cdots\to V\tensor_k\sO_{\P}(-1)\to 
\sO_{\P}\to 0.
\end{multline}
The surjective map \eqref{eqW} gives rise to a similar exact Koszul complex of locally free 
sheaves on $\P$:
\begin{multline}\label{koszulW} 
0\to  \ext^{n+2}W\tensor_k\sO_{\P}(-(n+2)(e-1))\to 
\ext^{n+1}W\tensor_k\sO_{\P}(-(n+1)(e-1))\to\\
\cdots\to 
W\tensor_k\sO_{\P}(-(e-1))\to \sO_{\P}\to 0.
\end{multline}

The composition 
\[
W\tensor_k\sO_{\P}(-(e-1))\stackrel{\tilde{\varphi}}{\longrightarrow}
V\tensor_k\sO_{\P}(-1)\to \sO_{\P}\]
is the surjection \eqref{eqW}. Hence the map $\tilde{\varphi}$
induces a map between the two Koszul complexes \eqref{koszulW} and 
\eqref{canonical},
\begin{multline*}
\scalebox{0.95}{$
\xymatrixcolsep{10pt}
\xymatrix{ 
0\ar[r] & \ext^{n+2}V\tensor_k\sO_{\P}(-n-2)\ar[r] 
&\ext^{n+1}V\tensor_k\sO_{\P}(-n-1)\ar[r] &\cdots\\
0\ar[r] & \ext^{n+2}W\tensor_k\sO_{\P}(-(n+2)(e-1))\ar[r]\ar[u]_{\varphi_{n+2}} 
& \ext^{n+1}W\tensor_k\sO_{\P}(-(n+1)(e-1))\ar[r]\ar[u]_{\varphi_{n+1}} &\cdots }
$}\\
\\
\scalebox{0.95}{$
\xymatrixcolsep{10pt}
\xymatrix{
\cdots\ar[r] & \ext^2V\tensor_k\sO_{\P}(-2)\ar[r] &
V\tensor_k\sO_{\P}(-1)\ar[r] &\sO_{\P}\ar[r] & 0\\
\cdots\ar[r]&\ext^2W\tensor_k\sO_{\P}(-2(e-1))\ar[r]\ar[u]_{\varphi_2} & W\tensor_k\sO_{\P}(-(e-1))\ar[r]\ar[u]_{\varphi_1}& 
\sO_{\P}\ar[r]\ar[u]_{=}& 0.} $}
\end{multline*}
where $\varphi_j=\ext^j\tilde{\varphi}$.  

Now apply the exact functor $\sH{om}(-,\omega_{\P})$ to the above exact diagram, 
to obtain a similar commutative diagram with exact rows:
\begin{multline}\label{diagram}
\scalebox{0.95}{$
\xymatrixcolsep{10pt}
\xymatrix{ 0\ar[r] & \omega_{\P}\ar[r]\ar@{=}[d]
&V^\vee\tensor_k\omega_{\P}(1)\ar[r]
\ar[d]^{\varphi_1^*} &\ext^2V^\vee\tensor_k\omega_{\P}(2)\ar[r]
\ar[d]^{\varphi_2^*} & \cdots \\
0\ar[r] & \omega_{\P}\ar[r] & W^\vee\tensor_k\omega_{\P}(e-1))\ar[r] & 
\ext^2W^\vee\tensor_k\omega_{\P}(2(e-1))\ar[r] & 
\cdots
}$}
\\\\
\scalebox{0.95}{$
\xymatrixcolsep{10pt}
\xymatrix{ 
\cdots\ar[r]&\ext^{n+1}V^\vee\tensor_k\omega_{\P}(n+1)\ar[r]\ar[d]^{\varphi_{n+1}^*} 
&\ext^{n+2}V^\vee\tensor_k\omega_{\P}(n+2)\ar[r]\ar[d]^{\varphi_{n+2}^*} & 0\\
\cdots\ar[r]& \ext^{n+1}W^\vee\tensor_k\omega_{\P}((n+1)(e-1))\ar[r]& 
\ext^{n+2}W^\vee\tensor_k\omega_{\P}((n+2)(e-1))\ar[r]& 0. }$}
\end{multline}
 To be precise, we apply $\sH{om}(-,\omega_{\P})$ without introducing signs in the differentials, and we put the term $V^\vee\tensor_k\omega_{\P}(1)$ in degree 0. The maps $\phi_i^*$ thus give us a map of resolutions of $\omega_\P$. For later reference, we record these two complexes separately.
 \begin{equation}\label{eqn:DualComplex1}
 0\to  \omega_{\P}\to V^\vee\tensor_k\omega_{\P}(1)\to \ext^2V^\vee\tensor_k\omega_{\P}(2)\to
 \cdots\to \ext^{n+2}V^\vee\tensor_k\omega_{\P}(n+2)\to 0
 \end{equation}

 \begin{multline}\label{eqn:DualComplex2}
 0\to\omega_{\P}\to W^\vee\tensor_k\omega_{\P}(e-1))\to
\ext^2W^\vee\tensor_k\omega_{\P}(2(e-1))\to 
\cdots\\\to \ext^{n+1}W^\vee\tensor_k\omega_{\P}((n+1)(e-1))\to
\ext^{n+2}W^\vee\tensor_k\omega_{\P}((n+2)(e-1))\to 0
 \end{multline}
 
Via the isomorphism $\can$ of \eqref{eqn:Canon}, the last nonzero term on the 
right in the complex \eqref{eqn:DualComplex1} is canonically identified with $\sO_{\P}$, so this complex gives an element
 \[
\Xi\in \Hom_{D^b(\P)}(\sO_\P[-n-1], \omega_\P)= H^{n+1}(\P, \omega_\P).  
\]

We also have the trace map 
\[
\Tr_{\P/k}:H^{n+1}(\P, \omega_\P)\xrightarrow{\sim}k.
\]
As there are several different normalizations of the trace map in the literature (see for example \cite{ST}), we make precise our normalization. We have the map
\[
d\log:\sO_\P^\times\to \Omega^1_\P
\]
inducing the 1st Chern class map
\[
c_1:\Pic(\P)=H^1(\P,\sO_\P^\times)\to H^1(\P, \Omega^1_\P)
\]
and giving us the (non-zero) element $c_1(\sO_\P(1))^{n+1}\in H^{n+1}(\P, \omega_\P)$. We normalize the trace map by requiring
\[
\Tr_{\P/k}(c_1(\sO_\P(1))^{n+1})=1.
\]

\begin{lemma}\label{lem:Tr1} $\Tr_{\P/k}(\Xi)=1$.
\end{lemma}

\begin{proof} As both $\Tr_{\P/k}$ and $\Xi$ are invariant under automorphisms of $V$, we may choose a basis $x_0,\ldots, x_{n+1}$ for $V$ to make the computation. As above, this gives $S(V)=k[X_0,\ldots, X_{n+1}]$, $\P=\P^{n+1}_k$,   and the Koszul complex \eqref{canonical} is identified with the Koszul complex for the surjection
\[
\oplus_{i=0}^{n+1}\sO_{\P^{n+1}}(-1)x_i\to \sO_{\P^{n+1}}
\]
sending $\ell\cdot x_i$ to $\ell\cdot X_i$, for $\ell$ a local section of $\sO_{\P^{n+1}}(-1)$.

For $I=(i_0,\ldots, i_r)$ with the $i_j$ integers and $0\le i_0<\ldots <i_r\le n+1$, we write $x_I$ for $x_{i_0}\wedge\ldots\wedge x_{i_r}$. We let $x^0,\ldots, x^{n+1}$ be the dual basis for $V^\vee$, $x^I:=x^{i_0}\wedge\ldots\wedge x^{i_r}$. The term $\ext^{r+1}V^\vee\otimes \omega_\P(r+1)$ in the dual complex is thus
\[
\ext^{r+1}V^\vee\otimes \omega_\P(r+1)=\oplus_{|I|=r+1}\omega_\P(r+1)\cdot x^I
\]
and the differential $d^r:\oplus_{|I|=r+1}\omega_\P(r+1)\cdot x^I\to \oplus_{|J|=r+2}\omega_\P(r+2)\cdot x^J$ is the sum of terms $d^r_{J,I}:\omega_\P(r+1)\cdot x^I\to \omega_\P(r+2)\cdot x^J$. Explicitly, if 
$I=(i_0,\ldots, i_r)$ as above and $J=(i_0,\ldots, i_{\ell-1}, j, i_{\ell},\ldots, i_r)$, then for a local section $\eta$
 of $\omega_\P(r+1)$ we have
 \[
 d^r_{J,I}(\eta\cdot x^I)=(-1)^\ell X_j\eta\cdot x^J.
 \]
 If $(I,J)$ is not of this form, then $d^r_{J,I}$ is zero.

 Let $\sU:=\{U_0,\ldots, U_{n+1}\}$ be the standard affine open cover of $\P^{n+1}$, with $U_i$ defined by $X_i\neq0$, and let $\sC^*(\sU, \omega_\P)$ be the ordered \v{C}ech sheaf complex,
 \[
  \sC^r(\sU, \omega_P)=\oplus_{|I|=r+1}j_{I*}j_I^*\omega_\P
  \]
  where for $I=(i_0,\ldots, i_r)$ an ordered index, $U_I:=\cap_{j=0}^r U_{i_j}$ with inclusion $j_I:U_I\to \P$. This gives us the augmented complex
  \[
  0\to \sC^{-1}(\sU, \omega_P)\xrightarrow{\epsilon} \sC^0(\sU, \omega_P)\to\ldots\to  \sC^{n+1}
  (\sU, \omega_P)\to 0
  \]
with $\sC^0 (\sU, \omega_P)$ in degree 0, and  $\sC^{-1}(\sU, \omega_P):=\omega_\P$ with $\epsilon$ the sum of restriction maps.

Define the map
 \[
 \alpha^r:\oplus_{|I|=r+1}\omega_\P(r+1)\cdot x^I\to \sC^r(\sU, \omega_P)
 \]
 by sending $\eta\cdot x^I$, for $\eta$ a section of $\omega_\P(r+1)$ over some open $V\subset \P$, to $(1/X^I)\cdot \eta_{|V\cap U_I}$, where $X^I=\prod_{j=0}^rX_{j_i}$; we set $\alpha^{-1}=\id_{\omega_\P}$.  It follows directly that the $\alpha^r$ define a map of complexes. Note that $\alpha^{n+1}(\can(1))$ is the \v{C}ech cochain 
\[
(\frac{\Omega}{X_0\cdots X_{n+1}}\in \omega_\P(U_{0,\ldots, n+1}))\in \sC^{n+1}(\sU, \omega_\P).
 \]
 
Thus the element $\Xi$ is represented in \v{C}ech cohomology by the \v{C}ech co-chain 
$(\Omega/\prod_{j=0}^{n+1}X_j\in \omega_\P(U_{0,\ldots, n+1}))$. Note that  $\sO(1)\in H^1(\P,\sO_\P^\times)$ is represented by the co-chain $(X_j/X_i\in \sO_\P(U_i\cap U_j))_{i,j}\in \sC^1(\sU, \sO_\P^\times)$, so $c_1(\sO(1))$ is represented by the co-chain 
\[
(d\log(X_j/X_i)\in \Omega^1_\P(U_i\cap U_j))_{i,j}\in \sC^1(\sU, \Omega^1_\P),
\]
and $c_1(\sO(1))^{n+1}$  is represented by the co-chain 
\[
(d\log(X_1/X_0)\wedge d\log(X_2/X_1)\wedge\ldots\wedge d\log(X_{n+1}/X_n)\in \omega_\P(U_{01\ldots n+1})\in \sC^{n+1}(\sU, \omega_\P).
\]
An elementary computation gives
\[
\frac{\Omega}{X_0\cdots X_{n+1}}=d\log(X_1/X_0)\wedge d\log(X_2/X_1)\wedge\ldots\wedge d\log(X_{n+1}/X_n).
\]
Thus
\[
\Xi=c_1(\sO(1))^{n+1}
\]
in $H^{n+1}(\P,\omega_\P)$, hence $\Tr_{\P/k}(\Xi)=1$, by our choice of normalization.
\end{proof}

The diagram \eqref{diagram} gives two complexes, both of which yield acyclic 
resolutions for $\omega_{\P}$, and both may be used to compute 
$H^{n+1}(\P,\omega_{\P})$.  The resolution \eqref{eqn:DualComplex2} of $\omega_\P$ gives us the map
\begin{equation}\label{eqn:MapToDb}
\Hom_{\sO_\P}(\sO_\P, \ext^{n+2}W^\vee\tensor_k\omega_{\P}((n+2)(e-1)))\to
\Hom_{D^b(X)}(\sO_\P[-n-1], \omega_\P)=H^{n+1}(\P,\omega_\P),
\end{equation}
and we have the isomorphism induced by $\can$
\begin{multline*}
 \Hom_{\sO_\P}(\sO_\P, \ext^{n+2}W^\vee\tensor_k\omega_{\P}((n+2)(e-1)))\\\cong
\Hom_{\sO_\P}( \ext^{n+2}V^\vee\tensor_k\omega_{\P}(n+2), \ext^{n+2}W^\vee\tensor_k\omega_{\P}((n+2)(e-1)))\\
\cong \Hom_{\sO_\P}( \sO_\P, \ext^{n+2}W^\vee\tensor_k \ext^{n+2}V\tensor_k\sO_\P((n+2)(e-2)))
\end{multline*}
Moreover, the map $\varphi_{n+2}^*$ induces the element $\Xi\in H^{n+1}(\P, \omega_\P)$.

We note that 
\[{\rm coker}\left(H^0(\P,\ext^{n+1}W^\vee\tensor_k\omega_{\P}((n+1)(e-1))\to 
H^0(\P,\ext^{n+2}W^\vee\tensor_k\omega_{\P}((n+2)(e-1))\right)\]
is identified via the above isomorphism with
\[\ext^{n+2}W^\vee\tensor_k\ext^{n+2}V\tensor_kR(W)_{(n+2)(e-2)},\]
and  the natural map \eqref{eqn:MapToDb} descends to an isomorphism 
\[
\rho:\ext^{n+2}W^\vee\tensor_k\ext^{n+2}V\tensor_kR(W)_{(n+2)(e-2)}\xrightarrow{\sim} H^{n+1}(\P,\omega_\P).
\]
Via this isomorphism the element 
\[
\bar{\vartheta}_{n+2}^*(\can(1))\in \ext^{n+2}W^\vee\tensor_k\ext^{n+2}V\tensor_k R(W)_{(n+2)(e-2)}
\]
 induced by 
$\varphi_{n+2}^*$ maps to the trace one element $\Xi\in H^{n+1}(\P,\omega_\P)$.

We thus obtain a commutative diagram of isomorphisms
\[
\xymatrix{
 H^0(\P,\sO_{\P})\ar[r]^-{\times\Xi}\ar[d]^{\bar{\vartheta}_{n+2}^*}& 
H^{n+1}(\P,\omega_{\P})\\
\ext^{n+2}W^\vee\tensor_k\ext^{n+2}V\tensor_kR(W)_{(n+2)(e-2)}\ar[ur] _-\rho
}
\]
Having chosen bases for $V$ and $W$, this gives the  isomorphism 
\[
\ext^{n+2}W^\vee\tensor_k\ext^{n+2}V\tensor_kR(W)_{(n+2)(e-2)}\cong R(W)_{(n+2)(e-2)},
\]
and our matrix $(a_{ij})\in M_{n-2, n-2}(k[X_0,\ldots, X_{n+1}]_{(n+2)(\dg-2)}$, with
we have 
\[
\bar{\varphi}_{n+2}^*(\can(1))=\det(a_{ij})\in R(W)_{(n+2)(e-2)}.
\]

In our case at hand, with $W=W(X,V)$, we have the isomorphism $\ev_F:V^\vee\to W$, transforming the above to the commutative diagram of isomorphisms
\[
\xymatrix{
 H^0(\P,\sO_{\P})\ar[r]^-{\times\Xi}\ar[d]_{\ext^{n+2}\ev_F^\vee\tensor\id\circ\bar{\vartheta}_{n+2}^*}& 
H^{n+1}(\P,\omega_{\P})\\
(\ext^{n+2}V)^{\otimes 2}\tensor_kJ(X,V)_{(n+2)(e-2)}\ar[ur] _-{\quad\rho\circ(\ext^{n+2}\ev_F^\vee\tensor\id)^{-1}}&
}
\]
Choosing a basis $x_0,\ldots, x_{n+1}$ for $V$ as above, we have $J(X,V)=J(F)$ and $1\in H^0(\P,\sO_{\P})$ maps to the   Scheja-Storch element for $F$.

\subsection{Trace on the hypersurface}
Let 
\[
\sL=\ext^{n+2}W^\vee\tensor_k\ext^{n+2}V\tensor_k\sO_{\P}((n+2)(e-2))).
\]
 We define an isomorphism of  the  complex \eqref{eqn:DualComplex2} with the complex \eqref{koszulW} twisted by $\sL$ as follows.
 The interior multiplication gives the isomorphism
 \[
 \iota_r:\ext^rW\otimes\ext^{n+2}W^\vee\to \ext^{n+2-r}W^\vee.
 \]
  If we give $W$ a basis $w_0,\ldots, w_{n+1}$ with dual basis $w^0,\ldots, w^{n+1}$ for $W^\vee$, then $\iota_r$ sends $w_I\otimes w^0\wedge\ldots\wedge w^{n+1}$ to $\sgn(I, I^c)w^{I^c}$. Here $I^c$ is the complement of $I$ in $\{0,\ldots, n+1\}$, with everything written in increasing order, and  $\sgn(I, I^c)\in\{\pm1\}$ is defined by
\[
w_I\wedge w_{I^c}= \sgn(I, I^c)w_0\wedge\ldots\wedge w_{n+1}.
\]
 If $I=(i_0,\ldots, i_{r-1})$, then 
  \[
 \sgn(I, I^c)=(-1)^{\sum_{j=0}^{r-1}i_j-j}
 \]
We let $\alpha_r:\ext^{n+2-r}W^\vee\to \ext^rW\otimes\ext^{n+2}W^\vee$ be the inverse of $\iota_r$. 

The canonical isomorphism $\omega_\P=\ext^{n+2}V\otimes\sO_\P(-n-2)$ gives us the canonical isomorphism
\[
\omega_\P((n+2)(e-1))=\ext^{n+2}V\otimes_k\sO_\P((n+2)(e-2))
\]
 
Define the sign $\epsilon_{n,r}\in \{\pm1\}$ by
\[
\epsilon_{n,r}=\begin{cases} (-1)^{(r-1)(r-2)/2}&\text{ for $n$ even,}\\
(-1)^{r(r-1)/2}&\text{ for $n$ odd.}
\end{cases}
\]

\begin{lemma}  The maps $\epsilon_{n, r}\alpha_{n+2-r}\otimes\id$, 
\begin{multline*}
\ext^{r}W^\vee\otimes\omega_\P(r(e-1))\\\xrightarrow{\epsilon_{n, r}\alpha_{n+2-r}\otimes\id}
\ext^{n+2-r}W\otimes \ext^{n+2}W^\vee\otimes \omega_\P(r(e-1))\\
\cong \ext^{n+2-r}W\otimes\sO_\P((r-n-2)(e-1))\otimes \sL
\end{multline*}
define an isomorphism of complexes
\begin{equation}\label{eqn:Alpha}
\alpha:\eqref{eqn:DualComplex2}\to\eqref{koszulW}\otimes\sL.
\end{equation}
\end{lemma}

\begin{proof} We claim that the diagram
\[
\xymatrix{
\ext^{r}W^\vee\otimes\omega_\P(r(e-1))\ar[d]^{\alpha_{n+2-r}\otimes\id}\ar[r]^d&
\ext^{r+1}W^\vee\otimes\omega_\P((r+1)(e-1))\ar[d]^{\alpha_{n+1-r}\otimes\id}\\
\ext^{n+1-r}W\otimes\sO_\P((r-n-1)(e-1))\otimes \sL\ar[r]^d&
\ext^{n-r}W\otimes\sO_\P((r-n)(e-1))\otimes \sL
}
\]
commutes up to the sign $(-1)^{n+1+r}$. As
\[
(-1)^{\frac{1}{2}(r(r-1)-(r-1)(r-2))}=(-1)^{r+1}
\]
and
\[
(-1)^{\frac{1}{2}((r+1)r-r(r-1))}=(-1)^r
\]
the result follows from this.

To make the computation, we use the basis $w_i:=\del F/\del X_i$ for $W$ and dual basis $w^0,\ldots, w^{n+1}$ for $W^\vee$. For an index $I=(i_0,\ldots, i_{r-1})$, $0\le i_0<\ldots<i_{r-1}\le n+1$, we let   $w^I$ denote the corresponding basis element for $\ext^rW^\vee$ and define the basis element $w_{I^c}$ for $\ext^{n+2-r}W$ similarly. In our computation, we omit the contribution of $\ext^{n+2}W^\vee$ and of local sections of $\omega_\P(r(e-1))$, $\sL$ etc., from the notation

For an element $w^I\in \ext^{r}W^\vee$, we have
\[
\alpha_{n+2-r}(w^I)=\sgn(I^c, I)w_{I^c}\in \ext^{n+2-r}W
\]
Take $j\in I^c$, with $i_{\ell-1}<j<i_\ell$ and let $J=(i_0,\ldots, i_{\ell-1}, j,i_\ell,\ldots, i_{r-1})$. Then the $w^J$ component of $d(w^I)$ is
\[
d(w^I)^J:=(-1)^\ell \frac{\del F}{\del X_j}\cdot w^J
\]
Similarly, $J^c=I^c\setminus \{j\}$, we have
\[
\alpha_{n+1-r}(w^J)=\sgn(J^c, J)w_{J^c}\in \ext^{n+1-r}W
\]
and  the $J^c$-component of $d(w_{I^c})$ is 
\[
d(w_{I^c})_{J^c}=(-1)^{j -\ell}\frac{\del F}{\del X_j}\cdot w_{J^c}
\]

Since
\begin{align*}
\sgn(I^c, I)&=(-1)^{r(n+2-r)}\sgn(I, I^c)\\
&=(-1)^{r(n+2-r)}(-1)^{\sum_{m=0}^{r-1}i_m-m}\\
\sgn(J^c, J)&=(-1)^{(r+1)(n-r+1)}\sgn(J, J^c)\\
&=(-1)^{(r+1)(n-r+1)}(-1)^{\sum_{m=0}^{\ell-1}i_m-m}\cdot(-1)^{j-\ell}\cdot (-1)^{\sum_{m=\ell}^{r-1}i_m-m-1},
\end{align*}
we have
\[
\sgn(I^c, I)\cdot \sgn(J^c, J)\cdot (-1)^\ell\cdot (-1)^{j-\ell}=(-1)^{n+1+r},
\]
which proves our claim.
\end{proof} 

It will be convenient to replace $\eqref{koszulW}\otimes\sL$ with another isomorphic complex. We have the isomorphism
\begin{multline*}
\sigma^r: \sO_\P(-r(e-1))\otimes\sL=
 \sO_\P(-r(e-1))\otimes\ext^{n+2}W^\vee\otimes \ext^{n+2}V\otimes \sO_\P((n+2)(e-2))\\\xrightarrow{\sim}
 \sO_\P(r)\otimes \sO_\P((n-r)X)\otimes  \omega_\P(X)^{\otimes 2}.
\end{multline*}
defined by using $F$ to identify $\sO_\P(je)$ with $\sO_P(jX)$, using the dual of the isomorphism $\ev_F:V^\vee\xrightarrow{\sim}W$ to identify $\ext^{n+2}W^\vee$ with $\ext^{n+2}V$,   and using the canonical isomorphism $\can:\ext^{n+2}V\otimes\sO_\P(-(n+2))\xrightarrow{\sim}\omega_\P$.  This gives us the isomorphisms
\[
\id_{\ext^rW}\otimes\sigma^r:\ext^rW\otimes \sO_\P(-r(e-1))\otimes\sL\xrightarrow{\sim}
\ext^rW\otimes\sO_\P(r)\otimes \sO_\P((n-r)X)\otimes  \omega_\P(X)^{\otimes 2}
\]

Using $F$ again, we have the isomorphism $\sigma_0:W\otimes \sO_\P(1)\otimes \sO_\P(-X)\to
W\otimes \sO_\P(1-e)$, so the surjection $W\otimes \sO_\P(1-e)\to \sO_\P$ used to define  \eqref{koszulW} defines a surjection $W\otimes \sO_\P(1)\otimes \sO_\P(-X)\to\sO_\P$, giving us another Koszul complex
\begin{equation}\label{koszulW2}
0\to \ext^{n+2}W\otimes\sO_\P(n+2)\otimes\sO_\P(-(n+2)X)\to\ldots\to W\otimes\sO_\P(1)\otimes\sO_\P(-X)\to \sO_\P\to 0
\end{equation}
Letting $\sL':= \sO_\P(nX)\otimes  \omega_\P(X)^{\otimes 2}$, the maps $\id_{\ext^rW}\otimes\sigma^r$ define an isomorphism of Koszul complexes
\begin{equation}\label{eqn:Sigma}
\sigma:\eqref{koszulW}\otimes\sL\to \eqref{koszulW2}\otimes\sL'
\end{equation}

We complete the picture by defining a map of complexes 
\[
\beta: \text{\eqref{koszulW2}}\otimes\sL'\to i_*[0\to \sT^n_X\otimes\omega_X^{\otimes 2}\xrightarrow{\epsilon_n\otimes\id_{\omega_X^{\otimes 2}}} \sD(n)[-n-1]\otimes \omega_X^{\otimes 2}].
\]

For this, we have the surjection $\Eul:V^\vee\otimes\sO_\P(1)\to \sT_\P$ in the Euler sequence. 
For each $r=0,\ldots, n+1$, we  have the map
\begin{multline*}
(\ext^r\Eul\otimes\id)\circ(\ext^r \ev_F^{-1}\otimes\id):\ext^rW\otimes \sO_\P(r)\otimes\sO_\P((n-r)X)\otimes\omega_\P(X)^{\otimes 2}\\\to
\sT^r_\P\otimes\sO_\P((n-r)X)\otimes\omega_\P(X)^{\otimes 2}
\end{multline*}
and then the map
\[
i^*\otimes(\res_X)^{\otimes 2}:
\sT^r_\P\otimes\sO_\P((n-r)X)\otimes\omega_\P(X)^{\otimes 2}
\to i^*\sT^r_\P\otimes\sO_X((n-r)X)\otimes\omega_X^{\otimes 2}.
\]
For $r=0,\ldots, n$, let
\begin{multline*}
\beta_r:\ext^rW\otimes \sO_\P(r)\otimes\sO_\P((n-r)X)\otimes\omega_\P(X)^{\otimes 2}\\\to
i_*[i^*\sT^r_\P\otimes\sO_X((n-r)X)\otimes\omega_X^{\otimes 2}]=i_*(\sD(n)[-n]\otimes\omega_X^{\otimes 2})^r
\end{multline*}
be the composition of these two.  For $r=n+1$, take the composition as above and compose with the isomorphism $\res_X^*:i^*\sT^{n+1}_\P(-X)\to \sT^n_X$, this being the inverse of the dual of the residue isomorphism $\res_X:i^*\omega_\P(X)\to \omega_X$. This gives the map
\[
\beta_{n+1}:\ext^{n+1}W\otimes \sO_\P(n+1)\otimes\sO_\P(-X)\otimes\omega_\P(X)^{\otimes 2}\to
\sT^n_X\otimes\omega_X^{\otimes 2}.
\]

\begin{lemma} The maps 
\begin{multline*}
(-1)^{(n+1)(r-n)}\beta_r:\ext^rW\otimes \sO_\P(r)\otimes\sO_\P((n-r)X)\otimes\omega_\P(X)^{\otimes 2}\\\to
i_*[i^*\sT^r_\P\otimes\sO_X((n-r)X)\otimes\omega_X^{\otimes 2}]=i_*(\sD(n)[-n-1]\otimes\omega_X^{\otimes 2})^{n+1-r}
\end{multline*}
$r=0,\ldots, n$, together with 
\[
\beta_{n+1}:\ext^{n+1}W\otimes \sO_\P(n+1)\otimes\sO_\P(-X)\otimes\omega_\P(X)^{\otimes 2}\to
\sT^n_X\otimes\omega_X^{\otimes 2},
\]
  define a map of complexes
\begin{equation}\label{eqn:Beta}
\beta: \eqref{koszulW2}\otimes\sL'\to i_*[0\to \sT^n_X\otimes\omega_X^{\otimes 2}\xrightarrow{\epsilon_X\otimes\id} \sD(n)[-n-1]\otimes \omega_X^{\otimes 2}] 
\end{equation}
\end{lemma}

\begin{proof} For a complex $A^*$ and integer $m$, let $A\{m\}^*$ denote the complex with $A\{m\}^j=A^{m+j}$ and with $d_{A\{m\}}^j=d_A^{m+j}$ (no signs are introduced in the differential). 

The maps $\beta_r$ define a map of Koszul complexes
\[
\eqref{koszulW2}\otimes\sL'\{-n-1\}\to 
i_*[0\to \sT^n_X\otimes\omega_X^{\otimes 2}\{-n-1\}\xrightarrow{\epsilon_X\{-n-1\}\otimes\id} \sD(n)\otimes \omega_X^{\otimes 2}] 
\]
so give a map of shifted Koszul complexes
\[
\eqref{koszulW2}\otimes\sL'\to 
i_*[0\to \sT^n_X\otimes\omega_X^{\otimes 2}\xrightarrow{\epsilon_X\otimes\id} \sD(n)\{-n-1\}\otimes \omega_X^{\otimes 2}] 
\]
The  factor $(-1)^{(n+1)(r-n)}$ for $\beta_r$ takes care of the of the factor $(-1)^{n+1}$ in the differential for the complex $\sD(n)[-n-1]$, as compared with $\sD(n)\{-n-1\}$, away from the  augmentation $r=n+1$. In the square involving the augmentation, there is no sign introduced in the target complex, and the sign for $\beta_n$ is $+1$,  so we need to use the map $\beta_{n+1}$ without sign.
\end{proof}

We truncate the complex \eqref{eqn:DualComplex2} in degree $\le 1$, giving the short exact sequence
\[
0\to \omega_\P\to W^\vee\otimes\omega_\P(e-1)\xrightarrow{\bar{d}^1} \sF\to 0
\]
Composing the maps  \eqref{eqn:Alpha},  \eqref{eqn:Sigma}, and  \eqref{eqn:Beta} gives us the map of complexes
\[
\beta\circ\sigma\circ\alpha:\eqref{eqn:DualComplex2}\to
 i_*[0\to \sT^n_X\otimes\omega_X^{\otimes 2}\xrightarrow{\epsilon_n\otimes\id_{\omega_X^{\otimes 2}}} \sD(n)[-n]\otimes \omega_X^{\otimes 2}]\{-1\}.
\]
The isomorphism   $\iota(n)_X\otimes\id:\sT^n_X\otimes\omega_X\otimes\omega_X\to \omega_X$ gives the map 
\[
(\iota(n)_X\otimes\id)\circ \overline{\beta^0\circ\alpha^0}:\sF\to i_*\omega_X
\]
Moreover, we have the map $\ev_W:W^\vee\otimes\omega_\P(e-1)\to \omega_\P(X)$ defined as the composition
\[
W^\vee\otimes\omega_\P(e-1)\xrightarrow{\ev_F^\vee\otimes\id} V\otimes\sO_\P(-1)\otimes\omega_\P(X)\xrightarrow{\ev\otimes\id}\omega_\P(X).
\]
\begin{lemma}\label{lem:OmegaComparison}
The diagram
\begin{equation}\label{eqn:LastDiagram}
\xymatrix{
0\ar[r]&\omega_\P\ar[r]\ar[d]^{(-1)^{n+1}\cdot \dg}&W^\vee\otimes\omega_\P(e-1)\ar[d]^{(-1)^{n+1}\ev_W}\ar[r]^-{\bar{d}^1}&\sF\ar[r]\ar[d]^{(\iota(n)_X\otimes\id)\circ\overline{\beta^0\circ\alpha^0}}&0\\
0\ar[r]&\omega_\P\ar[r]&\omega_\P(X)\ar[r]^\res&i_*\omega_X\ar[r]&0
}
\end{equation}
commutes.
\end{lemma}

\begin{proof} All the maps are canonically defined, so we may fix a basis $x_0,\ldots, x_{n+1}$ for $V$ to prove the result. As above, this gives the dual basis $x^0,\ldots, x^{n+1}$ for $V^\vee$, and the isomorphism $\ev_F:V^\vee\to W$ gives the corresponding bases $w_0,\ldots, w_{n+1}$ for $W$, with  $w_i=\del F/\del X_i$, and dual basis $w^0,\ldots, w^{n+1}$ for $W^\vee$. We let $x_\bullet:=x_0\wedge\ldots\wedge x_{n+1}$ and $w^\bullet:=w^0\wedge\ldots\wedge w^{n+1}$ denote the corresponding generators for $\ext^{n+2}V$ and $\ext^{n+2}W^\vee$. For $i\in \{0,\ldots, n+1\}$ let 
\[
w_{\hat{i}}=w_0\wedge\ldots\wedge w_{i-1}\wedge w_{i+1}\wedge\ldots\wedge w_{n+1}.
\]
and define $x^{\hat{i}}$,   $\del/\del X_{\hat{i}}$, etc., similarly.

We factor the map $(\iota_X(n)\otimes\id)\circ\beta^1\circ(\id\otimes\sigma)^1\circ\alpha^1:W^\vee\otimes\sO_\P(e-1)\to i_*\omega_X$ as the composition
\begin{multline*}
W^\vee\otimes\sO_\P(e-1)\xrightarrow{(\id\otimes\sigma^{n+1})\circ(\alpha_{n+1}\otimes\id)}
\ext^{n+1}W\otimes\sO_\P(n+1)\otimes\sO_\P(-X)\otimes \omega_\P(X)^{\otimes2}\\
\xrightarrow{\ext^{n+1}\Eul\otimes\id}
\sT_\P^{n+1}\otimes\sO_\P(-X)\otimes \omega_\P(X)^{\otimes2}\cong
\sT_\P^{n+1}\otimes\omega_\P\otimes \omega_\P(X)\\
\xrightarrow{\iota_\P(n+1)\otimes\id}\omega_\P(X)
\xrightarrow{\res_X}i_*\omega_X
\end{multline*}
Let 
\[
\gamma:W^\vee\otimes\sO_\P(e-1)\to \omega_\P(X)
\]
be the composition described above, so we have
\[
(\iota_X(n)\otimes\id)\circ\beta^1\circ(\id\otimes\sigma)^1\circ\alpha^1=
\res_X\circ \gamma.
\]

We compute $\gamma$ in local coordinates near $p\in \P$. For $l$ a section of $\sO_\P(1)$, non-vanishing at   $p$, we have the local generator $\Omega\cdot l^{-(n+2)}$ for $\omega_\P$ near $p$. 

The isomorphism $\alpha_{n+1}:W^\vee\to \ext^{n+1}W\otimes\ext^{n+2}W^\vee$ sends $w^i$ to $(-1)^{n+1-i}
w_{\hat{i}}\otimes w^\bullet$ and we have
\[
(\id\otimes\sigma)^1\circ\alpha^1(w^i\otimes \Omega\cdot l^{e-(n+3)})=
(-1)^{n+1-i}\cdot w_{\hat{i}}\cdot l^{n+1}\otimes\frac{F}{l^e}\otimes \frac{\Omega}{F}\cdot l^{e-(n+2)}\otimes  \frac{\Omega}{F}\cdot l^{e-(n+2)}.
\]
$w_{\hat{i}}$ maps to $x^{\hat{i}}\in \ext^{n+1}V$ via $\ext^{n+1}\ev_F^{-1}$, and $l^{n+1}x^{\hat{i}}$ maps to $l^{n+1}\del/\del X_{\hat{i}}$ in $\sT_\P^{n+1}$ via $\ext^{n+1}\Eul$. This gives us
\[
(-1)^{n+1-i}\cdot l^{n+1}\del/\del X_{\hat{i}}\otimes\frac{F}{l^e}\otimes \frac{\Omega}{F}\cdot l^{e-(n+2)}\otimes  \frac{\Omega}{F}\cdot l^{e-(n+2)}
\]
in $\sT_\P^{n+1}\otimes\sO_\P(-X)\otimes \omega_\P(X)^{\otimes2}$, which maps to 
\[
(-1)^{n+1-i}\cdot l^{n+1}\del/\del X_{\hat{i}}\otimes \Omega\cdot l^{-(n+2)}\otimes  \frac{\Omega}{F}\cdot l^{e-(n+2)}
\]
in $\sT_\P^{n+1}\otimes\omega_\P\otimes \omega_\P(X)$. Since the interior multiplication sends
$\del/\del X_{\hat{i}}\otimes \Omega$ to $(-1)^iX_i$, this element maps to
\[
(-1)^{n+1}X_i\cdot l^{-1}\otimes  \frac{\Omega}{F}\cdot l^{e-(n+2)}
\]
in $\omega_\P(X)$, and finally to 
\[
(-1)^{n+1}X_i\cdot l^{-1}\otimes \res_X(\frac{\Omega}{F}\cdot l^{e-(n+2)})
\]
in $i_*\omega_X$. 

On the other hand, the map $\ev_W$ sends $w^i\otimes \Omega\cdot l^{e-(n+3)}$ first to $x_i\cdot l^{-1}
\otimes \Omega/F\cdot l^{e-(n+2)}$ and then to $X_i\cdot l^{-1}\otimes (\Omega/F)\cdot l^{e-(n+2)}$. Thus
\[
\ev_W(w^i\otimes \Omega\cdot l^{e-(n+3)})=(-1)^{n+1} \gamma(w^i\otimes \Omega\cdot l^{e-(n+3)})
\]
so $\ev_W=(-1)^{n+1}\gamma$ and hence
\[
\res_X\circ\ev_W=(-1)^{n+1}\ev_X\otimes\id\circ \beta^1\circ \alpha^1
=(-1)^{n+1}\ev_X\otimes\id\circ \overline{\beta^0\circ \alpha^0}\circ\bar{d}^1
\]
This shows that the right-hand square in the diagram \eqref{eqn:LastDiagram} commutes.

For the left-hand square, take our local generating section $\Omega\cdot l^{-(n+2)}$ of $\omega_\P$. This maps to $\sum_iw^i\cdot \del F/\del X_i\otimes \Omega\cdot l^{-(n+2)}$ in $W^\vee\otimes\omega_\P(e-1)$. Using the Euler relation and our computations above, this maps to
\[
(-1)^{n+1}\sum_iX_i  \del F/\del X_i\otimes \frac{\Omega}{F}\cdot l^{-(n+2)}
=(-1)^{n+1}\dg\cdot  \Omega\cdot l^{-(n+2)},
\]
in $\omega_\P(X)$, so the left-hand square commutes as well.
\end{proof}

Via the isomorphism $\iota(n)_X\otimes\id:\sT^n_X\otimes\omega_X^{\otimes 2}\to \omega_X$ and the resolution $\sT^n_X\otimes\omega_X^{\otimes 2}\to \sD(n)[-n]\otimes\omega_X^{\otimes 2}$, 
the map 
\[
\xi:\sO_\P\to i_*(\sD(n)[-n]\otimes\omega_X^{\otimes 2})^n
\]
defined as the composition
\begin{multline*}
\sO_\P\cong \ext^{n+2}V^\vee\otimes\omega_\P(n+2)\xrightarrow{\phi_{n+2}^\vee}
\ext^{n+2}W^\vee\otimes\omega_\P((n+2)(e-1))\\
\xrightarrow{(-1)^{n(n+1)/2}\alpha_0\otimes\id}\sL
\xrightarrow{\sigma_0}\sO_\P(nX)\otimes\omega_\P(X)^{\otimes 2}\\
\xrightarrow{(-1)^{n^2}\beta_0}i_*(\sD(n)[-n]\otimes\omega_X^{\otimes 2})^n
\end{multline*}
defines an element $\tau\in H^n(X, \omega_X)$. 

\begin{proposition} $\Tr_{X/k}(\tau)= (-1)^{n+1}\cdot \dg$.
\end{proposition}

\begin{proof} We may extend the resolution 
\[
0\to \omega_X[-1]\xrightarrow{\epsilon_n[-1]\circ \iota_X(n)^{-1})\otimes\id} \sD(n)[-n-1]\otimes\omega_X^{\otimes 2}
\]
to a resolution of $\omega_\P$ by patching in the exact residue sequence. This gives the resolution of $\omega_\P$
\[
0\to \omega_\P\to \omega_\P(X)\to i_*(\sD(n)[-n-1]\otimes\omega_X^{\otimes 2}).
\]
It follows from Lemma~\ref{lem:Tr1} and  Lemma~\ref{lem:OmegaComparison} that the element $\tau'\in  H^{n+1}(\P, \omega_\P)$ corresponding to the map 
$\xi:\sO_\P\to i_*(\sD(n)[-n-1]\otimes\omega_X^{\otimes 2})^{n+1}$ defined above   has
\[
\Tr_{\P/k}(\tau')=(-1)^{n+1}\cdot \dg.
\]

Let $\delta:H^n(X, \omega_X)\to H^{n+1}(\P,\omega_\P)$ be the boundary map in the residue sequence. It  follows from Remark~\ref{rem:DerivedFunctorDeltaFunctor} and Lemma~\ref{lem:OmegaComparison} that $\tau'=\delta(\tau)$. By Lemma~\ref{lem:ResiduePushforward}, we have
\[
i_*(\tau)=\tau'
\]
Thus
\[
\Tr_{X/k}(\tau)=\Tr_{\P/k}(i_*(\tau))=\Tr_{\P/k}(\tau')=(-1)^{n+1}\cdot \dg.
\]
\end{proof} 

With this in hand, we can complete our second proof of Theorem~\ref{thm:Main}. 

\subsection{The proof of Theorem~\ref{thm:Main}} We fix a basis $x_0,\ldots, x_{n+1}$ for $V$, giving the canonical isomorphisms $\ext^{n+2}V\cong k$, and
\[
J(X,W)\otimes\ext^{n+2}V\cong J(X,W)\otimes(\ext^{n+2}V)^{\otimes 2}\cong J(F). 
\]
Take  $A\in J(F)_{(q+1)\dg-n-2}$, $B\in J(F)_{(p+1)\dg-n-2}$ with $p+q=n$. Letting $\omega_A=\psi_q(A)\in H^q(X, \Omega_X^p)$, $\omega_B=\psi_p(B)\in H^p(X, \Omega_X^q)$, we have by Lemma~\ref{lem:Products},
\[
\Tr_{X/k}(\omega_A\cup\omega_B)=(-1)^{n(n-1)/2}\cdot \Tr_{X/k}(H^n(\iota_X(n)\otimes\id)(\phi'_n(AB))).
\]

On the other hand, for $n$ even, the map $\xi: \sO_\P\to i_*(\sD(n)[-n]\otimes\omega_X^{\otimes 2})^n$ defined above sends $1$ to $(-1)^{n(n+1)/2}\cdot e_F$, with $e_F\in J(F)_{(\dg-2)(n+2)}$ the Scheja-Storch element. For $n$ odd, we have $(-1)^{(n+1)(n+2)/2}\cdot e_F$, but as this differs from $(-1)^{n(n+1)/2}\cdot e_F$ by the factor $(-1)^{n+1}$, we may write this as  $(-1)^{n(n+1)/2}\cdot e_F$ for odd $n$ as well.

Thus, 
\begin{multline*}
\Tr_{X/k}(H^n(\iota_X(n)\otimes\id)(\phi'_n(e_F)))=(-1)^{n(n+1)/2}\cdot\Tr_{X/k}(\tau)\\
=(-1)^{n(n+1)/2}\cdot (-1)^{n+1}\dg=(-1)^{n(n-1)/2}\cdot (-\dg).
\end{multline*}
Therefore, if $AB =\lambda\cdot e_F$, then
\[
\Tr_{X/k}(\omega_A\cup\omega_B)=
(-1)^{n(n-1)/2}\cdot \Tr_{X/k}(H^n(\iota_X(n)\otimes\id)(\phi'_n(\lambda\cdot e_F)))=
(-\dg)\cdot \lambda
\]
which gives the desired formula.

\end{document}